\documentclass[10pt,a4paper,reqno]{amsart}
\usepackage{amsfonts,amsmath}
\usepackage[pdf]{pstricks}
\usepackage[pdftex]{graphicx}
\usepackage[utf8]{inputenc}
\usepackage[english]{babel}
\usepackage{amsthm}
\usepackage{amsmath}
\usepackage{caption}
\usepackage{subcaption}
\usepackage{booktabs}
\usepackage{fancyhdr}
\usepackage{mathrsfs}
\usepackage{amssymb}
\usepackage{graphicx}
\usepackage{wasysym}
\usepackage{enumerate}
\usepackage{mathtools}
\graphicspath{ {./images/} }
\usepackage{wrapfig}
\usepackage{hyperref}
\usepackage{blindtext}
\usepackage{xcolor}
\usepackage{tcolorbox}
\usepackage{stmaryrd}
\usepackage[a4paper,left=2.5cm,right=2.5cm]{geometry}
\usepackage{tikz}
\usetikzlibrary{patterns}
\usepackage{pgfplots}

\usepackage[numbers,sort,compress]{natbib}
\usepackage[nameinlink]{cleveref}

\definecolor{my-blue}{cmyk}{0.16, 0.1, 0, 0, 0}
\definecolor{my-yellow}{cmyk}{0,0,0.14,0}
\definecolor{my-purple}{cmyk}{0,0.08,0.08,0}
\definecolor{my-green}{cmyk}{0.18,0.04,0.11,0}
\definecolor{my-green2}{cmyk}{0.25,0.02,0.19,0}
\tcbset{colframe=white}
\DeclareMathOperator*{\esssup}{ess\,sup}

\newcommand{\norm}[1]{\|#1\|}
\newtheorem{theo}{Theorem}[section]
\newtheorem{definition}[theo]{Definition}
\newtheorem{lemma}[theo]{Lemma}

\newtheorem{prop}[theo]{Proposition}
\newtheorem{remark}[theo]{Remark}
\makeatletter
\newcommand{\newitem}[1]{%
\item[#1]\protected@edef\@currentlabel{#1}%
}
\makeatother

\definecolor{sorange}{RGB}{203,75,22}
\definecolor{sblue}{RGB}{38,139,210}
\definecolor{smagenta}{RGB}{211,54,130}

\hypersetup{
  colorlinks = true,
  linkcolor  = sorange,
  citecolor  = sblue,
  urlcolor   = sblue,
}

\numberwithin{equation}{section}
\crefname{theo}{Thm.}{Thm.}
\crefname{definition}{Def.}{Def.}
\crefname{lemma}{Lem.}{Lem.}
\crefname{remark}{Rem.}{Rem.}
\crefname{equation}{}{}
\crefname{section}{Sec.}{Sec.}
\crefname{figure}{Fig.}{Fig.}
\crefname{table}{Tab.}{Tab.}

\usepackage[pagewise]{lineno}
\definecolor{sgrey}{RGB}{153,153,153}

\usepackage[utf8]{inputenc}
\usepackage{fancyhdr}
\pgfplotsset{compat=1.16}

\title[ESFEM parabolic transmission]{Evolving finite elements for advection diffusion with an  evolving interface}
\author{C.M. Elliott \and  T. Ranner \and P. Stepanov }
\address{
	C. M. Elliott  \hfill\break 
	Mathematics Institute, University of Warwick \\
	Zeeman Building, Coventry CV4 7AL, UK}
	\email{C.M.Elliott@warwick.ac.uk}

\address{T. Ranner \hfill\break
	School of Computing, University of Leeds,
	  Leeds, LS2 9JT, UK}
	\email{T.Ranner@leeds.ac.uk}

\address{P. Stepanov \hfill\break 
	Mathematics Institute, University of Warwick \\
	Zeeman Building, Coventry CV4 7AL, UK } \email{P.Stepanov@warwick.ac.uk}

\begin{document}
\maketitle
\begin{abstract}
 The aim of this paper is to develop a numerical scheme to approximate evolving interface problems for parabolic equations based on the abstract evolving finite element framework proposed in \cite{EllRan21}. An appropriate weak formulation of the problem is derived for the use of   evolving finite elements designed to accommodate a moving  interface. Optimal order error bounds are proved for arbitrary order evolving  isoparametric finite elements. The paper concludes with numerical results for a model problem verifying orders of convergence. 
\end{abstract}

\section{Introduction}
The model studied in the paper is the following: Let $\Omega$ be a stationary domain with a moving interface $\Gamma(t)$ that encloses a subdomain $\Omega_1(t)$ and let $\Omega_2(t) = \Omega \setminus \overline{\Omega}_1(t)$.
We denote by $\nu_\Gamma$ the outward pointing normal to $\Omega_1(t)$.
Let, for $i=1$ and $2$,  $\mathcal{A}_i$ be a diffusion tensor field, $\mathcal{B}_i$ be  a vector field and $\mathcal{C}_i$ be  a  scalar field, each  continuous on  $\Omega_i(t)$. 
 Let $f_1, f_2$ and $g$ be time-dependent functions on $\Omega_1(t)$, $\Omega_2(t)$ and $\Gamma(t)$ respectively.
More precise definitions of the problem data are given in \cref{theorem:existence}.
We are interested in  well posedness and a suitable finite element scheme for the solutions of the following problem:
Find scalar fields $u_1$ on the subdomain $\Omega_1(t)$ and $u_2$ on the subdomain $\Omega_2(t)$,  which satisfy:
\begin{subequations}
\label{eq:strong}
\begin{align}
    \partial_t u_i - \nabla\cdot (\mathcal{A}_i(t;x) \nabla u_i) +\mathcal{B}_i(t;x) \cdot \nabla u_i + \mathcal{C}_i(t;x)u_i &= f_i(t;x)&&\quad \text{in} \; \Omega_i(t),\\
    u_2 &= 0&&\quad \text{on} \; \partial\Omega,\\
    u_1 -u_2 &= 0 &&\quad \text{on} \; \Gamma(t), \\
    \mathcal{A}_1(t;x)\frac{\partial}{\partial \nu_{\Gamma}} u_1 \bigg|_{\Gamma(t)}- \mathcal{A}_2(t;x)\frac{\partial}{\partial\nu_{\Gamma}}u_2\bigg|_{\Gamma(t)} &= g(t;x)&&\quad \text{on} \; \Gamma(t),\\
    u_i(0) &= u_i^0 && \quad \text{on} \; \Omega_i(0).
\end{align}
\end{subequations}
Such equations can arise as subproblems when modelling the transport and diffusion of  the concentration of a  dissolved chemical species in evolving spatial domains.
In particular, we mention applications in fluid dynamics \cite{schramm20032,Barrett_2014c,Abels2017}, materials science \cite{Barrett2013,Gurtin_1988} and cell biology \cite{Werner_2022,Ryder_2019,Hakkinen_2019}.

%Introduce the domain and the PDE
%Outline other papers and schemes (ALE framework, or maybe talk about Edelman's approach, Trace-FEM, phase-fields)
There are two main difficulties concerning this problem. The first is the evolution of the subdomains and the second is the presence of a flux jump across the interface. One common approach to moving domains is the \textit{ALE} (Arbitrary Eulerian Lagrangian) method, see \cite{DonHuePon04,ale1,ale3}. This involves having a parametrisation of the evolving region. The flow associated with this parametrisation could be physical or could be made to fit a specific purpose such as in \cite{edelmann2020finite} where the flow is chosen  using  knowledge of the surface velocity to construct  a harmonic extension. 
Another common method is  to use a discontinuous or immersed Galerkin method \cite{dg1, dg2, dg3}. In this paper we propose and analyse  an ALE approach using evolving finite elements on  an evolving fitted mesh allowing the use of isoparametric elements that accurately approximate the boundary and result in higher order error estimates. The underlying parametrisation is assumed given.

The key contributions of this work are:
\begin{itemize}
    \item We provide a functional analytic setting to show well posedness of the continuous problem, \cref{eq:strong}.
    \item We provide an ALE approach based on evolving isoparametric finite element spaces attached to evolving subdomains. The  evolving mesh  is based on moving  the Lagrange nodes with a given known  smooth velocity. Achieving a higher order method requires   a good initial mesh.
    \item We provide a robust error bound which demonstrates the error in an $L^2$ norm is bounded, up to a constant, by $h^{k+1}$, where $h$ represents the mesh size and $k$ is the degree of polynomials used both for the discretisation of the domain and the solution. This is the same order error as if we interpolated a known smooth solution.
    \item Numerical results and the simulation code are provided both to demonstrate the results and to allow others to use the implementation.
\end{itemize}

This work builds on the framework of \cite{EllRan21} by applying it to propose and analyse an isoparametric finite element method for a two-sided interface problem and build the necessary theory to achieve an optimal order convergence result.

\begin{remark}
  \label{rem:smooth-mesh}
  \begin{enumerate}
  \item The method requires that we are given a global, smooth velocity field $\mathbf{w}$.  Furthermore, the analysis (\cref{thm:error}) requires that the velocity field is such that moving the nodes of the mesh with the $\mathbf{w}$ preserves the regularity of the mesh over time. The velocity may be derived from problem-dependent considerations or otherwise an arbitrary velocity may be constructed in order to define a well-behaved numerical scheme. We do not address how to achieve such a velocity in this work.
There are methods in the literature to prevent mesh deformation, which involve re-parametrising the flow responsible for the movement of the interface into a more suitable flow, see, for example, \cite{MR2406538, EllFri16,MR3649420}.

\item The method can also be applied to a situation where the evolution of the domain is the solution of a partial differential equation coupled to the solution of the equations we consider in this work. However, our current techniques do not allow us to analyse such a problem. Related results for surface only problems can be found in \cite{KovLiLub17,HuLi22}.
\end{enumerate}
\end{remark}

\subsection{Outline}

\Cref{sec:well-posedness} gives a well posedness analysis of the continuous equations along with the necessary functional analysis setting.
The finite element construction is in \cref{sec:fem} and the finite element scheme is in \cref{sec:fem-scheme}.
An optimal order error bound is shown in \cref{sec:error} under smoothness assumptions on the domain and its evolution and the solution.
\Cref{sec:numerics} includes a time discretisation of the finite element scheme along with numerical experiments demonstrating the error bounds are tight.
The Appendix includes further details of the proof of the well posedness of the continuous scheme.

\section{Evolving space formulation and well posedness}
\label{sec:well-posedness}
In this paper, $c$ will be used as a generic constant that depends on no quantity of particular  importance. We use $(\cdot, \cdot)_{H}$ for an inner product on a Hilbert space $H$ and $\langle \cdot, \cdot \rangle_{X}$ as the dual pairing between a Banach space $X$ and its topological dual $X^*$.  
\subsection{Evolving Hilbert Spaces}
\label{sec:evolving-sobolev}
%Setting up the abstract Sobolev Spaces

We set up the necessary tools from the theory of evolving Sobolev spaces which were introduced and developed in \cite{AlpEllSti15a, AlpEllSti15b, diogo}. We will only concern ourselves with the Hilbert case. A more general theory is developed in \cite{diogo} concerning general Banach spaces. Let $I = [0,T]$ be a closed time interval and let $\{X(t)\}_{t \in I}$ be a family of Hilbert spaces equipped with norm $\|\cdot\|_{X(t)}$. Assume that there exists a linear map $\phi_t:X(0) \to X(t)$ satisfying  the following properties:
\begin{enumerate}
\newitem{$B1$} \label{B1}
  The map $\phi_t$ is invertible for all $t \in I$ with inverse denoted by $\phi_{-t}$ and $\phi_0$ being  the identity.
\newitem{$B2$} \label{B2}
  There exists a constant $C$ independent of time such that $\norm{\phi_t \eta}_{X(t)} \leq C\norm{\eta}_{X(0)}$, $\norm{\phi_{-t} \widetilde{\eta}}_{X(0)} \leq C\norm{\widetilde{\eta}}_{X(t)}$, for all $\eta \in X(0)$ and $\widetilde{\eta} \in X(t)$, for all $t \in I$.
\newitem{$B3$} \label{B3}
  The map $t \mapsto \norm{\phi_t \eta}_{X(t)}$ is measurable  for all $\eta \in X(0)$.
\end{enumerate}
Here and elsewhere we use the notation $\phi_t \eta$ to denote the map $\phi_t$ applied to $\eta$.
 If such a map $\phi_t$ exists then we call it the \textit{flow map} and the pair $(X(t), \phi_t)_{t \in I}$ a \textit{compatible pair}. Given a compatible pair, define the \textit{Hilbert moving spaces} as:
\begin{align}\label{eq:moving-def}
  L^2_{X} :&= \bigg \{ \eta: I \to \bigcup_{t \in I} X(t) \times \{t\}, \quad t \mapsto (\widehat{\eta}(t),t)
             \,|\, \phi_{-t}\widehat{\eta}(t) \in L^2(I; X(0)) \bigg\},
\end{align}
and the uniformly bounded equivalent:
\begin{align}
    L^\infty_{X} :&= \bigg \{ \eta: I \to \bigcup_{t \in I} X(t) \times \{t\}, \quad t \mapsto (\widehat{\eta}(t),t) \,|\, \phi_{-t}\widehat{\eta}(t) \in L^\infty(I; X(0)) \bigg\}.
\end{align}
We identify $\eta(t) = (\widehat{\eta}(t), t)$ with $\widehat{\eta}(t)$. The spaces  $L^p_X$ are equipped with the  norm:
\[
\norm{\eta}_{L^p_X} :=
\begin{cases}
\left({\int_0^T \norm{\hat\eta(t)}_{X(t)}^2}\right)^{\frac{1}{2}} &\text{for } p =2,\\ 
\esssup_{t \in [0,T]} \norm{\eta(t)}_{X(t)} &\text{for } p=\infty.
\end{cases}
\]
$L^2_X$ is indeed a Hilbert space, see \cite[Thm.~3.4]{diogo}.
The analogues of the spaces of continuous functions and of compactly supported smooth functions are defined as:
\begin{align*}
  C^k_X &:= \bigg \{ \eta : I \to \bigcup_{t \in I} X(t) \times \{t\}, \quad t \to (\widehat{\eta}(t),t)
          \,|\, \phi_{-t}\widehat{\eta}(t) \in C^k(I; X(0)) \bigg\}, \\
  \mathcal{D}_X &:= \bigg \{ \eta : I \to \bigcup_{t \in I} X(t) \times \{t\}, \quad t \to (\widehat{\eta}(t),t)
                  \,|\, \phi_{-t}\widehat{\eta}(t) \in \mathcal{D}(I; X(0)) \bigg\}.
\end{align*}
%The spaces $L^p_X$ do not depend on the choice of map $\phi_t$ but $C^k_X$ and $\mathcal{D}_X$ do.

\begin{remark}
\normalfont
\begin{enumerate}
\item
The use of a Cartesian product $\times$  inside the union in \cref{eq:moving-def} rather than just taking the union of $\{X(t)\}_{t \in I}$ by itself is in order to guarantee a disjoint union which is crucial to identify the function point-wise.
\item
Note that the spaces $L^p_X$ do not depend on the choice of the map $\phi_t$. 
\end{enumerate}
\end{remark}
 The \textit{strong material} derivative in the evolving Hilbert space setting is defined as follows:
\begin{align*}
    \partial_t^\bullet \eta := \phi_t \partial_t( \phi_{-t} \eta), \quad \eta \in C^1_X.
\end{align*}

\begin{lemma}[\cite{diogo}, Thm.~2.4]\label{lem:isomo}
    Given a compatible pair $(X(t), \phi_t)_{t \in I}$, the maps $\phi_t : L^2(I; X(0)) \to L^2_X$ and $\phi_{-t}: L^2_X \to L^2(I;X(0))$ define continuous linear isomorphism to their respective spaces.
\end{lemma}
Now assume $\{X(t)\}_{t \in I}$, $\{Y(t)\}_{t \in I}$ and $\{X^*(t)\}_{t \in I}$ are families of Hilbert spaces, with $X^*(t)$ the dual of $X(t)$ for all $t \in I$ (crucially, $X(t)$ and $X^*(t)$ are not identified). Assume further that for all $t \in I$,  $X(t) \subset Y(t) \cong Y^*(t) \subset X^*(t)$ constitutes a Hilbert triple (in the sense that $X(t)$ is densely and continuously embedded into $Y(t)$ and $Y(t)$ is identified with its dual via Riesz representation). It is also assumed that there exists a map $\phi_t:Y(0) \to Y(t)$ with $\phi_{t}|_{X(0)}: X(0) \to X(t)$ with adjoint flow $\phi_{-t}^*:X^*(0) \to X^*(t)$, 
\begin{align*}
    \langle \phi_{-t}^* f, v \rangle_{X(t)} := \langle f, \phi_{-t} v\rangle_{X(0)}, \quad f \in X^*(0), \, v \in X(t),
\end{align*}
such that $(X(t), \phi_t|_{X(t)})_{t \in I}$, $(Y(t), \phi_{t})_{t \in I}$ and $(X^*(t), \phi_{-t}^*)_{t \in I}$ all define compatible pairs and therefore we can define the spaces $L^2_X, \; L^2_Y, \; L^{2}_{X^*}$ with their respective flows. In this case, just as for Bochner spaces, we have that $(L^2_X)^*$ is isometrically isomorphic to $L^{2}_{X^*}$ \cite[Thm.~3.7]{diogo}. Moreover, the Hilbert triple structure is preserved: $L^2_X \subset L^2_Y \subset L^{2}_{X^*}$. Note that $L^2_Y$ remains a Hilbert space with a natural inner product structure, see \cite[Rem.~3.9]{diogo}. In order to generalise the concept of a ``weak time derivative'' to the evolving space, we first assume the following:
\begin{enumerate}
  \newitem{$D1$} \label{D1}
  The map $t \mapsto \langle \phi_t w_0, \phi_t v_0 \rangle_{X(t)} = (\phi_t w_0, \phi_t v_0)_{Y(t)}$ is continuously differentiable for fixed $w_0, v_0 \in X_0$.
  \newitem{$D2$} \label{D2}
  For all $t \in I$, the map:
    \begin{align*}
        [w_0, v_0] \mapsto \frac{d}{dt} (\phi_t w_0, \phi_t v_0)_{Y(t)},
    \end{align*} %\textcolor{magenta}{[Is this really a partial derivative? also in the next couple of lines]}
    for $[w_0, v_0] \in X(0) \times X(0)$ is continuous.
    \newitem{$D3$} \label{D3}
    There exists a constant $C_\lambda$ independent of time such that, for almost all $t \in I$ and $w_0, v_0 \in X(0)$, we have:
    \begin{align*}
        \bigg| \frac{d}{dt} (\phi_t w_0, \phi_t v_0)_{Y(t)} \bigg| \leq C_\lambda \norm{w_0}_{Y(0)} \norm{v_0}_{Y(0)}.
    \end{align*}
\end{enumerate}
\begin{definition}\label{def:lambda}
Let Ass.~\ref{D1} to \ref{D3} hold and label:
\begin{align*}
    \lambda(t; w,v):= \bigg[\frac{d}{dt} (\phi_t w_0, \phi_t v_0)_{Y(t)}\bigg]\bigg|_{(w_0,v_0) = (\phi_{-t} w, \phi_{-t} v)}, \quad w,v \in X(t).
\end{align*}
Then $\lambda(t; \cdot, \cdot): Y(t) \times Y(t) \to \mathbb{R}$ is a continuous, symmetric, bounded  bilinear form for almost all $t \in I$. We say $w \in L^1_X$ has a weak material derivative if there exists $v \in L^1_{X^*}$ such that:
\begin{align*}
    \int_0^T ( w(t), \partial ^\bullet_t \eta)_{Y(t)}\, dt = \int_0^T \langle v(t), \eta \rangle_{X(t)} + \lambda(t; w, \eta)\,dt,
\end{align*}
for all $\eta \in \mathcal{D}_X$. We label the weak material derivative  $v = \partial_t^\bullet w$.
\end{definition}
This definition satisfies all properties we expect from a weak derivative, such as being equivalent to the strong material derivative if the function is regular enough.

This allows us to define the equivalent of the Bochner solution space. 
\begin{definition}We define $W(X,Y):= \left\{v \in L^2_X, \; \partial_t^\bullet v \in L^2_{Y}\right\}$ with the norm:
\begin{align*}
    \norm{v}^2_{W(X,Y)}: = \norm{v}^2_{L^2_X} + \norm{\partial_t^\bullet v}^2_{L^{2}_{Y}}
\end{align*}
and  the solution space $W(X,X^*):= \left\{v \in L^2_X, \; \partial_t^\bullet v \in L^2_{X^*}\right\}$ with the norm:
\begin{align*}
    \norm{v}^2_{W}: = \norm{v}^2_{L^2_X} + \norm{\partial_t^\bullet v}^2_{L^{2}_{X^*}}.
\end{align*}
\end{definition}
%The last assumption needed is the moving space equivalence:
\begin{definition}\label{def:moving-equivalence}The  space $W(X,Y)$ is said to satisfy  moving space equivalence if:
\begin{align*}
    v \in W(X,Y) \iff \phi_{-t} v \in \mathcal{W}^{2,2}(X(0), Y(0)),
\end{align*}
where $\mathcal{W}^{2,2}(X(0),Y(0)) = \{v_0 \in L^2(I;X(0)), \; \partial_t v \in L^{2}(I; Y(0))\}.$
\end{definition}
\begin{theo}[The Transport Theorem] Assume $v,w \in W(X,X^*)$ and the moving space equivalence is satisfied, then, the map $t \mapsto (v,w)_{Y(t)}$ is uniformly continuous and for almost all $t \in I$, and the following holds:
\begin{align*}
    \frac{d}{dt} (v,w )_{Y(t)} = \langle \partial_t^\bullet v,w\rangle_{X(t)} +\langle \partial_t^\bullet w,v\rangle_{X(t)} + \lambda(t;w,v).
\end{align*}
Moreover, $C^0_Y \hookrightarrow  W(X,X^*)$.
\end{theo}
See \cite[Sec.~4.5]{diogo} for proofs.
\begin{lemma}[Characterisation of Material Derivative]\label{lem:equiv-material-derivative}
    Let the moving space equivalence hold for $W(X,X)$, then for $v \in W(X,X)$, we have $v \in C_X^0$ and there exists a function $v_0 \in \mathcal{W}^{2,2}(X(0),X(0))$ such that $v = \phi_t v_0$. Moreover, $C^1_X$ is dense in $W(X,X)$.
\end{lemma}
This follows from \cite[Lem.~3.20]{diogo}. Importantly, this implies $\partial_t^\bullet v=0$ if and only if $v = \phi_t v_0$ for some $v_0 \in X(0)$.
\subsection{Setting up the Domain}\label{sec:settingdom}
\begin{figure}[th]
\fontsize{0.9cm}{0}\resizebox{140mm}{!}{\hbox{\hspace{13em} %% Creator: Inkscape 1.1.2 (08b2f3d93c, 2022-04-05), www.inkscape.org
%% Accompanies image file 'domain.pdf' (pdf, eps, ps)
\begingroup%
  \makeatletter%
  \providecommand\color[2][]{%
    \errmessage{(Inkscape) Color is used for the text in Inkscape, but the package 'color.sty' is not loaded}%
    \renewcommand\color[2][]{}%
  }%
  \providecommand\transparent[1]{%
    \errmessage{(Inkscape) Transparency is used (non-zero) for the text in Inkscape, but the package 'transparent.sty' is not loaded}%
    \renewcommand\transparent[1]{}%
  }%
  \providecommand\rotatebox[2]{#2}%
  \newcommand*\fsize{\dimexpr\f@size pt\relax}%
  \newcommand*\lineheight[1]{\fontsize{\fsize}{#1\fsize}\selectfont}%
  \ifx\svgwidth\undefined%
    \setlength{\unitlength}{680.23745079bp}%
    \ifx\svgscale\undefined%
      \relax%
    \else%
      \setlength{\unitlength}{\unitlength * \real{\svgscale}}%
    \fi%
  \else%
    \setlength{\unitlength}{\svgwidth}%
  \fi%
  \global\let\svgwidth\undefined%
  \global\let\svgscale\undefined%
  \makeatother%
  \begin{picture}(1,0.48806137)%
    \lineheight{1}%
    \setlength\tabcolsep{0pt}%
    \put(0,0){\includegraphics[width=\unitlength,page=1]{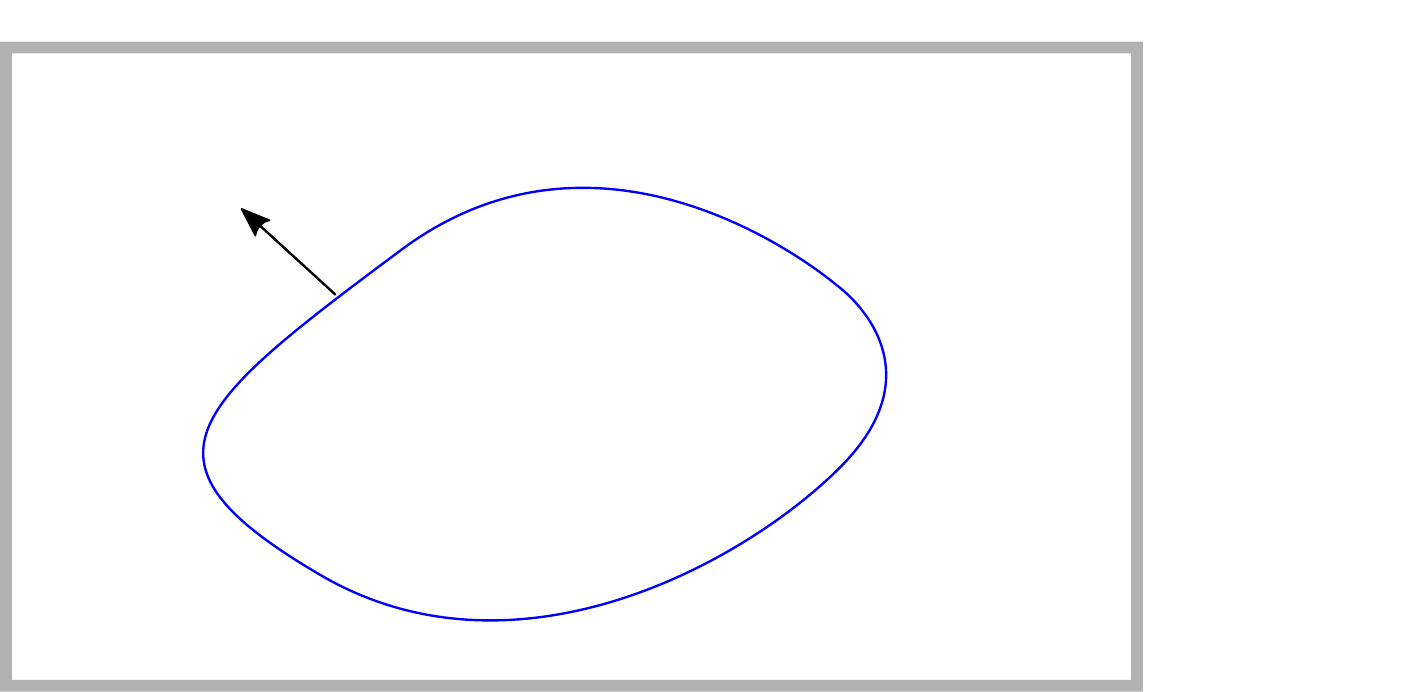}}%
    \put(0.37588239,0.19267144){\color[rgb]{0,0,0}\makebox(0,0)[lt]{\lineheight{1.25}\smash{\begin{tabular}[t]{l}$\Omega_1(t)$\end{tabular}}}}%
    \put(0.52362125,0.07248429){\color[rgb]{0,0,0}\makebox(0,0)[lt]{\lineheight{1.25}\smash{\begin{tabular}[t]{l}$\Gamma(t)$\end{tabular}}}}%
    \put(0.20283995,0.32919895){\color[rgb]{0,0,0}\makebox(0,0)[lt]{\lineheight{1.25}\smash{\begin{tabular}[t]{l}$\nu(t)$\end{tabular}}}}%
    \put(0.67104309,0.36513718){\color[rgb]{0,0,0}\makebox(0,0)[lt]{\lineheight{1.25}\smash{\begin{tabular}[t]{l}$\Omega_2(t)$\end{tabular}}}}%
    \put(0.00416714,0.45455411){\color[rgb]{0,0,0}\makebox(0,0)[lt]{\lineheight{1.25}\smash{\begin{tabular}[t]{l}$\partial \Omega$\end{tabular}}}}%
  \end{picture}%
\endgroup%
}}
\caption{An example configuration of the domain.}
\label{fig:domain}
\end{figure}
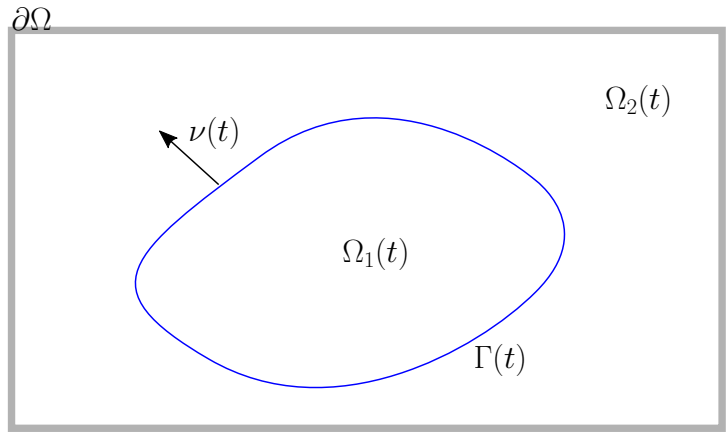
Let $\Omega$ be a stationary domain in $\mathbb{R}^d$, $d=2,3$, with piecewise linear boundary and let $\{\Gamma(t), \; t \in I\}$ be a family of closed compact connected $C^{2+k}$ ($k\geq 0$) hypersurfaces with $\Gamma(t) \subset \Omega$.
Let $\Omega_1(t)$ be a domain in $\Omega$ without boundary $\partial \Omega_1(t) = \Gamma(t)$ for all $t \in I$.
Let $\Omega_2(t): = \Omega \setminus \Omega_1(t)$ and assume that $\Gamma(t) \cap \partial \Omega = \varnothing$ for all $t \in I$, then:
\begin{align*}
    \overline{\Omega} = \overline{\Omega}_1(t) \cup \overline{\Omega}_2(t), \; \overline{\Omega}_1(t) \cap \overline{\Omega}_2(t) = \Gamma(t), \, \partial\Omega_2(t) = \Gamma(t) \cup \partial \Omega.
\end{align*}
A sketch of the domains is shown in \cref{fig:domain}.

\begin{remark}
\normalfont
The assumption that the outer boundary is piecewise linear is made to avoid having to analyse perturbation of the domain for Dirichlet boundary conditions, however the presented  method and analysis can easily be altered if one removes this assumption.
\end{remark}

We label the outer normals of $\Omega_1(t)$ and $\Omega_2(t)$ by $\nu_{\Gamma(t)}$ and $\nu_{\partial \Omega_2(t)}$ respectively.
Let:
\begin{align*}
    \mathcal{Q}_i := \bigcup_{t \in I} \Omega_i(t) \times \{t\}, \; \mathcal{Q}:= \Omega \times I.
\end{align*}
Furthermore, we assume there exists a given, global velocity field $\mathbf{w}$ transporting $\Omega_1(t)$ and $\Omega_2(t)$, i.e $ \mathbf{w} \cdot \nu_{\Gamma(t)}|_{\Gamma(t)}= V_{\Gamma}$ where $V_\Gamma$ is the normal velocity of $\Gamma(t)$, $\mathbf{w} \cdot \nu_{\partial \Omega_2(t)}|_{\partial \Gamma(t)} = -V_{\Gamma}$ and $\mathbf{w} \cdot \nu_{\partial \Omega_2(t)}|_{\partial \Gamma} = 0$. Throughout the paper this velocity is assumed to be of regularity $\mathbf{w} \in C(I; C(\overline{\Omega}; \mathbb{R}^d))$ with $\mathbf{w}_i(t; \cdot) \in C^{2}(\overline{\Omega}_i(t); \mathbb{R}^d)$. Let $\mathbf{\Phi}_i(t;x): \overline{\Omega}_i(0) \to \overline{\Omega}_i(t)$ be the solution to the ordinary differential equations:
\begin{align}\label{eq:flow-def}
    \frac{d}{dt} \mathbf{\Phi}_i(t;x) &=\mathbf{w}(t; \mathbf{\Phi}_i(t;x)) \qquad x \in \Omega_i(0),\\
    \mathbf{\Phi}_i(0; x) &= x.\nonumber
\end{align}
From the regularity of $\mathbf{w}$, we know that $\mathbf{\Phi}_i$ exist and are of regularity $\mathbf{\Phi}_i \in C^{1}(\overline{\mathcal{Q}}_i; \mathbb{R}^d)$ with $\mathbf{\Phi}_i(t;\cdot)\colon \overline{\Omega}_i(0) \to \overline{\Omega}_i(t)$ and $\mathbf{\Phi}_i(t;\cdot) \in C^{2}(\Omega_i(0); \mathbb{R}^d)$, see \cite[Thm.~1.45]{Rou05} and \cite[Thm.~II.1.1, Sec.~V]{hartman1982ordinary} for the necessary additional conditions. Furthermore, both $\mathbf{\Phi}_i(t; \cdot)$ are invertible diffeomorphisms for all $t \in I$ with $\mathrm{Im}(\mathbf{\Phi}_i(t;\cdot)) = \Omega_i(t)$.
We denote by $\Phi_i(-t, \cdot)$ the inverse of $\Phi_i(t, \cdot)$.
Since we assumed $\mathbf{w} \in C(I; C(\overline{\Omega}; \mathbb{R}^d))$, it follows that $\mathbf{\Phi}_1(t;x)|_{\Gamma(0)} = \mathbf{\Phi}_2(t;x)|_{\Gamma(0)}$.

\begin{remark}
\normalfont
For the abstract formulation of the problem, it is only required to assume that the velocity field $\mathbf{w}$ is of sufficient regularity. However, for the purpose of evolving the mesh later, we will assume this velocity field is known explicitly. 
\end{remark}

Let $J^t_i$ denote the determinant of Jacobian matrix, $J^t_i = \text{det} [\nabla \mathbf{\Phi}_i(t;x)]$. The prior assumptions imply $J^t_i(\cdot) \in C^{1}(\overline{\Omega}_i(0); \mathbb{R})$ and there exists $C_\Omega$ independent of time and space such that:
\begin{align*}
    \frac{1}{C_\Omega} \leq |J_i^t | \leq C_\Omega,
\end{align*}
and $J^{-t}_i$ denotes its inverse.
\begin{remark}\label{rem:normal-vel}
\normalfont
Note that we have assumed the global parametric velocity $\mathbf{w}$ is given.
The solution of the partial differential equation system is independent of $\mathbf{w}$ apart from the requirements that $\mathbf{w}_\Gamma \cdot \nu_{\Gamma(t)} = V_\Gamma$ and $\mathbf{w} \cdot \nu_{\partial \Omega} = 0$. However the evolving mesh does depend on $\mathbf{w}$ hence the discrete solution depends on the full parametric velocity.
\end{remark}

Let $d_\Gamma(t;x)$ be the signed distance function to $\Gamma(t)$:
\begin{align*}
   d_\Gamma(t;x) = \begin{cases}
   -\inf\{|x - y|: y \in \Gamma(t)\}, \quad &\text{for} \, x \in \overline{\Omega}_1(t), \\
   \inf\{|x - y|: y \in \Gamma(t)\}, \quad &\text{for} \, x \in \Omega_2(t) .
   \end{cases} 
\end{align*}
Then, since the interface is of class $C^{2}$, there exists a constant  $\delta >0$ such that if $ x \in \mathcal{N}_{\Gamma(t)} := \{x \in \Omega. \; | d_{\Gamma}(t;x)| \leq \delta \}$, it can be uniquely decomposed as:
\begin{align}\label{eq:tubular-decom}
    x = d_{\Gamma}(t;x) \nu_{\Gamma(t)}(\Pi_{t}(x)) + \Pi_t(x),
\end{align}
where $\Pi_t(\cdot)$ is the nearest point on $\Gamma(t)$, i.e, $\Pi_t(x) := \inf \{|y - x| : y \in \Gamma(t)\}$ (see \cite[Sec.~2.3]{pruss2016moving}). We refer to the set $\mathcal{N}_{\Gamma(t)}$ as the \textit{tubular neighbourhood} of $\Gamma(t)$. Note that $\delta$ can be chosen independently of time by the fact that $I$ is compact. Moreover, via the assumption that $\Gamma(t) \cap \partial \Omega = 0$ for all $t \in I$, $\delta$ can be chosen small enough such that $\mathcal{N}_{\Gamma(t)} \cap \partial \Omega = \varnothing$ for all $t \in I$.

For the error analysis, we require further  regularity of the flow to yield the results collected in the following lemma.  
\begin{lemma}\label{lem:geom-regularity}
Let $\Theta \in \mathbb{N}$. Assume further regularity on the flow map $\mathbf{\Phi} \in C^{2 +\Theta}\left(I; C^{2 +\Theta}\left(\overline{\Omega}_i(0); \mathbb{R}^d\right)\right)$ and the initial surface $\Gamma(0)$ is class $C^{2 + \Theta}$, then the following geometric quantities have additional regularity:
\begin{itemize}
    \item $\Gamma(t)$ is of class $C^{2 + \Theta}$,
    \item $J^t_i(\cdot) \in C^{1 + \Theta}(\overline{\Omega}_i(0); \mathbb{R})$,
    \item $d_\Gamma(t;\cdot) \in C^{2+\Theta}(\mathcal{N}_{\Gamma(t)}; \mathbb{R})$,
    \item $\Pi_t(\cdot) \in C^{1+\Theta}(\mathcal{N}_{\Gamma(t)}; \mathbb{R}^d)$.
\end{itemize}
\end{lemma}
See \cite{signed} and \cite[Lem.~14.16]{GilTru98}.
%Diagram, talk about the flow, velocity of the interface and set up the appropriate regularity for the flow., introduce the tubular neighbourhood theorem and the assumptions required for appropriate regularity. 

\subsection{Realisation}\label{sec:realisation}
%Setting the particular ones we need (H^1, H^{-1}, L^2, H^{1/2}, etc)
For a given function $v$ acting on $\Omega$, we decompose it as:
\begin{align*}
    v_1  := v \chi_{\overline{\Omega_1(t)}}, \quad v_2 :=  v \chi_{\overline{\Omega_2(t)}},
\end{align*}
where $\chi_{\overline{\Omega_i(t)}} = 1$ if $x \in \overline{\Omega_i(t)}$ and zero otherwise. A function $v$ on $\Omega$ will be identified as the pair $v = (v_1, v_2)$. The jump operator $\llbracket \cdot \rrbracket_{\Gamma(t)}: C(\overline{\Omega_1(t)}; \mathbb{R}) \times C(\overline{\Omega_2(t)};\mathbb{R}) \to C(\Gamma(t); \mathbb{R})$ % \textcolor{magenta}{[I think this is a funny way to write this, can we just write the $H^1$ version?]} 
as:
\begin{align*}
    \llbracket v \rrbracket_{\Gamma(t)} := [v_1 - v_2]|_{\Gamma(t)}.
\end{align*}
This functional has a natural extension on the Cartesian product of standard Sobolev spaces $H^1(\Omega_1(t)) \times H^1(\Omega_2(t))$ via use of the trace maps; $T_i(t):H^1(\Omega_i(t)) \to H^{1/2}(\partial \Omega_i(t))$, $ i \in \{1,2\}$, (see \cite[Sec.~7.2.5]{MR2028503} for an extensive definition of the trace map) as:
\begin{align*}
      \llbracket v \rrbracket_{\Gamma(t)} := [T_1(t) v_1 - T_2(t) v_2]|_{\Gamma(t)}.
\end{align*}
We define the following spaces:
\begin{align*}
    H(t) &:= L^2(\Omega_1(t)) \times L^2(\Omega_2(t)), \quad \norm{v}_{H(t)}^2 := \sum_{i = 1}^2 \norm{v_i}^2_{L^2(\Omega_i(t))},  \\
    V(t) &:= \{v \in H^1(\Omega_1(t)) \times H^1(\Omega_2(t)),\, \llbracket v \rrbracket_{\Gamma(t)} = 0, \; v_2|_{\partial \Omega} = 0\}, \quad \norm{v}_{V(t)}^2 := \sum_{i = 1}^2 \norm{v_i}^2_{H^1(\Omega_i(t))}\\
     Z_k(t) & := \{w \in V(t)|\; w_i \in H^{1+k}(\Omega_i(t))\}, \; \norm{w}^2_{Z_k(t)} := \sum_{i = 1}^2 \norm{w_i}^2_{H^{1+k}(\Omega_i(t))}.
\end{align*}
 
Note that due to the continuity of the trace operators $T_i(t)$ on $H^1(\Omega_i(t))$, $V(t)$ defines a closed subspace of $H^1(\Omega_1(t)) \times H^1(\Omega_2(t))$ and contains $H_0^1(\Omega_1(t)) \times H_0^1(\Omega_2(t))$, hence is dense within $H(t)$ . For an element $v = (v_1, v_2) \in V(t)$, we will identify $v|_{\Omega_1(t)} = v_1$ and $v|_{\Omega_2(t)} = v_2$.
We also define the interface space:
\begin{align*}
    H^{1/2}(\Gamma(t)) &= \{ v \in L^2(\Gamma(t)), \, |v|_{H^{1/2}(\Gamma(t))} < \infty\}, \; |v|_{H^{1/2}(\Gamma(t))} := \int_{\Gamma(t)}\int_{\Gamma(t)} \frac{|v(x) - v(y)|^2}{|x - y|^d }  \, dxdy,  \\
  \intertext{with norm given by:}
    \norm{v}^2_{H^{1/2}(\Gamma(t))}:&= \norm{v}^2_{L^2(\Gamma(t))} + |v|^2_{H^{1/2}(\Gamma(t))}.
\end{align*}
Then, $H^{1/2}(\Gamma(t))$ is a Hilbert space and moreover is dense and compactly embedded in $L^2(\Gamma(t))$ (see \cite[Sec.~2]{MIKHAILOV2011324}). For consistency of notation, let $\mathcal{V}_\Gamma(t) = H^{1/2}(\Gamma(t))$ and $\mathcal{H}_\Gamma(t) = L^2(\Gamma(t))$, and identify the Hilbert  triple $\mathcal{V}_\Gamma(t) \subset \mathcal{H}_\Gamma(t) \subset \mathcal{V}_\Gamma^*(t)$. 

Now for a function $v \in H(t)$ and $w \in \mathcal{H}_\Gamma(t)$, the respective flows are defined as:
\begin{align*}
    \phi_t v := (v_1(t;\mathbf{\Phi}_1(-t;x)), v_2(t;\mathbf{\Phi}_2(-t;x))), \quad \phi_t w := w(t;\mathbf{\Phi}_1(-t;x)).
\end{align*}

\begin{lemma}\label{lem:pf-compat}
The pairs $(V(t), \phi_t)|_{t \in I}$, $(H(t), \phi_t)|_{t \in I}$, $(V^*(t), \phi_{-t}^*)|_{t \in I}$, $(\mathcal{V}_\Gamma(t), \phi_t)|_{t \in I}$, $(\mathcal{H}_\Gamma(t), \phi_t)|_{t \in I}$ and $(\mathcal{V}_\Gamma^*(t), \phi_{-t}^*)|_{t \in I}$ are all compatible. 

Assuming the added regularity $\mathbf{\Phi}_i(t; \cdot) \in C^{1 +k}(\Omega_i(0); \mathbb{R}^d))$, 
the pair $(Z_k(t), \phi_t)|_{t \in I}$ is compatible. \end{lemma}
\begin{proof}
Ass.~\ref{B1} to \ref{B3} need to be checked. This will be checked only for $V(t)$ as a similar logic can be employed for the remaining spaces. % and by using the same argument as Proposition 2.17 in \cite{djurdjevac2021non} \textcolor{magenta}{[TR: this phrase is missing a verb, I think, but not sure how to fix]}.
 \ref{B1} follows from both $\mathbf{\Phi}_1(-t;x)$, $\mathbf{\Phi}_2(-t;x)$ being invertible diffeomorphisms.
For \ref{B2}, via simple manipulation:
\begin{align}\label{compatibility}
   \norm{\phi_t v}^2_{V(t)}  &= \sum_{i = 1}^2 \int_{\Omega_i(t)} |v_i(t;\mathbf{\Phi}_i(-t;x))|^2 + |\nabla v_i(t;\mathbf{\Phi}_i(-t;x))|^2, \nonumber \\
   &= \sum_{i = 1}^2 \int_{\Omega_i(0)}[ \,|v_i(t;x)|^2 + |[\nabla \mathbf{\Phi}_i(-t;y)]^T|_{y = \mathbf{\Phi}(t;x)} \nabla v_i(t;x)|^2]J^t_i \\
   &\leq c(|J^t_i|_{L^\infty(\Omega_i(0))}, |\nabla \mathbf{\Phi}_i(-t,x)|_{L^\infty(\Omega_i(t))}) \norm{v}^2_{V(0)}, \nonumber
\end{align} 
 The bound follows from the assumption on the regularity of the velocity field. The same method shows a similar bound for $\norm{\phi_{-t} \widetilde{v}}_{V(0)} \leq c\norm{\widetilde{v}}_{V(t)}$ for all $\widetilde{v} \in V(t)$. To show measurability, \ref{B3}, note that the second equality in \cref{compatibility} is continuous. For the compatibility of the boundary spaces $\mathcal{V}_\Gamma(t)$, $\mathcal{H}_\Gamma(t)$ and $\mathcal{V}^*_\Gamma(t)$, see \cite[Sec.~4 and 5]{AlpEllSti15b}.
 Under the added regularity the compatibility of $(Z_k(t), \phi_t)|_{t \in I}$ follows similarly.
\end{proof}
We will identify both $V(t) \subset H(t) \subset V^*(t)$ and $\mathcal{V}_\Gamma(t) \subset \mathcal{H}_\Gamma(t) \subset \mathcal{V}_\Gamma^*(t)$ with the structure $X(t) \subset Y(t) \subset X^*(t)$ developed in \cref{sec:evolving-sobolev}.
\begin{lemma}The moving space equivalence is satisfied between $W(V,V^*)$ and $\mathcal{W}(V(0), V^*(0))$.
\end{lemma}
\begin{proof}
The proof follows similarly from the one given for \cite[Prop.~7.4]{diogo} as by assumption the Jacobian determinate is at least of regularity $J^t_i \in C^1 (I; C^1(\overline{\Omega}_i(0); \mathbb{R}))$.
\end{proof}
\begin{remark}
\normalfont
It does not matter which of the flows $\mathbf{\Phi}_i(t;x)$ is used to define $\mathcal{H}_\Gamma(t)$ as $\mathbf{\Phi}_1(t;x)|_{\Gamma(0)} = \mathbf{\Phi}_2(t;x)|_{\Gamma(0)}$. Moreover, it can be shown that $v = (v_1, v_2) \in V(t)$ if, and only if, $v_1\chi_{\Omega_1(t)} + v_2(1 - \chi_{\Omega_1(t)}) \in H^1_0(\Omega)$ with equivalent norms, hence the space $V(t)$ can be thought as an identification of the components of functions in $H^1_0(\Omega)$. 
\end{remark}

We may define both moving space triples $L^2_V \subset L^2_H \subset L^2_{V^*}$ and $L^2_{\mathcal{V}_\Gamma} \subset L^2_{\mathcal{H}_\Gamma} \subset L^2_{\mathcal{V}^*_\Gamma}$.
\begin{theo}[Reynolds' Transport Theorem] Let $g_i \in C^1(\mathcal{Q}_i; \mathbb{R})$, then:
\begin{align*}
    \frac{d}{dt} \sum_{i = 1}^2\int_{\Omega_i(t)} g_i = \sum_{i = 1}^2\int_{\Omega_i(t)} \partial_t g_i +  \mathbf{w} \cdot \nabla g_i  + g_i \nabla \cdot \mathbf{w}.
\end{align*}
\end{theo}
\begin{proof}
We use another version of Reynolds' Transport Theorem given in \cite[Sec.~2.5]{pruss2016moving}. For $g = (g_1, g_2) \in C^1(\mathcal{Q}_1; \mathbb{R}) \times C^1(\mathcal{Q}_2; \mathbb{R})$, then:
\begin{align*}
     \frac{d}{dt} \int_{\Omega \setminus \Gamma(t)} g =  \int_{\Omega \setminus \Gamma(t)} \partial_t g - \int_{\Gamma(t)} \llbracket g\rrbracket_{\Gamma(t)} V_\Gamma.
\end{align*}
Note that here $\Omega \setminus \Gamma(t) = \Omega_1(t) \cup \Omega_2(t)$, and:
\begin{align*}
    \sum_{i = 1}^2 \int_{\Omega_i(t)}  \mathbf{w} \cdot \nabla g_i  + g_i \nabla \cdot \mathbf{w} =  \sum_{i = 1}^2  \int_{\Omega_i(t)} \nabla \cdot( \mathbf{w} g_i) &= - \int_{\partial \Omega} [\mathbf{w} g_i] \cdot \nu_{\partial \Omega} - \int_{\Gamma(t)} \llbracket \mathbf{w} g \rrbracket_{\Gamma(t)} \cdot \nu_{\Gamma(t)},\\
    &=  - \int_{\Gamma(t)} \llbracket  g \rrbracket_{\Gamma(t)}  V_\Gamma.
    \qedhere
\end{align*}
\end{proof}
Note that via use of the chain rule and the definition of $\mathbf{\Phi}(t; \cdot)$ \cref{eq:flow-def}, for a function $\eta_i \in C^1(\mathcal{Q}_i; \mathbb{R})$:
\begin{align*}
  \bigg[\frac{d}{dt}(\eta_i(t;\mathbf{\Phi}_i(t;x))\bigg]\bigg|_{x = \mathbf{\Phi}_i(-t;y)} &= \partial_t \eta_i(t; y) + \partial_t(\mathbf{\Phi}_i(t;x))|_{x = \mathbf{\Phi}_i(-t;y)} \cdot \nabla \eta_i(t;x),\\
  &= \partial_t \eta_i(t; y)+  \mathbf{w}(t;y)\cdot \nabla \eta_i(t;x).
\end{align*}
Giving us back the classical definition for the material derivative, see \cite[Sec~1.1.1]{GrosReus}. For a function $v \in C^1_V$, we define:
\begin{align*}
    \partial_t^\bullet v = \phi_t \frac{d}{dt}(v_1(t;\mathbf{\Phi}_1(t;x)), v_2(t;\mathbf{\Phi}_2(t;x))) = ([\partial_t +\mathbf{w} \cdot \nabla]v_1(t;x), [\partial_t +\mathbf{w} \cdot \nabla]v_2(t;x)) =: (\partial_t^\bullet v_1, \partial_t^\bullet v_2).
\end{align*}
One can check that due to the regularity of the flow, the assumptions \ref{D1} to \ref{D3} are satisfied on the triple $V(t) \subset H(t) \subset V^*(t)$, moreover, via Reynold's transport theorem, one can check that the bilinear form $\lambda$ introduced \cref{sec:evolving-sobolev} in this case becomes:
\begin{align}\label{eq:lambda-1}
    \lambda(t;v, \eta) = (\nabla \cdot \mathbf{w} v, \eta)_{H(t)} = \sum_{i = 1}^2 \int_{\Omega_i(t)}[\nabla \cdot \mathbf{w} ]v_i\,\eta_i.
\end{align}
See \cite[Lem.~6.3]{diogo} for more details.

\subsection{The Weak Formulation}
%derive the weak form + remarks on the B term. introduce all the require bilinear form derivatives.
% \textcolor{magenta}{[Check against implementation;\\
%   Notes about weak formulations: I would recommend this formulation to be consistent with big paper
%   \begin{align*}
%     \partial_t u - \nabla \cdot ( A \nabla u ) + B \cdot \nabla u + c u & = f && \mbox{ on } \Omega \\
%     [A \nabla u \cdot \nu ] & = G && \mbox{ in } \Gamma \\
%     \text{etc}.
%   \end{align*} 
%   Weak form:
%   \begin{align*}
%     \int_\Omega \partial^\bullet u v + u v \nabla \cdot w + A \nabla u \cdot \nabla v + (B - w) \cdot \nabla u v + (c - \nabla \cdot w) u v
%     & = \int_\Omega f v + \int_\Gamma g v && \mbox{ for all } v \\
%     m(t; \partial^\bullet u, v) + \lambda(t; u, v) + a(t; u, v) & = l(t; v).
%   \end{align*}
%   Conservative weak form:
%  \begin{align*}
%   \frac{d}{dt} \left(\int_\Omega u v \right)
%   + \int_\Omega A \nabla u \cdot \nabla v + (B - w) \cdot \nabla u v + (c - \nabla \cdot w) u v
%   & = \int_\Omega u \partial^\bullet v +  f v + \int_\Gamma g v && \mbox{ for all } v \\
%   \frac{d}{dt} m(t; u, v) + a(t; u, v) = m(t; u, \partial^\bullet v) + l(t; v).
%   \end{align*}
% ]}
Taking the strong problem \cref{eq:strong}, assuming there exists a regular enough solution $u$, we can rewrite the partial differential equation as:
\begin{align}
   \partial_t^{\bullet} u_i + \nabla \cdot \mathbf{w}u_i - \nabla \cdot (\mathcal{A}_i (t;x) \nabla u_i)  +[\mathcal{B}_i(t;x)- \mathbf{w}] \cdot \nabla u_i +[\mathcal{C}_i(t;x) - \nabla \cdot \mathbf{w}]u_i=  f_i.
\end{align}
Here the term $\nabla \cdot \mathbf{w}$, corresponding to the previously identified bilinear form $\lambda(t; \cdot,\cdot)$, \cref{eq:lambda-1}, is introduced to get the equation in a more convenient form. Then testing with a function $v \in L^2_V$ and using the interface condition, we arrive at the following variational problem:
\begin{align*}
    \int_0^T \langle \partial_t^{\bullet} u, v \rangle_{V(t)}\,dt +\int_0^T \lambda(t;u,v) + \underbrace{\sum_{i = 1}^2\int_{\Omega_i(t)} \mathcal{A}_i \nabla u_i \cdot \nabla v_i + [\mathcal{B}_i - \mathbf{w}] \cdot \nabla u_i \,v_i + [\mathcal{C}_i - \nabla\cdot\mathbf{w}]u_i\, v_i}_{=:a(t;u,v)}  \,dt\\
    = \int_0^T \underbrace{\langle f, v\rangle_{V(t)} +\langle g,v \rangle_{\mathcal{V}_\Gamma(t)}}_{=:l(t;v)}\,dt. 
\end{align*}
Note that here the Hilbert triple structure is used for the duality pairings; $\langle f, v\rangle_{V(t)} = (f,v)_{H(t)}$ and $\langle g,v \rangle_{\mathcal{V}_\Gamma(t)} = (g,v)_{\mathcal{H}_\Gamma(t)}$, and we have the initial condition $u(0) = u_0$. This gives us the weak formulation:
\begin{align}\label{variational}
     \int_0^T \langle \partial_t^\bullet u, v \rangle_{V(t)} +a(t;u,v) +\lambda(t;u,v) \,dt = \int_0^T l(t;v)\,dt,
\end{align}
for all $v \in L^2_V$. Moreover if, instead $v \in W(V,V^*)$, we get the equivalent formulation via the transport theorem (for notational convenience later on, we will label the inner product $(\cdot, \cdot)_{H(t)} =: m(t;\cdot,\cdot)$):
\begin{align*}
    \frac{d}{dt} m(t; u, v) + a(t; u, v) = m(t; u, \partial^\bullet_t v) + l(t; v).
\end{align*}
If $\partial_t^\bullet u \in L^2_H$, then via identification of the Hilbert triple, we have:
\begin{align*}
    \langle \partial_t^\bullet u ,v \rangle_{V(t)} = m(t; \partial_t^\bullet u, v),
\end{align*}
and the problem can be restated abstractly in this case as $u \in W(V,H)$ being the solution to:
\begin{align}\label{eq:abstract}
    m(t; \partial_t^\bullet u, v) + a(t;u,v) + \lambda(t; u,v)= l(t;v),
\end{align}
for almost all $t \in I$, and all $v \in L^2_V$.

\subsection{Well Posedness} ~

\begin{theo}\label{theorem:existence}
Assume the following:
\begin{enumerate}
    \newitem{$A1$} \label{A1}
    The coefficients $\mathcal{A}_i \in C(\overline{\mathcal{Q}}_i; \mathbb{R}^{d \times d})$, $\mathcal{B}_i \in C(\overline{\mathcal{Q}}_i; \mathbb{R}^d)$ and $\mathcal{C}_i \in C(\overline{\mathcal{Q}}_i; \mathbb{R})$;
        \newitem{$A2$} \label{A2}
        There exists a constant $\gamma >0$ such that:
    \begin{align}
      \label{eq:Agamma}
    \inf_{t \in I} \inf_{x \in \Omega_i(t)} \mathcal{A}_i(t;x) \xi \cdot \xi^T \geq \gamma |\xi|^2 \quad \forall \xi \in \mathbb{R}^d;
    \end{align}
        \newitem{$A3$} \label{A3}
        $(u_0, f, g, \mathbf{w}_i, \mathbf{w}) \in H(0) \times L^2_{V^*} \times L^2_{V_\Gamma^*} \times C^1(\overline{\mathcal{Q}}_i, \mathbb{R}^d) \times C(I \times \overline{\Omega}; \mathbb{R}^d))$,
\end{enumerate}
 then there exists a unique solution $u \in W(V,V^*)$ to \cref{variational} with inequality:
\begin{align*}
    \norm{u}_{W} \leq C \bigg(\norm{f}_{L^2_{V^*}} + \norm{g}_{L^2_{\mathcal{V}_\Gamma^*}} +\norm{u_0}_{H(0)}\bigg).
\end{align*}
Furthermore, if it holds that:
\begin{enumerate}
    \newitem{$A4$} \label{A4}
    $(u_0,f,g, \mathcal{A}_i) \in V(0) \times  L^2_H  \times W(\mathcal{V}_\Gamma, \mathcal{V}^*_\Gamma) \times C^1(\overline{\mathcal{Q}_i}; \mathbb{R}^{d \times d})$, and $\mathcal{A}_i$ is symmetric,
\end{enumerate}
then the solution is of additional regularity $u \in W(V,H)$ with bound:
\begin{align*}
    \norm{u}_{W(V,H)} \leq C\bigg(\norm{f}_{L^2_H} +\norm{u_0}_{V(0)} +\norm{g}_{W(\mathcal{V}_\Gamma, \mathcal{V}^*_\Gamma)} \bigg).
\end{align*}
\end{theo}
\begin{proof}
The existence and uniqueness follows from a standard application of the Babuska-Lax-Milgram theorem in conjunction with Poincar\'e's inequality, detailed in \cite[Thm.~3.6]{AlpEllSti15a}. The proof of additional regularity under Ass.~\ref{A4} is given in the Appendix (\cref{lem:extra-reg}).
\end{proof}
Furthermore, in order to analyse the error in the finite element approximation of the material derivative,  it is convenient to define notation for the derivative of the bilinear form $a(t;\cdot, \cdot)$ to be:
\begin{align}\label{eq;derive}
    b(t;v,w) := \frac{d}{dt} [a(t;v,w)] - a(t; \partial^\bullet_t v,w) - a(t;v, \partial^\bullet_t w), \quad \forall v,w \in W(V,V).
\end{align}
Then, assuming furthermore that $\mathcal{A}_i \in C^1(\overline{\mathcal{Q}}_i; \mathbb{R}^{d \times d})$, $\mathcal{B}_i \in C^1(\overline{\mathcal{Q}}_i; \mathbb{R}^d)$ and $\mathcal{C}_i \in C^1(\overline{\mathcal{Q}}_i; \mathbb{R})$, the bilinear form $b(t; \cdot, \cdot)$ exists and can be explicitly calculated as:
\begin{multline}
\label{eq:diff}
  b(t;v,w) = \sum_{i = 1}^2 \int_{\Omega_i(t)} \mathcal{D}^\mathcal{A}_i(\mathbf{w}, \mathcal{A}_i, v_i,w_i) + \mathcal{D}^\mathcal{B}_i(\mathbf{w},\mathcal{B}_i, v_i,w_i) \\
  \qquad\qquad\qquad\qquad +  v_i\,w_i  \partial^\bullet_t [\mathcal{C}_i - \nabla \cdot \mathbf{w}] +  \nabla \cdot \mathbf{w}[\mathcal{C}_i - \nabla \cdot \mathbf{w}] v_i\,w_i,
\end{multline}
where
\begin{align*}
  \mathcal{D}^\mathcal{A}_i(\mathbf{w}, \mathcal{A}_i, v_i, w_i) &= (\partial^\bullet_t \mathcal{A}_i(t;x)+ \nabla \cdot \mathbf{w}\mathcal{A}_i(t;x))\nabla v_i \cdot \nabla w_i  - 2 D_i(\mathbf{w}, \mathcal{A}_i) \nabla v_i \cdot \nabla w_i,\\
  \mathcal{D}^\mathcal{B}_i(\mathbf{w},\mathcal{B}_i,v_i,w_i) &= \partial^\bullet_t [\mathcal{B}_i(t;x) - \mathbf{w}] \cdot \nabla v_i\, w_i + [\mathcal{B}_i(t;x) - \mathbf{w}]\cdot\nabla v_i\, w_i \nabla \cdot \mathbf{w} \\
  & \qquad\qquad - \sum_{j,k = 1}^{d} [\mathcal{B}_i - \mathbf{w}]_j  (\nabla_j \mathbf{w}_k)\nabla_k v_i\,w_i,\\
  [D_i(\mathbf{w}, \mathcal{A}_i)]_{jl}
  &= \frac{1}{2} \sum_{r = 1}^d [\mathcal{A}_i(t;x)]_{jr} \nabla_r \mathbf{w}_l
  + [\mathcal{A}_i(t; x)]_{lr} \nabla_r \mathbf{w}_j.
\end{align*}
Note that the derivative of the bilinear form $m(t;\cdot,\cdot)$ is already assumed to exist and equals $\lambda(t;\cdot,\cdot)$ introduced in \cref{sec:evolving-sobolev}.
%existence and uniqueness, assume regularity

\section{Evolving finite elements}\label{sec:fem}
\begin{figure}[ht]
\fontsize{8}{10}\selectfont
\centering{
\resizebox{150mm}{!}{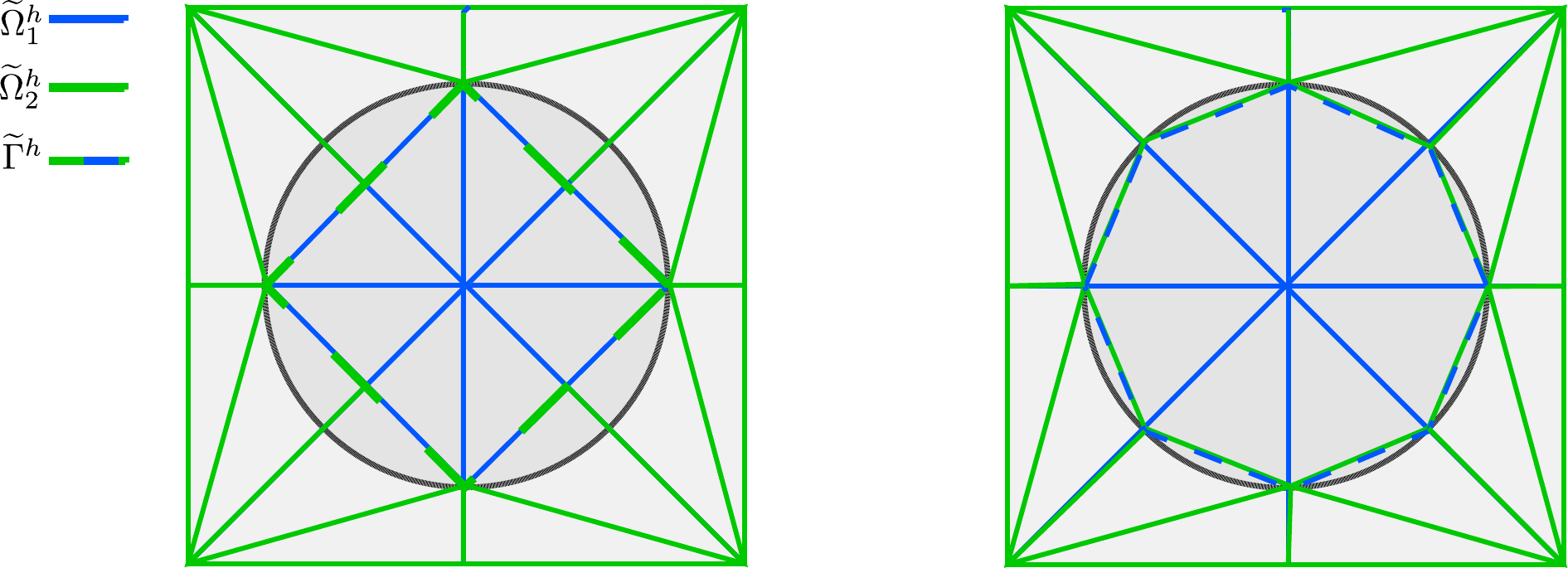}
\caption{Showing the difference between a non viable initial mesh and an adequate one for a circle enclosed in a square. The left initial mesh breaks condition \ref{M5} as not all vertices lie on the interface. The right initial mesh satisfies \ref{M5}.}
\label{fig:topView}
}
\end{figure}

From this point on, we assume the additional  geometric regularity as described in \cref{lem:geom-regularity}. We begin by detailing the initial triangulation of the domain and follow with the construction of the evolving mesh. In order to relate discrete and continuous functions we introduce the concept of a  {\it lift} mapping and then finally define evolving finite element spaces.

Throughout the remainder of the paper, we will denote by $k \ge 1$ the degree of the finite element spaces used both for the function space approximations and the approximation for the computational domain.

\subsection{Construction of the Initial Domain}
%Define the base mesh for both \Omega_1 and \Omega_2. (polyhedral boundary (\partial \Omega) (triangulation of \Partial \Omega is exact).  

%Dimensions = d, Number of elements = M_1 (\Omega_1), M_2 (\Omega_2)
 
\subsubsection*{Initial Mesh Construction/Assumption:}
\begin{enumerate}
    \newitem{$M1$} \label{M1} We first perform a partition into $d$-dimensional simplices corresponding to a polyhedral approximation $\widetilde{\Omega}^h_1$ of the interior domain $\overline{\Omega}_1(0)$, $\widetilde{\Omega}^h_1 = \cup_{j = 1}^{M_1} \widetilde{K}^j_1$, $\widetilde{\mathcal{J}}^h_1:= \{\widetilde{K}^j_1\}_{j = 1}^{M_1}$, where $\widetilde{K}^j_1$ are the simplicial elements of positive diameter, bounded by some $\widetilde{h}$, and $M_1$ is the number of elements. 
    \newitem{$M2$} \label{M2}
    The set $\widetilde{\Omega}^h_2 := \overline{\Omega \setminus \widetilde{\Omega}^h_1}$ is polyhedral and we construct a partition into $d$-dimensional simplices $\widetilde{\mathcal{J}}^h_2 := \{\widetilde{K}^j_2\}_{j=1}^{M_2}$ with maximum diameter $\widetilde{h}$. Let $\widetilde{\mathcal{J}}^h = \cup_{i = 1}^2 \widetilde{\mathcal{J}}^h_i$ assume that all partitions $\{\widetilde{\mathcal{J}}_1^h, \widetilde{\mathcal{J}}_2^h, \widetilde{\mathcal{J}}^h\}$ are admissible, shape regular and quasi-uniform in $\{\widetilde{\Omega}_1^h, \widetilde{\Omega}^h_2, \overline{\Omega}\}$ respectively, see \cite[Def.~5.1]{braess2007finite}.
    \newitem{$M3$} \label{M3} Each element $\widetilde{K}$ contains $d+1$ facets labelled $\{\widetilde{E}^j\}_{j = 1}^{d+1} \subset \widetilde{K}$. We refer to the set of all facets of all elements in $\widetilde{\mathcal{J}}^h$ by $\widetilde{\mathcal{J}}^h_\partial$. 
    \newitem{$M4$} \label{M4} For $\widetilde{E} \in \widetilde{\mathcal{J}}^h_\partial$, if there exists $\widetilde{K}_1 \in \widetilde{\mathcal{J}}_1^h$ and $\widetilde{K}_2 \in \widetilde{\mathcal{J}}^h_2$ such that $\widetilde{E} = \widetilde{K}_1 \cap  \widetilde{K}_2$, then we call $\widetilde{E}$ an \textit{interface facet} and label the collection of those facets $\widetilde{\mathcal{J}}_\Gamma^h$ and the union of interface facets $\widetilde{\Gamma}^h$. If for a given $\widetilde{E}$, there is only one element $\widetilde{K}\in \widetilde{\mathcal{J}}^h_2$ such that $\widetilde{E} \subset \widetilde{K}$, then such a facet is called a \textit{boundary facet}.
    \newitem{$M5$} \label{M5} We restrict the vertices of interface facets to be on $\Gamma(0)$, i.e, if $\widetilde{E}$ is an interface facet, and $\{\widetilde{a}_{\widetilde{E}}^j\}_{j = 1}^{d-1}$ are the vertices of $\widetilde{E}$, then $\{\widetilde{a}_{\widetilde{E}}^j\}_{j = 1}^{d-1} \subset \Gamma(0)$. Conversely, we will assume that if a facet has all its vertices on the interface, then it is an interface facet. See \cref{fig:topView} for an example.
    \newitem{$M6$} \label{M6} Let $\widehat{K}_{\mathrm{ref}}$ be the reference element of $\widetilde{\mathcal{J}}$ (i.e for all $\widetilde{K} \in \widetilde{\mathcal{J}}$, there exists an invertible affine map $F_{\widetilde{K}}$ such that $F_{\widetilde{K}}(\widehat{K}_{\mathrm{ref}}) = \widetilde{K}$). The reference element is then equipped with the standard $k^{th}$ ($k \in (0, \Theta)$) Lagrangian element triple $(\widehat{K}_{\mathrm{ref}}, \widehat{P}^k, \widehat{\Sigma}^k)$ (see \cite[Sec.~3.2]{MR2373954}) where $\widehat{P}^k$ is the set of $k^{th}$ order Lagrange polynomials and $\widehat{\Sigma}^k$ is the dual basis of $\widehat{P}^k$, which in this case takes the form $\widehat{\Sigma}= \{ \chi \to \chi(\widehat{\alpha}), \; \widehat{\alpha} \in N(\widehat{K}_{\mathrm{ref}})\}$, where $N(\widehat{K}_{\mathrm{ref}})$ is the set of Lagrangian nodes in $\widehat{K}_{\mathrm{ref}}$. Let $(\widetilde{K}, \widetilde{P}^k,\widetilde{\Sigma}^k)$ and $(\widetilde{K}', \widetilde{P}'^k, \widetilde{\Sigma}'^k)$ be two adjacent elements in $\mathcal{\widetilde{J}}$, the following assumption is made
    \begin{align*}
        \bigg(\bigcup_{\alpha \in N(\widetilde{K})} \widetilde{\alpha} \bigg) \cap \widetilde{K}' = \bigg(\bigcup_{\alpha' \in N(\widetilde{K}')} \widetilde{\alpha}' \bigg) \cap \widetilde{K},
    \end{align*}
    i.e the Lagrangian nodes are shared between two adjacent elements.
\end{enumerate}

Note that via construction, $\widetilde{E} \in \widetilde{\mathcal{J}}_\Gamma^h$, if and only if there exists an element $\widetilde{K}_1 \in \widetilde{\mathcal{J}}^h_1$ and $\widetilde{K}_2 \in\widetilde{\mathcal{J}}^h_2$ with $\widetilde{E} = \widetilde{K}_1 \cap \widetilde{K}_2$ and hence $\widetilde{\Gamma}^h_0 = \widetilde{\Omega}^h_1(0) \cap \widetilde{\Omega}^h_2(0)$. This construction defines Lagrangian triangulated bulk domains $(\widetilde{\Omega}^h_1, \widetilde{\Omega}^h_2, \widetilde{\Omega}^h)$, and $\widetilde{\Gamma}^h_0$ defines a triangulated hypersurface, see \cite[Def.~4.14 and 6.14]{EllRan21}.
\begin{table}
\begin{center}
\renewcommand{\arraystretch}{1.3}
\begin{tabular}{|l | l|}
\hline
$\Theta$ & Degree of additional geometric regularity assumed in \cref{lem:geom-regularity}. \\
$\left(\widehat{K}_{\mathrm{ref}}, \widehat{P}^k, \widehat{\Sigma}^k\right)$& Standard $k-th$ order Lagrangian reference element, with $k \in (0, \Theta)$.\\
$\widehat{\alpha}$ & Lagrange node of the reference element. \\
$i\in \{1,2\}$ & As a subscript, will always only refer to which of the domains\\
 & the quantity appertains.\\
$\widetilde{\Omega}_i^h, \widetilde{\Gamma}^h$& Initial Triangulation of the domains and interface(at $t = 0$).\\
$\widetilde{\mathcal{J}}_i^h, \widetilde{\mathcal{J}}^h, \widetilde{\mathcal{J}}^h_{\Gamma}$& Partitions of $\widetilde{\Omega}_i^h, \overline{\Omega}$ and $\widetilde{\Gamma}^h$ respectively.\\
$\widetilde{K},\widetilde{E}, \widetilde{\alpha}, \widetilde{a}$& Element/Facet/Lagrangian node/vertex appertaining to $\widetilde{\mathcal{J}}^h$.\\ 
$\mathbf{\Psi}^h$& Diffeomorphism map $\mathbf{\Psi}^h\colon \widetilde{\Omega}_i^h \to \overline{\Omega}_i$.\\
$\Pi_0(x)$& Minimal distance projection onto $\Gamma(0)$.\\
$\Omega^h_i(0), \Gamma^h(0)$& Triangulated bulk domains (hypersurface) approximating $\overline{\Omega}_i(0), \Gamma(0)$.\\
$\mathcal{J}^h_i, \mathcal{J}^h, \mathcal{J}^h_\Gamma$& Partition of isoparametric element of $\Omega^h_i(0), \overline{\Omega}, \Gamma^h(0)$.\\
\hline 
\end{tabular}
\caption{List of symbols}
\label{list_of_symbols}
\end{center}
\end{table}
After the initial triangulation, we define the isoparametric version using the same method as \cite[Sec.~8.5]{EllRan21} which we detail in the following. Let $\widehat{K}_{\mathrm{ref}}$ be the reference element of the partition $\widetilde{\mathcal{J}}^h$, with reference map $F_{\widetilde{K}}:\widehat{K}_{\mathrm{ref}} \to \widetilde{K}$.

 Let $\{\widetilde{a}_{\widetilde{K}}^j\}_{j  =1}^{d+1}$ be the vertices of an element $\widetilde{K} \in \widetilde{\mathcal{J}}^h$. If two or more of the vertices are on the interface $\Gamma(0)$, then the element is referred to as an \textit{interface element}. Let $\widetilde{\mathcal{F}}$ be the set of all interface elements and define the following function $\mathbf{\Psi}^h: \Omega \to \Omega$ element-wise as follows. If $\widetilde{K} \notin \widetilde{\mathcal{F}}$, then $\mathbf{\Psi}^h(x) = x$ for $x \in \widetilde{K}$. If instead $\widetilde{K} \in \widetilde{\mathcal{F}}$, then expand $x \in \widetilde{K}$ into barycentric coordinates:
\begin{align*}
    x = \sum_{j = 1}^{d+1} \mu_j(x) \widetilde{a}_{\widetilde{K}}^j.
\end{align*}
Let $L_K$ be the number of vertices in $\widetilde{K}$ that lie on $\Gamma(0)$ ($L_K \geq 2$ by assumptions) and assume that the vertices are ordered so that the first $L_K$ lie on $\Gamma(0)$. Let:
\begin{align*}
    \widetilde{\mu}_K(x):= \sum_{j = 1}^{L_K} \mu_j(x), \quad \sigma_{\widetilde{K}} := \{ x \in \widetilde{K}, \; \widetilde{\mu}_K(x) = 0 \}.
\end{align*}
From the properties of barycentric coordinates, $\widetilde{\mu}_K$ can be seen as the distance from the discrete interface, with $\widetilde{\mu}_K(x) = 1$ when $x$ is on a facet between vertices on the interface, and $\widetilde{\mu}_K(x) = 0$ when $x$ is on the facet spanned by non-interface vertices.

Let
\begin{align}\label{eq:norm-proj-y}
    y(x) = \sum_{j = 1}^{L_K} \frac{\mu_j(x)}{\widetilde{\mu}_K(x)} \widetilde{a}_{\widetilde{K}}^j.
\end{align}
Note that $y(x) \in \widetilde{K}$ since $0 \leq \mu_j(x) \leq \widetilde{\mu}_K(x)$. Hence define:
\begin{align}\label{eq:push-forward}
    \mathbf{\Psi}^h|_{\widetilde{K}}(x) := \begin{cases} 
 x + (\widetilde{\mu}_K(x))^{k+2}(\Pi_0(y(x)) - y(x)) \quad &\text{if} \, x \notin \sigma_{\widetilde{K}}, \\
     x \quad &\text{otherwise},
   \end{cases}
\end{align}
where $\Pi_0$ is the nearest point projection on $\Gamma(0)$, introduced in \cref{eq:tubular-decom}.
Here, and also in \cref{eq:lift-define}, we raise the barycentric coordinate to the power $k+2$. Whilst simply not raising to any power (i.e., raising to the power 1), would suffice for ensuring that both maps as bijections, raising to a higher power allows us to ensure higher order smoothness of each map.

We summarise the properties of this map in the following theorem. For the definition of \textit{triangulated bulk domain} and \textit{k-bulk finite element}, see \cite[Def.~4.14 and 4.5]{EllRan21}. 
We denote by $\tilde{I}^h$ Lagrangian interpolation into the space of polynomials of degree $k$ over a single element $\tilde{K}$.

\begin{theo}[\cite{EllRan21}, Lem.~4.8 and 8.8]\label{theo:main-theo}
For $\widetilde{h}$ small enough, the map $ \mathbf{\Psi}^h|_{\widetilde{K}} \in C^{k+1}(\widetilde{K}; \mathbb{R}^d)$and is invertible for each $\widetilde{K} \in \widetilde{\mathcal{J}}^h$ and $ \mathbf{\Psi}^h: \widetilde{\Gamma}^h_0 \to \Gamma(0)$. Define the following:
\begin{align*}
    F_{K} &:= [\widetilde{I}^h  \mathbf{\Psi}^h](F_{\widetilde{K}}), \\
    K&:= F_K(\widehat{K}_{\mathrm{ref}}), \\
    P^k&:= \{ \widehat{\chi}_k \circ F^{-1}_K: \widehat{\chi}_k \in \widehat{P}^k \}, \\
    \Sigma^k&:= \{\chi \mapsto \widehat{\sigma}(\chi \circ F_K): \, \widehat{\sigma} \in \widehat{\Sigma}^k\},
\end{align*}
then the triplet $(K, P^k, \Sigma^k)$ with reference map $F_{K}$ defines a $k$-bulk finite element triplet (\cite[Def.~4.5]{EllRan21}). Let $\mathcal{J}^h_i = \{[\widetilde{I}^h  \mathbf{\Psi}^h](\widetilde{K}_i), \; \widetilde{K}_i \in \widetilde{\mathcal{J}}^h_i\}$, $\mathcal{J}^h_\Gamma = \{ [\widetilde{I}^h\mathbf{\Psi}^h]_1(E), \, E \in \widetilde{\mathcal{J}}^h_\Gamma\}$ (here $[\widetilde{I}^h\mathbf{\Psi}^h]_1$ refers to taking the interpolation with the adjacent element in $\widetilde{\mathcal{J}}_1$),  then $\{\mathcal{J}^h _1, \mathcal{J}^h _2, \mathcal{J}_\Gamma^h \}$ are conforming admissible sub-divisions. Furthermore, let:
\begin{align*}
    \Omega^h_i(0) := \bigcup_{K_i \in \mathcal{J}^h_i} K_i, \, \Omega^h := \bigcup_{K \in \mathcal{J}^h} K, \, \Gamma^h(0) := \bigcup_{E\in \mathcal{J}^h_\Gamma} E,
\end{align*}
then $(\Omega^h_1(0),\Omega^h_2(0))$ define triangulated bulk domains approximating $(\overline{\Omega}_1(0), \overline{\Omega}_2(0))$, $\Gamma^h(0)$ a triangulated hypersurface approximating $\Gamma(0)$.
\end{theo}
Now since we are dealing with an interface problem, we require additionally to check if $\mathcal{J}^h = \mathcal{J}^h_1 \cup \mathcal{J}^h_2$ forms a conforming admissible sub-division of the whole domain $\overline{\Omega}$.
\begin{lemma}\label{lemma:bulkdef} $\mathcal{J}^h$ forms a conforming admissible sub-division of the whole domain $\overline{\Omega}$, moreover, interface facets are mapped to their isoparametric equivalent in such a way that:
\begin{align*}
    \Gamma^h(0) = \bigcup_{E\in \mathcal{J}^h_\Gamma} E = \Omega_1^h \cap \Omega_2^h.
\end{align*}
\end{lemma}
\begin{proof}
Since $\mathcal{J}^h$ is the union of two admissible conforming subdivision, it only remains to check that if we are given two elements $K_i \in \mathcal{J}^h_i$, then  $K_1^\circ \cap K_2^\circ = \varnothing$. It suffices to show the invertibility of the map $\tilde{I}^h \mathbf{\Psi}^h$ on $\overline{\Omega}$. For an interface facet $\widetilde{E}$ with two adjacent element $\widetilde{K}_i$ we require continuity across $E$: $\mathbf{\Psi}^h|_{\widetilde{K}_1}(\widetilde{E}) = \mathbf{\Psi}^h|_{\widetilde{K}_2}(\widetilde{E})$. By construction of the mesh, $x \in \widetilde{E}$, $\mu_{K_1}(x) = \mu_{K_2}(x) = 1$. This implies $y(x) = x$ in \cref{eq:norm-proj-y} and hence both maps $\mathbf{\Psi}^h|_{\widetilde{K}_i}(x) = \Pi_0(x)$ from \cref{eq:push-forward}. Any Lagrangian node $\widetilde{\alpha}_i$ on $\widetilde{E}$ will be mapped by both maps to $\alpha_i := \Pi_0(\widetilde{\alpha}_i)$. Since each interface facet contains the exact amount of nodes to uniquely define a polynomial on the facet, which must equal the restriction on the interface element of the Lagrangian polynomial on the full element (see \cite[Rem.~5.4]{braess2007finite}), hence for $x \in \widetilde{E}$:
\begin{align*}
    [\tilde{I}^h \mathbf{\Psi}^h]|_{\widetilde{K}_1}(x) 
    &= \sum_{\{j : \chi^j \in P^k \}} \mathbf{\Psi}^h(\widetilde{\alpha}_{\widetilde{K}_1}^j) \chi^j(x)
    = \sum_{\{j : \chi^j \in P^k \}} \Pi_0(\widetilde{\alpha}^j_{\widetilde{K}_1}) \chi^j(x),\\
    &= \sum_{\{j : \chi^j \in P^k \}} \Pi_0(\widetilde{\alpha}^j_{\widetilde{K}_2}) \chi^j(x)
    = \sum_{\{j : \chi^j \in P^k \}} \mathbf{\Psi}^h(\widetilde{\alpha}_{\widetilde{K}_2}^j) \chi^j(x)
    =  [\tilde{I}^h \mathbf{\Psi}^h]|_{\widetilde{K}_2}(x).
\end{align*}
For $\widetilde{E} \in \widetilde{\mathcal{J}}^h_\Gamma$, let $\widetilde{K}_i \in \widetilde{\mathcal{J}}^h_i$ be the adjacent elements to $E$ and $K_i = I^h \mathbf{\Psi}^h(\widetilde{K}_i)$. Since the map $I^h \mathbf{\Psi}^h$ is invertible onto its image for each $\widetilde{K} \in \widetilde{\mathcal{J}}^h$ and is continuous across the intersection $\widetilde{E}$, it holds that $I^h \mathbf{\Psi}^h$ is invertible on $\widetilde{K}_1 \cup \widetilde{K}_2$, since both elements are closed, hence:
\begin{align*}
    E := [I^h \mathbf{\Psi}^h]|_{\widetilde{K}_1}(\widetilde{E}) = I^h \mathbf{\Psi}^h(\widetilde{K}_1 \cap \widetilde{K}_2) = I^h \mathbf{\Psi}^h(\widetilde{K}_1) \cap I^h \mathbf{\Psi}^h(\widetilde{K}_2) = K_1 \cap K_2.
\end{align*}
Therefore the image of an interface facet remains an interface facet. Moreover, this shows that $K^\circ_1 \cap K^\circ_2 = \varnothing$ for any $K_i \in \mathcal{J}^h_i$ and hence $\mathcal{J}^h$ is a conforming admissible sub-division.
\end{proof}

\begin{figure}[th]
\fontsize{0.6cm}{0}\resizebox{250mm}{!}{\hbox{\hspace{3em} 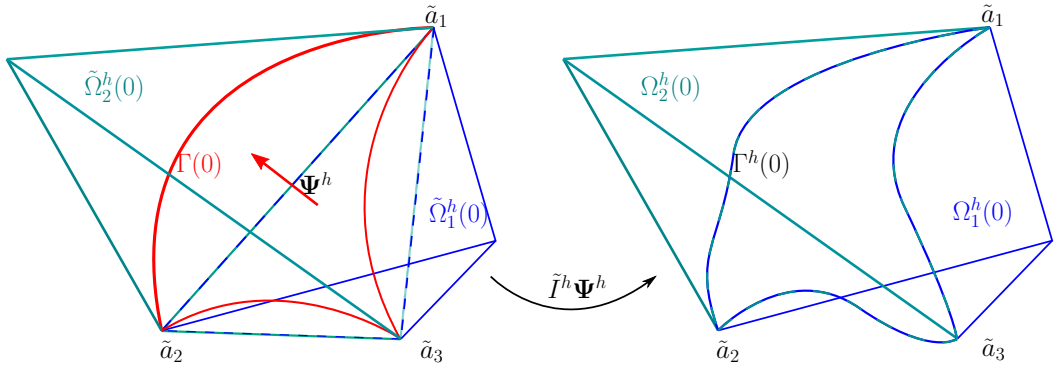}}
\caption{\textit{The intersection of two interface elements (teal and blue respectively) of different domains. The shared interface facet is then pushed by the map $\mathbf{\Psi}^h$ to become a piece of $\Gamma(0)$. The map $\widetilde{I}^h \mathbf{\Psi}^h$ maps the original mesh to an isoparametric mesh approximating the interface.}}
\label{fig:intermesh}
\end{figure}

\Cref{fig:intermesh} shows how the map $\tilde{I}^h \mathbf{\Psi}^h$ deforms the original mesh. We are initially given two tetrahedral elements of the initial meshes, one in $\widetilde{\Omega}^h_1$ and one in $\widetilde{\Omega}^h_2$, intersecting on an interface element. Applying the map $\tilde{I}^h \mathbf{\Psi}^h$ to this yields isoparametric elements whose intersection is the image of the interface element under $\tilde{I}^h \mathbf{\Psi}^h$.

Let $\alpha_K^j|_{j = 1}^{N(k)}$ be the Lagrangian nodes on an element $K \in \mathcal{J}^h$, which by construction are defined as $\alpha_K^j = F_K(\widehat{\alpha}_{\text{ref}}^j)$, the corresponding interpolation operator for a given function $\eta \in C(\overline{\Omega}; \mathbb{R})$ is given by:
\begin{align*}
    [I^h \eta]|_{K} = \sum_{\{ j : \chi^j \in P^k \}} \eta(\alpha_K^j) \chi^j.
\end{align*}
\subsection{Time Dependent Mesh}\label{section:time}
Define the flow:
\begin{align}\label{eq:discrete-flow}
    \mathbf{\Phi}^h_i(t; \cdot)|_{K_i(0)} :=
    I^h_{K_i(0)} [ \mathbf{\Phi}_i(t; \mathbf{\Psi}^h \circ (\widetilde{I}^h \mathbf{\Psi}^h)^{-1}(\cdot))].
\end{align} 
We denote by $\mathbf{\Phi}^h_i(-t, \cdot)|_{K_i(t)} \colon K_i(t) \to K_i(0)$ the space-only inverse of $\mathbf{\Phi}^h_i(t; \cdot)|_{K_i(0)}$.

\begin{remark}
\normalfont
The flow is defined this way such that it evolves the parametric meshes $\Omega_i^h(0)$. Indeed, decomposing $ \mathbf{\Phi}^h_i(t; \cdot)$ into its components:
\begin{align*}
    (\widetilde{I}^h \mathbf{\Psi}^h)^{-1}: \Omega^h_i(0) \to \widetilde{\Omega}^h_i, \; \mathbf{\Psi}^h: \widetilde{\Omega}^h_i \to \overline{\Omega}_i(0), \, \mathbf{\Phi}_i(t; \cdot):\overline{\Omega}_i(0) \to  \overline{\Omega}_i(t).
\end{align*}
Hence, the flow $\mathbf{\Phi}^h_i(t; \cdot)$ is a polynomial function approximating the evolution of the domains. Moreover, this does define a proper flow map, as, at $t = 0$:
\begin{align*}
    \mathbf{\Phi}^h_i(0; x)|_{K_i(0)} =  I^h_{K_i(0)} [ \mathbf{\Psi}^h \circ (\widetilde{I}^h \mathbf{\Psi}^h)^{-1}(x)] &= \sum_{\chi^j \in P^k} \mathbf{\Psi}^h \circ [(\widetilde{I}^h \mathbf{\Psi}^h)^{-1}(\alpha^j_{K_i})] \chi^j(x), \\
    &= \sum_{\chi^j \in P^k}  \mathbf{\Psi}^h (\widetilde{\alpha}^j_{K_i}) \chi^j(x) = \sum_{\chi^j \in P^k}  \alpha^j_{K_i} \chi^j(x)  = x.
\end{align*}
The composition property for $\mathbf{\Phi}^h_i(t +s;\cdot) = \mathbf{\Phi}^h_i(t; \mathbf{\Phi}^h_i(s; \cdot))$ will be shown following the next lemma.
\end{remark}
%We can add a note that later, we will define the lift $\Lambda^h$ which will satisfy $\Lambda^h(0, \cdot) = \mathbf{\Psi}^h \circ (\widetilde{I}^h \mathbf{\Psi}^h)^{-1}(\cdot)$. I've also checked and the previous proofs assumed this corrected version.

For $\widetilde{h}$ small enough, this map is an invertible diffeomorphism on each element $K_0 \in \mathcal{J}^h$. As before, we summarise the construction in the following lemma:
\begin{lemma}For $\widetilde{h}$ small enough, map $\mathbf{\Phi}^h_i(t; \cdot)|_{\widetilde{K}} \in C^{k+1}(\widetilde{K}; \mathbb{R}^d)$ and is invertible onto its image for each $\widetilde{K} \in \widetilde{\mathcal{J}}^h$. Moreover, define the following:
\begin{align*}
    F_{K(t)}(\cdot) &:=  \mathbf{\Phi}^h_i(t; F_{\widetilde{K}}(\cdot)), \\
    K(t)&:= F_{K(t)}(\widetilde{K}), \\
    P^k(t)&:= \{ \chi_k \circ F^{-1}_{K(t)}: \chi_k \in \widehat{P}^k \}, \\
    \Sigma^k(t)&:= \{\chi \mapsto \sigma(\chi \circ F_{K(t)}): \, \sigma \in \widehat\Sigma^k\},
\end{align*}
then the triplet $({K(t)}, {P^k(t)}, {\Sigma^k(t)})$ with reference map $ F_{K(t)}$ defines a bulk evolving finite element triplet. Let $\mathcal{J}^h_i(t) = \{\mathbf{\Phi}^h_i(t;\widetilde{K}_i), \; \widetilde{K}_i \in\mathcal{J}^h_i\}$, $\mathcal{J}^h_\Gamma(t) = \{\mathbf{\Phi}^h_1(t; E),\, E \in \mathcal{J}^h_\Gamma\}$ and $\mathcal{J}^h(t) = \mathcal{J}^h_1(t) \cup \mathcal{J}^h_2(t)$, then $\{\mathcal{J}^h _1(t), \mathcal{J}^h_2(t), \mathcal{J}^h(t) \}$ are evolving conforming admissible sub-divisions (see \cite[Def.~4.32]{EllRan21}). Furthermore, let:
\begin{align*}
    \Omega^h_i(t) = \bigcup_{K_i(t) \in \mathcal{J}^h_i(t)} K_i(t), \, \Omega^h(t) = \bigcup_{K \in \mathcal{J}^h} K(t), \, \Gamma^h(t) = \bigcup_{E \in \mathcal{J}^h_\Gamma} E= \Omega^h_1(t) \cap \Omega^h_2(t),
\end{align*}
then $(\Omega^h_1(t),\Omega^h_2(t))$ define triangulated bulk domains approximating $(\overline{\Omega}_1(t), \overline{\Omega}_2(t))$, $\Gamma^h(t)$ is a triangulated hypersurface approximating $\Gamma^h(t)$, and $\Omega^h(t)$ defines a triangulated bulk domain that is an exact partition of $\overline{\Omega}$
\end{lemma}
\begin{proof}
The proof follows the same way as \cref{lemma:bulkdef}.
\end{proof}

For each $K(t) \in \mathcal{J}^h(t)$, let $h_{K(t)}$ be the diameter of the flat simplex whose vertices match $K(t)$.
We define $h := \max_{t \in I} \max_{K(t) \in \mathcal{J}^h(t)} \mathrm{diam}(K(t))$ to be the maximum mesh diameter, where $\mathrm{diam}(K(t))$ is the diameter of the affine element whose vertices match $K(t)$ (see \cite[Lem.~4.9]{EllRan21}).

\begin{figure}[th]
\fontsize{0.8cm}{0}\resizebox{140mm}{!}{\hbox{\hspace{10em} %% Creator: Inkscape 1.2 (56b05e47e7, 2022-06-09, custom), www.inkscape.org
%% Accompanies image file 'flowmap.pdf' (pdf, eps, ps)
\begingroup%
  \makeatletter%
  \providecommand\color[2][]{%
    \errmessage{(Inkscape) Color is used for the text in Inkscape, but the package 'color.sty' is not loaded}%
    \renewcommand\color[2][]{}%
  }%
  \providecommand\transparent[1]{%
    \errmessage{(Inkscape) Transparency is used (non-zero) for the text in Inkscape, but the package 'transparent.sty' is not loaded}%
    \renewcommand\transparent[1]{}%
  }%
  \providecommand\rotatebox[2]{#2}%
  \newcommand*\fsize{\dimexpr\f@size pt\relax}%
  \newcommand*\lineheight[1]{\fontsize{\fsize}{#1\fsize}\selectfont}%
  \ifx\svgwidth\undefined%
    \setlength{\unitlength}{695.04927027bp}%
    \ifx\svgscale\undefined%
      \relax%
    \else%
      \setlength{\unitlength}{\unitlength * \real{\svgscale}}%
    \fi%
  \else%
    \setlength{\unitlength}{\svgwidth}%
  \fi%
  \global\let\svgwidth\undefined%
  \global\let\svgscale\undefined%
  \makeatother%
  \begin{picture}(1,0.33444047)%
    \lineheight{1}%
    \setlength\tabcolsep{0pt}%
    \put(0,0){\includegraphics[width=\unitlength,page=1]{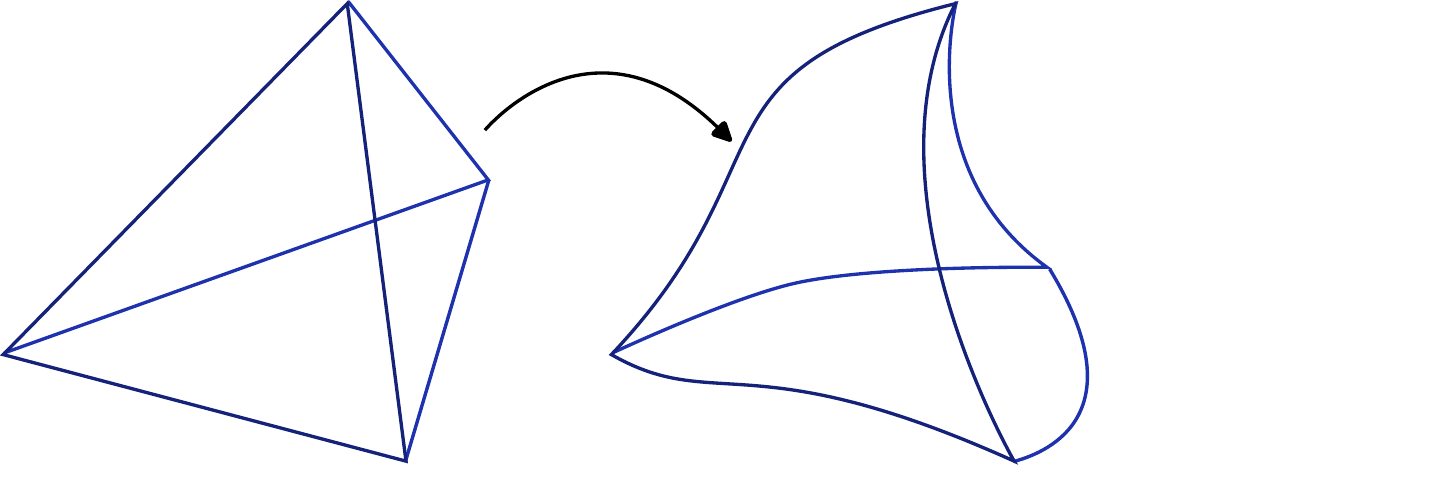}}%
    \put(0.09905352,0.16581211){\color[rgb]{0,0,0}\makebox(0,0)[lt]{\lineheight{1.25}\smash{\begin{tabular}[t]{l}$K(0)$\end{tabular}}}}%
    \put(0.49419041,0.16003603){\color[rgb]{0,0,0}\makebox(0,0)[lt]{\lineheight{1.25}\smash{\begin{tabular}[t]{l}$K(t)$\end{tabular}}}}%
    \put(0,0){\includegraphics[width=\unitlength,page=2]{flowmap.pdf}}%
    \put(0.39252476,0.29886067){\color[rgb]{0,0,0}\makebox(0,0)[lt]{\lineheight{1.25}\smash{\begin{tabular}[t]{l}$\mathbf{\Phi}^h(t; \cdot)$\end{tabular}}}}%
  \end{picture}%
\endgroup%
}}
\caption{\textit{Example of the temporal deformation of an interior element in three space dimensions.}}
\label{fig:flowmesh}
\end{figure}

\begin{remark}
\normalfont
This allows us to move the Lagrangian nodes via $\alpha_{K_i(t)}^j = \mathbf{\Phi}_i^h(t; \alpha^j_{K_i})$. The Lagrangian interpolation operator, $I^h|_{K(t)}$, is then defined in the canonical way. Moreover for $x \in K_i(0)$:
\begin{align*}
    \partial_t  \mathbf{\Phi}_i^h(t; x) &= \sum_{\, \substack{\{j:\widetilde{\chi}^j \in \widetilde{P}^k\}}} \partial_t  \mathbf{\Phi}_i(t; \mathbf{\Psi}^h(\alpha^j_{\widetilde{K}_i})) \widetilde{\chi}^j(x)= \sum_{\substack{\{j:\,\widetilde{\chi}^j \in \widetilde{P}^k\}}}\mathbf{w}(t;  \mathbf{\Phi}_i\circ\mathbf{\Psi}^h(\alpha^j_{\widetilde{K}_i}))\widetilde{\chi}^j(x) \\
    &= \sum_{\substack{ \{j:\,\chi^j \in P^k(t)\}}}\mathbf{w}(t; \alpha^j_{K_i(t)})\chi^j(t;\mathbf{\Phi}_i^h(t; x)) =: \mathbf{w}^h(t;  \mathbf{\Phi}_i^h(t; x)),
\end{align*} 
where one sees that $\mathbf{w}^h$ is the interpolated velocity with respect to the moving nodes:
\begin{align}\label{eq:discrete-velocity-2}
    \mathbf{w}^h(t;\cdot)|_{K_i(t)} = I^h|_{K_i(t)}[\mathbf{w}(t;\cdot)].
\end{align}
Hence, element-wise, the discrete flow satisfies ODE:
\begin{align}\label{eq:discrete-flow-ode}
   \frac{d}{dt} \mathbf{\Phi}^h_i(t;x) &= \mathbf{w}^h(t; \mathbf{\Phi}_i^h(t;x)), \quad x \in K_i(0),\\
    \mathbf{\Phi}^h_i(0;x) &= x \nonumber,
\end{align}
and therefore satisfies the composition property $\mathbf{\Phi}^h_i(t +s; \cdot) = \mathbf{\Phi}^h_i(t;\mathbf{\Phi}^h_i(s;\cdot))$, see \cite{hartman1982ordinary}.
\end{remark}
It will be assumed that the mesh remains \textit{uniformly quasi-uniform} in time, see \cite[Def.~4.35]{EllRan21} as the discrete flow $\mathbf{\Phi}^h$ can deform the mesh significantly. An example of the temporal deformation of an evolving element is shown in \cref{fig:flowmesh}. Despite interior elements of the initial partition being linear, since the velocity used to displace the elements is a polynomial interpolant of the velocity, the resulting element might not remain linear and can be deformed.
An alternative construction, for which interior elements remain affine, is given in \cite{Li2022}.

The discrete spaces are then defined as:
\begin{align*}
    H^h(t)& :=   L^2(\Omega^h_1(t))  \times  L^2(\Omega^h_2(t)), \\
V^h(t)&:= \{ \eta^h \in W^{1,2}(\Omega^h_1(t)) \times W^{1,2}(\Omega^h_2(t)), \eta_1^h - \eta_2^h|_{\Gamma^h(t)} = 0 ,\; \text{and} \; \eta^h_2|_{\partial \Omega^h} = 0\},
\end{align*}
equipped with the norms:
\begin{align*}
  \norm{\cdot}^2_{V^h(t)}:= \sum_{i = 1}^2 \norm{\cdot}^2_{W^{1,2}(\mathcal{J}^h_i(t))}, \quad
  \norm{\cdot}^2_{H^h(t)}:= \sum_{i = 1}^2 \norm{\cdot}^2_{L^2(\Omega^h_i(t))}.
\end{align*}

Define the map $\phi^h_t:H^h(0) \to H^h(t)$ element-wise as:
\begin{align*}
    (\phi^h_t v^h)(x_1, x_2) := (v^h_1(\mathbf{\Phi}^h_1(-t;x_1)),
    v^h_2(\mathbf{\Phi}^h_2(-t;x_2)))
    \qquad \text{for} \quad
    (x_1, x_2) \in K_1(t) \times K_2(t) \subset \Omega^h_1(t) \times \Omega^h_2(t).
\end{align*}
That is
\begin{align*}
    (\phi^h_t v^h)(x_1, x_2) := (v^h_1(y_1), v^h_2(y_2))
    \qquad \text{for} \quad
    y_i \in K_i(0), y_i = \mathbf{\Phi}^h_i(-t; x_i), x_i \in K_i(t).
\end{align*}

\begin{lemma}
 $(V^h(t), \phi^h_t)_{t \in I}$ and $(H^h(t), \phi^h_t)_{t \in I},$ are compatible pairs.
\end{lemma}
\begin{proof}
This follows by the regularity of the map $\mathbf{\Phi}^h$ and $\mathcal{J}^h(t)$, see \cite[Lem.~4.36]{EllRan21}.
\end{proof}
Hence the moving spaces $L^2_{V^h}$ and $L^2_{H^h}$ are well defined. Denote the discrete material derivative by:
\begin{align}\label{eq:discrete-material}
    \partial^h_t \eta := \phi^h_{t} \frac{\partial}{\partial t} \phi_{-t}^h \eta,
\end{align}
for $\eta \in C^1_{H^h}$. The bilinear form $\lambda^h(t; \cdot, \cdot)$ of \cref{def:lambda} associated with this material derivative is:
\begin{align*}
    \lambda^h(t;\eta^h,v^h) = (\nabla \cdot \mathbf{w}^h \eta^h, v^h)_{H^h(t)},
\end{align*}
where $\mathbf{w}^h$ is the previously defined discrete velocity from \cref{section:time} (see \cite[Lem.~8.10]{EllRan21} for derivation). This allows us to define, just as before, the discrete space:
\begin{align*}
    W(V^h, V^h) = \left\{ v^h \in L^2_{V^h}, \; \partial_t^h v_h \in L^2_{V^h}\right\}.
\end{align*}

\subsection{The lift}
The last mesh related concept needed is the \textit{lift map} (see \cite[Sec.~8.6]{EllRan21}).
The lift is used to bring functions defined on the computational domain to the continuous domain of the smooth problem. Our error estimate compares the exact solution to the composition of the discrete solution with the lift.
For interior elements with one or no vertices on the interface, the lift map is the identity.
Otherwise, for interface elements, the lift uses the normal projection operator and barycentric coordinates on the reference element in a one-element-wide region on each side of the interface. 

Fix $t \in I$, if $K(t) \in \mathcal{J}^h(t)$ is an interior element, then define the lift $\Lambda^h(t;\cdot)$ as:
\begin{align*}
   \Lambda^h(t;x) = x, \; \text{for} \; x \in K(t).
\end{align*}
If instead $K(t)$ is an interface element, we first pull-back the reference map to $\widehat{x} \in \widehat{K}_{\text{ref}}$ such that $x = F_{K(t)}(\widehat{x})$, then decomposing $\widehat{x}$ into barycentric coordinates with respect to the vertices $\widehat{a}_{K_{\text{ref}}}^j$ of $\widehat{K}_{\text{ref}}$, we have:
\begin{align*}
    \widehat{x} = \sum_{j = 1}^{d+1} \mu_j(\widehat{x}) \widehat{a}^j_{K_{\text{ref}}},
\end{align*}
and once again, let $L_K$ be the number of vertices on the interface and assume the vertices are ordered so that $\widehat{a}^j_{K_{\text{ref}}}$, $j=1, \ldots, L_K$, get mapped on to $\Gamma(t)$, then we introduce the interface distance and the singular set analogously:
\begin{align*}
    \widetilde{\mu}(\widehat{x}) = \sum_{j = 1}^{L_K}\mu_j (\widehat{x}), \; \sigma = \{\widehat{x} \in \widehat{K}_{\text{ref}}|\; \widetilde{\mu}(\widehat{x}) = 0\}.
\end{align*}
The projection is now defined on the reference element:
\begin{align*}
    \widehat{y}(\widehat{x}) = \sum_{j = 1}^{L_K} \frac{\mu_i(\widehat{x})}{\widetilde{\mu}(\widehat{x})}, \; y(t;x):= F_{K(t)}( \widehat{y}(\widehat{x})).
\end{align*}
Hence the lift operator can now be defined on interface elements as:
\begin{align}
  \label{eq:lift-define}
    \Lambda^h(t; x)|_{K(t)} = \begin{cases} 
 x + (\widetilde{\mu}(\widehat{x}))^{k+2}(\Pi_t(y(t;x)) - y(t;x)) \quad &\text{if} \, \widehat{x} \notin \sigma,\\
     x \quad &\text{otherwise}.
   \end{cases}
\end{align}

\begin{figure}[ht]
\fontsize{40}{10}\selectfont
\centering{
\resizebox{150mm}{!}{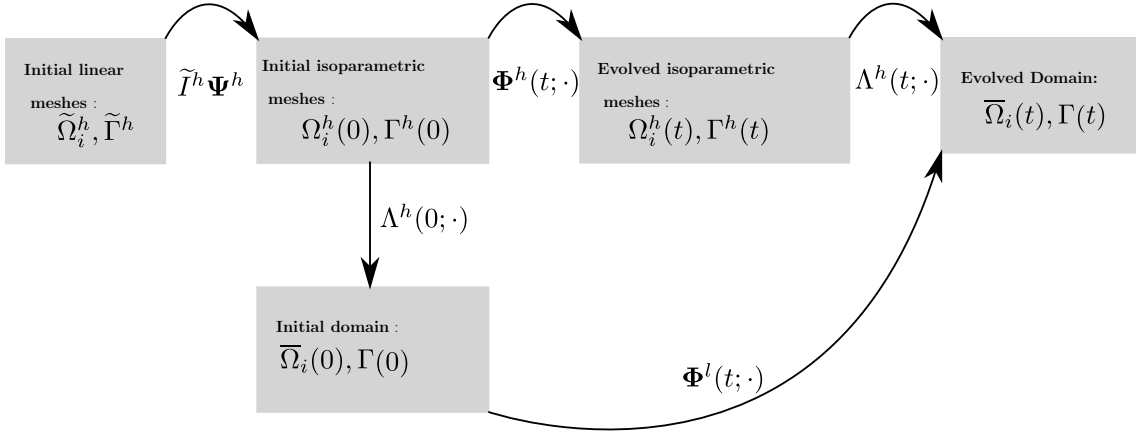}
\caption{Schematic of the setup used. $\Lambda^h(t;x)$ might be needed, depending on the problem, to define the discrete data. However once the discrete problem is known, only the knowledge of $\Omega^h_i(0)$, $\Gamma^h(0)$ and $\mathbf{\Phi}^h(t;\cdot)$ are needed to calculate the discrete solution $U^h(t;\cdot)$. $\mathbf{\Phi^l}$ is only needed in the analysis of theoretical error estimates.
We note that our constructions imply that this diagram is commutative.}
\label{fig:chart-flows}
}
\end{figure}

%\textcolor{magenta}{[The proper definitions for here are given later - they should be moved up!]}
%For a function $\zeta \in C^1(\cup_{t \in I} \{t\} \times K_i(t); \mathbb{R}^d)$, where $K_i(t) \in \mathcal{J}^h_i(t)$, $K_i(t) = \mathbf{\Phi}^h_i(t;\cdot)|_{\widetilde{K} }$ for all $t \in I$, for some $\widetilde{K} \in \widetilde{\mathcal{J}}^h$, define for the moment $(\phi_{-t}^h \zeta(t))(x) = \zeta(t;\mathbf{\Phi}_i(t; x))$ where $x \in K_i(0)$, and similarly $(\phi_{t}^h \eta(t))(x) = \eta(t;\mathbf{\Phi}_i(-t; x))$ for $\eta \in  C^1(I \times K_i(0); \mathbb{R}^d)$ and $ x \in K_i(t) \in \mathcal{J}^h_i(t)$. Let
%\begin{align*}
%    \partial^h_t \zeta:= \left.\left[\frac{d}{dt}(\zeta(t; \mathbf{\Phi}_i^h(t;y))\right]\right|_{y = \mathbf{\Phi}_i(-t;x)} = \phi^h_t \frac{d}{dt} \phi^h_{-t}\eta.
%\end{align*}
% \textcolor{magenta}{[end]}
Then, computing component wise, we see:
\begin{align*}
    \partial_t^h \Lambda^h(t; x)|_{K(t)} = \begin{cases} 
 \mathbf{w}^h(t;x) +\partial_t^h [ (\widetilde{\mu}(\widehat{x}))^{k+2} (\Pi_t(y(t;x)) - y(t;x)) ]\quad &\text{if} \, \widehat{x} \notin \sigma,\\
     \mathbf{w}^h(t;x) \quad &\text{otherwise}.
   \end{cases}
\end{align*}
Let $z = \mathbf{\Phi}_i^h(-t; x)$ (depending on whether $x \in \Omega_i^h(t)$, we also use the shorthand $\mathbf{\Phi}_{-t}(x) = \mathbf{\Phi}_i^h(-t; x)$ and $\mathbf{\Phi}_{t}(z) = \mathbf{\Phi}_i^h(t; z)$, depending on whether $z \in \Omega^h_i(0)$):
\begin{align*}
    \partial_t^h [ (\widetilde{\mu}(\widehat{x}))^{k+2} (\Pi_t(y(t;x)) - y(t;x)) ] &= \phi_t^h \frac{d}{dt}  [ (\widetilde{\mu}(F^{-1}_{\hat{K}} z))^{k+2} (\Pi_t(y(t;\mathbf{\Phi}^h_t(z)) - y(t;\mathbf{\Phi}^h_t(z)) ], \\
    &= \widetilde{\mu}(\hat{x})^{k+2} \phi_t^h \frac{d}{dt}[\Pi_t(y(t;\mathbf{\Phi}^h_t(z))) - y(t;\mathbf{\Phi}^h_t(z)) ],\\
    &=   \widetilde{\mu}(\hat{x})^{k+2}\left[\partial_t\Pi_t(y(t;x))  +\partial^h_t y(t;x)\cdot[\nabla \Pi_t (y(t;x))] -\partial^h_t y(t;x)\right].
\end{align*}
The formula for $\partial^h_t y(t;x)$ can also be explicitly found:
\begin{align*}
    \partial^h_t y(t;x) = \phi^h_t \frac{d}{dt} F_{K(t)}( \widetilde{y}(t; F^{-1}_{\hat{K}} z))= \phi^h_t\mathbf{w}^h(t; F_{K(t)}\widetilde{y}(t; F^{-1}_{\hat{K}} z)) = \mathbf{w}^h(t; y(t;x)),
\end{align*}
by use of \cref{eq:discrete-flow} and the definition of $F_{K(t)}$. Hence:
\begin{align}\label{eq:lift-deriv}
    & \partial_t^h \Lambda^h(t; x)|_{K(t)} \\
    &= \begin{cases} 
 \mathbf{w}^h(t;x) + \widetilde{\mu}(\hat{x})^{k+2}\left[\partial_t\Pi_t(y)  +\mathbf{w}^h(t; y)\cdot[\nabla \Pi_t (y)] -\mathbf{w}^h(t; y)\right]\quad &\text{if} \, x \notin \sigma \nonumber\\
     \mathbf{w}^h(t;x) \quad &\text{otherwise}
   \end{cases}\\
   &= \begin{cases} 
 \mathbf{w}^h(t;x) -\widetilde{\mu}(\hat{x})^{k+2}\left[ \left(\mathbf{w}^h(t;y) - \mathbf{w}(t;\Pi_t(y))\right)\cdot \nu_\Gamma(\Pi_t(y)) \nu_\Gamma(\Pi_t(y)) + d_\Gamma(t;y)T(y)\right]\quad &\text{if} \, x \notin \sigma,\\
     \mathbf{w}^h(t;x) \quad &\text{otherwise},
      \end{cases}
      \nonumber
\end{align}
with $T(x): = \partial_t[\nu_\Gamma(\Pi_t(x))] + \mathbf{w}^h(t;x) \cdot\nabla [\nu_\Gamma(\Pi_t(x))]$. We have used the tubular neighbourhood decomposition \cref{eq:tubular-decom} and the following formulae:
\begin{align*}
    \partial_t d_\Gamma(t;x) = -\mathbf{w}(t;\Pi_t(x)) \cdot \nu_\Gamma(\Pi_t(x)), \; \nabla d_\Gamma(t;x) = \nu_\Gamma(\Pi_t(x)), \; x\in \mathcal{N}_{\Gamma},
\end{align*}
see \cite[Sec.~2]{MR2868564}. This gives us the following lemma:
\begin{lemma}For $h$ small enough, the map $ \Lambda^h(t; \cdot)|_{K_i(t)}$ is a $ C^{k+1}(K_i(t); \mathbb{R}^d)$ element-wise diffeomorphism with image $ \Lambda^h(t; \Omega^h_i(t)) = \overline{\Omega}_i(t)$. Moreover, define the following:
\begin{align*}
   \mathcal{J}_i^l(t) := \{ \Lambda^h(t; K_i(t))|\; K_i(t) \in \mathcal{J}^h_i(t) \}, \quad \mathcal{J}^l(t) := \mathcal{J}_1^l(t) \cup \mathcal{J}_2^l(t).
\end{align*}
Then $\mathcal{J}_1^l(t), \mathcal{J}_2^l(t), \mathcal{J}^l(t)$ define a uniform $k$-regular  evolving subdivision of $\overline{\Omega}_1(t), \overline{\Omega}_2(t), \overline{\Omega}$, respectively.
\end{lemma}
This follows from \cite[Lem.~8.12]{EllRan21} and the fact that facets are mapped to their evolving equivalent can be shown in the exact same way as in \cref{lemma:bulkdef}. A chart representing the full set-up is given in \cref{fig:chart-flows}. 

For a function $v^h \in H^h(t)$, the lift is denoted by $(\cdot)^l:H^h(t) \to H(t)$ and defined as follows: 
\begin{align*}
    v^{h,l}(x) := \left(v^h_1\left(t;[\Lambda^h(t;x)]^{-1}\right),v^h_2\left(t;[\Lambda^h(t;x)]^{-1}\right)\right).
\end{align*}
Its inverse will be labelled by $(\cdot)^{-l}$, i.e $(v^{h,l})^{-l}  =v^h$. Since $\Lambda^h(t;\cdot)|_{\overline{\Omega}^h_i(t)} \in W^{k+1, \infty}_T(\mathcal{J}^h_i(t); \mathbb{R}^d)$, with norm uniformly bounded in $h$, and invertible, via a similar change of variable method as \cref{lem:pf-compat} we have:
\begin{align*}
    c_1 \norm{v^{h,l}}_{H(t)} & \leq \norm{v^h}_{H^h(t)} \leq c_2 \norm{v^{h,l}}_{H(t)}
    && \text{for } v^h \in H^h(t) \\
    c_1 \norm{v^{h,l}}_{V(t)} & \leq \norm{v^h}_{V^h(t)} \leq c_2 \norm{v^{h,l}}_{V(t)}
    && \text{for } v^h \in V^h(t).
\end{align*}

We define the analogous flow $\mathbf{\Phi}^l_t: \overline{\Omega}_i(0) \to \overline{\Omega}_i(t)$ defined via the equation: $\mathbf{\Phi}^l_i(t;\Lambda^h(0;x)) = \Lambda^h(t;\mathbf{\Phi}^h_i(t;x))$. By the invertibility of $\Lambda^l$, this defines a flow, for which we can associate a push-forward map $\phi^l_t$ and inverse $\phi^l_{-t}$ as before. Note that this flow satisfies all properties \ref{B1} to \ref{B3} and \ref{D1} to \ref{D3} on the triplet $V(t) \subset H(t) \subset V^*(t)$ and therefore can be equipped with its own material derivative $\partial^l_t$:
\begin{align*}
    \partial_t^l \zeta := \phi^l_t \frac{\partial}{\partial t} \phi_{-t}^l \zeta,
\end{align*}
for $\zeta \in C^1_{(H, \phi^l_t)}$ (we make the flow $\phi_t^h$ explicit in the label for the space $ C^1_{(H, \phi^l_t)}$, so as to distinguish the space from $C^1_H$). Moreover, it is shown in \cite[Lem.~3.5]{EllRan21}:
\begin{align}\label{eq:commute}
  \partial^l_t \eta^{h,l} = (\partial^h_t \eta^h)^l,
\end{align}
for $\eta^h \in C^1_{H^h}$.

%Define the isoparametric mesh on the stationary domain (skim a little) emphasis that the boundaries match. Lagrangian nodes 

%The Evolving triangulation, flow + velocity, outline the property (quasi-uniform). 

%Define an interface element, (Remark, one can check that this will not cause issues, refer to the construction of the iso parametric elements previously), lift define. 
\subsection{Finite Element Spaces}

Let $\alpha(t)$ be a Lagrangian node and $\mathcal{J}_i(\alpha(t))$ be the set of elements in $K_i(t) \in \mathcal{J}_i(t)$ such that $\alpha(t) \in K_i(t)$, and let $\mathcal{N}^h_i(t)$ be the global set of all Lagrangian nodes in $\mathcal{J}_i^h(t)$. We introduce the finite dimensional subspace:
\begin{multline*}
  \mathcal{S}_i^h(t)
  := \biggl\{
    \chi^h_i = (\chi^h_i)_{K_i(t) \in \mathcal{J}^h_i(t)} \in \prod_{K_i(t) \in \mathcal{J}^h_i(t)}
    \{
    \widehat{\chi} \circ F^{-1}_{K_i(t)} : \, \widehat{\chi} \in \widehat{P}_k
    \} : \\
  \chi_i^h|_{K}(\alpha(t)) = \chi^h_i|_{K'}(\alpha(t)) \quad \text{for all} \, K_i(t), K_i'(t) \in \mathcal{J}_i(\alpha(t)), \forall \alpha(t) \in \mathcal{N}^h_i(t)
    \biggr\}.
\end{multline*}
Combining two copies of the space yields the adequate solution space:
\begin{multline*}
    \mathcal{S}^h (t):= \bigl\{\eta^h = (\eta^h_1, \eta^h_2) \in \mathcal{S}^h_1(t) \times \mathcal{S}^h_2(t) |\; \chi_1^h(\alpha(t)) = \chi_2^h(\alpha(t)) \\ \text{for all } \alpha(t) \in \Gamma^h(t) \cap \mathcal{N}^h_1(t) \; \text{and} \; \chi_2^h (\alpha(t)) = 0 \quad \forall \alpha(t) \in \partial \Omega \cap \mathcal{N}^h_2(t)\bigr\},
\end{multline*}
and we equip $\mathcal{S}^h(t)$ with the same norm as $V^h(t)$.
\begin{lemma}
 $(\mathcal{S}^h(t), \phi^h_t)|_{t \in I}$ form a compatible pair.
\end{lemma}
\begin{proof}
Since both the Lagrangian nodes and polynomials are evolved via $\mathbf{\Phi}_t^h$, one has by the definition of $\mathcal{S}(t)$, $\phi_t(\mathcal{S}(0)) = \mathcal{S}(t)$. Showing the remaining criterion for compatibility can be done in the same way as in \cref{lem:pf-compat}.
\end{proof}
Hence the moving spaces $L^2_{\mathcal{S}^h}$ is well defined.

The lifted solution space can now be defined as:
\begin{align*}
    \mathcal{S}^l(t) := \{ \chi^{h,l}|\; \chi^h \in \mathcal{S}^h(t) \}.
\end{align*}
The interpolation operator onto $\mathcal{S}^l(t)$, $I^l:C(\Omega) \to \mathcal{S}^l(t) $ can also be defined in a similar way:
\begin{align*}
    I^l(\eta)|_{K(t)} := \sum_{\{j : \chi^j \in P^k(t)\}} \eta(\alpha^{j,l}_{K(t)}) \chi^{j,l}.
\end{align*}
where $\{\alpha^{j,l}_{K(t)}\}_{j = 1}^{N(k)}$ are the lifted Lagrangian Nodes.

% Let:
% \begin{align*}
%     Z_k(t) := \{w \in V(t)|\; w_i \in H^{1+k}(\Omega_i(t))\}, \; \norm{w}^2_{Z_k(t)} := \sum_{i = 1}^2 \norm{w_i}^2_{H^{1+k}(\Omega_i(t))}.
% \end{align*} 
% Assume that the pair $(Z_k(t), \phi_t)|_{t \in I}$ is compatible (which follows from the flow being of added regularity $\mathbf{\Phi}_i(t; \cdot) \in C^{1 +k}(\Omega_i(0); \mathbb{R}^d))$. Then

The following variant of the approximation lemma holds:
\begin{lemma}[\cite{EllRan21}, Lem.~8.21]\label{lem:interp}
We have the estimates:
\begin{align*}
  \norm{w - I^l w}_{H(t)} +h \norm{w - I^lw}_{V(t)} & \leq c h^{k+1} \norm{w}_{Z_{k}(t)}, && \text{for } w \in Z_k(t), \\
  \norm{w - I^l w}_{H(t)} + h \norm{w - I^lw}_{V(t)} & \leq c h^2 \norm{w}_{Z_1(t)}, && \text{for } w \in Z_1(t).
\end{align*}
\end{lemma}
\section{Evolving finite element method}
\label{sec:fem-scheme}

\subsection{Scheme}
%introduce the finite dimensional problem, f,G such that the differentence between f^h, G^h acting on a lifted element is null.
For any $U^h, \, \zeta^h \in V^h(t)$, let:
\begin{align*}
    m^h(t; U^h, \zeta^h) &:= \sum_{i = 1}^2 \int_{\Omega^h_i (t)} U^h_i \zeta^h_i,\\
    a^h(t; U^h, \zeta^h) &:= \sum_{i = 1}^2 \int_{\Omega^h_i (t)} \mathcal{A}_i(t;\Lambda^{h}_l(t;x))\nabla U^h_i \cdot \nabla  \zeta^h_i + [\mathcal{B}(t;\Lambda^{h}_l(t;x)) - \mathbf{w}^h]\cdot \nabla U^h\, \zeta^h\\
    & \qquad\qquad\qquad\qquad + [\mathcal{C}(t; \Lambda^{h}_l(t;x)) - \nabla \cdot \mathbf{w}^h(t; x)] U^h \, \zeta^h,  \\
    l^h(t; \zeta^h) & := (f^{-l}J^h, \zeta^h)_{H^h(t)} + (g^{-l}\mu^h, \zeta^h)_{L^2(\Gamma^h(t))},
\end{align*}
where $J^h$, $\mu^h$ are the discrete Jacobians with respect to the lift maps $\Lambda^h(t;x)|_{\Gamma^h(t)}$, and by regularity of $\Lambda^h$, are of class $C^k(K_i(t); \mathbb{R}^d)$, $C^k(E(t); \mathbb{R}^{d}), \; \forall K_i(t)\in \mathcal{J}^h_i(t), \; \forall E(t) \in \mathcal{J}^h_\Gamma(t)$, respectively.

The finite element method seeks $U^h(t) \in \mathcal{S}^h(t)$ satisfying  the discrete variational problem:
\begin{align}\label{variational-eq}
    m^h(t; \partial^h_t U^h, \zeta^h) + a^h(t; U^h, \zeta^h) +\lambda^h(t; U^h, \zeta^h)&= l^h(t; \zeta^h) \quad \forall \zeta^h \in L^2_{\mathcal{S}^h}, \forall t \in I, \\
    \nonumber
    U^h(0) =U^h_0
     &:= \sum_{j = 1}^{dim(\mathcal{S}^h)} (u_0, \chi_j^h)_{H(0)} \chi_j^h.
\end{align}

\begin{remark}
  \label{rem:integrals-exact}
\normalfont
It might not be practical to calculate $l^h(t;\cdot)$ for an arbitrary pair $(f,g) \in L^2_H \times L^2_{\mathcal{V}_\Gamma}$ as it would have to be calculated via numerical integration. We write all integrals exactly to avoid additional technical challenges associated with analysing quadrature errors.
See \cite{Ciarlet1972THECE} for numerical integration on curved domains.
\end{remark} 

This formulation can be rearranged to a more useful form via the transport theorem with respect to the form $m^h(t; \cdot, \cdot)$:
\begin{equation}
    \label{eq:weak-form}
    \frac{d}{dt} m^h(t; U^h, \zeta^h) - m^h(t; U^h, \partial^h_t \zeta^h) +a^h(t; U^h, \zeta^h)= l^h(t; \zeta^h),
\end{equation}
for $\zeta^h \in C^1_{\mathcal{S}^h}$. Moreover, by construction of the $l^h(t;\cdot)$ term, for a function $\eta^h \in H^h(t)$:
\begin{align*}
    l(t; \eta^{h,l}) = l^h(t; \eta^{h}).
\end{align*}

\subsection{Well posedness of the finite element scheme}
%Show well posedness (duh)
~

\begin{theo}There exists a unique solution to \cref{variational-eq} with continuous bound:
\begin{align*}
    \sup_{t \in I}\norm{U^h}^2_{H^h(t)} +\int_0^T\norm{U^h}^2_{V^h(t)}\leq C(T)\bigg(\norm{U^{h}_0}^2_{H^h(0)} + \norm{f}_{L^2_H} +\norm{g}_{L^2_{\mathcal{V}_\Gamma^*}} \bigg).
\end{align*}
\end{theo}
\begin{proof}
Substituting the Ansatz:
\begin{align*}
    U^h(t;x) = \sum_{j= 1}^{\text{dim}(\mathcal{S}_h)} \alpha_j(t) \chi^h_j(t;x),
\end{align*}
where $\{\chi^h_j(t;x)\}^{\text{dim}(\mathcal{S}_h)}_{j = 1}$ are the basis functions of the evolving solution space $L^2_{\mathcal{S}^h}$. We refer to \cite[Lem.~3.1]{EllRan21} for a proof of the \textit{transport property}:
\begin{align*}
    \partial^h_t \chi^h_j = 0 \quad \forall j \in {1,...., N(k)}.
\end{align*}
Then the problem can be restated as the finite dimensional problem:
\begin{align*}
   \frac{d}{dt} \bigl(\mathbf{M}(t) \mathbf{\alpha}(t)\bigr) +\mathbf{A}(t)\mathbf{\alpha}(t) & = \mathbf{L}(t) \\
   \mathbf{\alpha}(0) & = \mathbf{\alpha}_0
\end{align*} 
where:
\begin{align*}
  & \mathbf{\alpha}(t) = (\alpha_1(t), \ldots, \alpha_{N(k)}(t)), \;
  && [\mathbf{M}(t)]_{j,k} = m^h(t;\chi^h_j,\chi^h_k), \; \\
  & [\mathbf{A}(t)]_{j,k} = a^h(t; \chi_j^h, \chi_k^h), \;
  && \mathbf{L}(t) = (l^h(t; \chi^h_1),...., l^h(t; \chi^h_{N(k)})).
\end{align*}
Note that $\mathbf{M}(t)$ is a Gram matrix (and hence invertible). Hence, by use of standard ODE theory (see \cite[Sec.~1.6]{Rou05}), there exists a solution $\alpha(t) \in \mathcal{W}^{1,1}(\mathbb{R}; \mathbb{R}^{N(k)})$. The uniform bound  and uniqueness follows from testing with  $U^h$ and using the transport theorem.
\end{proof}

\section{Error bound}
\label{sec:error}

%Put in the original assumptions between the discrete and full problem. Define the Ritz projection, state the results, especially Ritz

The main result of this article is the following optimal order error bound.
This result relies on structural assumptions on the partial differential equation, together with assumptions on the evolution of the physical and computational domains listed throughout this paper.
% The key assumption is that the given velocity is sufficiently nice to ensure that the mesh is uniformly quasi-uniform; in addition, Ass.~\ref{A1}, \ref{A2}, \ref{A3}, \ref{A4}, \ref{M1}, \ref{M2}, \ref{M3}, \ref{M4}, \ref{M5}, \ref{M6} are required. See \cref{rem:smooth-mesh} for some alternative approaches when this assumption might not hold.
% In addition, the following theorem requires uniform quasi-uniformity but this is not crucial for our scheme, but allows us to easily apply interpolation estimates throughout our analysis.
% See \cref{rem:smooth-mesh} for some alternative approaches when this assumption might not hold.

\begin{theo}\label{thm:error}
  Let Ass.~\ref{A1}, \ref{A2}, \ref{A3}, \ref{A4}, \ref{M1}, \ref{M2}, \ref{M3}, \ref{M4}, \ref{M5}, \ref{M6} hold together with the assumptions of \cref{lem:geom-regularity}. Further, we assume that the velocity is such that the triangulation is uniformly quasi-uniform.
If the solution to \cref{variational-eq} is of regularity $u \in W(Z_k, Z_k) \cap L^\infty_{Z_k}$ with uniform bound:
\begin{align*}
    \norm{u}_{L^\infty_{Z_k}} + \norm{u}_{W(Z_k, Z_k)} \leq C_u,
\end{align*}
then there exists a constant $\mathcal{C}$ depending on $C_u$ such that the following holds:
\begin{align*}
    \sup_{t \in I} \norm{u - U^{h,l}}^2_{H(t)} + h^2\int_0^T \norm{u - U^{h,l} }^2_{V(t)} \leq c\norm{u_0 - u_0^{h,l}}_{H(0)}^2 +ch^{2k+2} (\mathcal C).
\end{align*}
\end{theo}
Note that under the assumption of there existing a moving space equivalence on $W(Z_k,Z_k)$, $u \in L^\infty_{Z_k}$ automatically (see \cref{lem:equiv-material-derivative}) and hence it only suffices to assume $u \in W(Z_k,Z_k)$. In the next section we set out preliminary approximation results and then prove the error bound in the subsequent section.

In the next two subsections we introduce necessary tools from \cite[Sec.~3.3]{EllRan21} in order to obtain suitable orders of convergence. We will assume that for each space $W(Z_k,Z_k)$, there exists a moving space equivalence with $\mathcal{W}(Z_k(0), Z_k(0))$, see \cref{def:moving-equivalence}. This only requires the flow map $\mathbf{\Phi}_t$ to be regular enough and in particular is guaranteed if $\mathbf{\Phi}_t$ is smooth.

\subsection{Geometric perturbations}
 %Since we have two possible flows on $V(t)$, $\phi$ and $\phi^l$, we need to guarantee some compatibility condition that are reflected on the problem, to that end, note that we can define the derivative of the bilinear forms with respect to $\mathbf{\Phi}^l$ instead. 
 Let:
\begin{align*}
    b^l(t; w,v) &:= \frac{d}{dt}[a(t;w,v)] - a(t; \partial_t^l w,v) - a(t;w, \partial_t^l v),  
    && w,v \in W(V,V),\\
    b^h(t; w^h,v^h) &:= \frac{d}{dt}[a^h(t;w^h,v^h)] - a^h(t; \partial_t^h w^h,v^h) - a^h(t;w^h, \partial_t^l v^h),
    && w^h,v^h \in W(V^h,V^h),
\end{align*}
a.e for $t \in I$. Similarly as \cref{eq;derive}, these bilinear forms can be calculated explicitly and satisfy:
\begin{align*}
    |b^l(t; w,v)| \leq c \norm{w}_{V(t)} \norm{v}_{V(t)}, \; |b^h(t; w^h,v^h)| \leq c \norm{w^h}_{V(t)} \norm{v^h}_{V(t)}, \; \forall v,w \in V(t), \, \forall v^h,w^h \in V^h(t),
\end{align*}
for some constant $c$ independent of $t$ and $h$. Define $\lambda^l(t;\cdot , \cdot)$ to be the bilinear form of \cref{def:lambda} with respect to the flow $\mathbf{\Phi}^l$, which can be calculated to be:
\begin{align}\label{eq:lambda-disc-lif}
    \lambda^l(t;v, w) = \sum_{i = 1}^2 \int_{\Omega_i(t)} \nabla \cdot  \widetilde{\mathbf{w}} v\, w,
\end{align}
where $\widetilde{\mathbf{w}}$ is defined as:
\begin{align}\label{eq:vel-equiv}
   \widetilde{ \mathbf{w}}(t;x) = \frac{\partial}{\partial t} \mathbf{\Phi}^l(t;y)|_{y = \mathbf{\Phi}^l(-t;x)}.
\end{align}
See \cite[Lem.~8.15]{EllRan21}.
\begin{lemma}[\cite{EllRan21}, Lem.~8.16]\label{lem:jacobian}
    The lift satisfies the following:
    \begin{align*}
        \sup_{t \in I} \norm{\nabla \Lambda^h(t; \cdot) - \mathbf{I}}_{L^\infty(\Omega^h_i(t))} &\leq ch^k, \\
        \sup_{t \in I}\norm{\partial^h_t \nabla \Lambda^h(t; \cdot)}_{L^\infty(\Omega^h_i(t))} &\leq ch^k,
    \end{align*}
    and the Jacobian $J^h := \sqrt{\det{[\nabla \Lambda^h]^T\nabla \Lambda^h}}$ satisfies:
    \begin{align*}
        \sup_{t \in I}\norm{J^h(t; \cdot) - 1}_{L^\infty(\Omega^h_i(t))} &\leq ch^k .
    \end{align*}
\end{lemma}
Then the following holds for the bilinear forms introduced in \cref{eq:abstract} and \cref{eq;derive}:
\begin{prop}%[Lemma 8.23 $\&$ 8.24 \cite{EllRan21}]%
There exists a constant $c >0$ such that for almost all $t \in I$ and for all $w^h,v^h \in V^h(t)$, $w^{h,l}, v^{h,l} \in V(t)$ the following error bounds hold:
\begin{align}
\label{eq:m-error} \tag{P1}
|m( t; w^{h,l}, v^{h,l} ) - m^h( t; w^h, v^h ) |& \leq c h^{k+1} \norm{w^{h,l}}_{V(t)} \norm{v^{h,l}}_{V(t)}, \\
\label{eq:c-error} \tag{P2}
|\lambda( t; w^{h,l}, v^{h,l} ) - \lambda^h( t; w^h, v^h )| &\leq c h^{k+1} \norm{w^{h,l}}_{V(t)}\norm{v^{h,l}}_{V(t)}, \\
\label{eq:ct-error} \tag{P3}
|\lambda^l( t; w^{h,l}, v^{h,l} ) - \lambda( t; w^{h,l}, v^{h,l} )| &\leq c h^{k} \norm{w^{h,l}}_{V(t)}\norm{v^{h,l}}_{V(t)}, \\
\label{eq:a-error} \tag{P4}
|a( t; w^{h,l}, v^{h,l} ) - a^{h}( t; w^h, v^h )|& \leq c h^{k} \norm{w^{h,l}}_{V(t)} \norm{v^{h,l}}_{V(t)}, \\
\label{eq:b-error} \tag{P5}
|b^l( t; w^{h,l}, v^{h,l}) - b^h(t;w^h,v^h)|&\leq ch^k \norm{w^{h,l}}_{V(t)} \norm{v^{h,l}}_{V(t)}, \\
\label{eq:bt-error} \tag{P6}
|b^l( t; w^{h,l}, v^{h,l} ) - b( t; w^{h,l}, v^{h,l} )|& \leq c h^{k} \norm{w^{h,l}}_{V(t)} \norm{v^{h,l}}_{V(t)}.
\end{align}
  For $\eta, \zeta \in Z_1(t)$ with inverse lifts $\eta^{-l}, \zeta^{-l}$:
  \begin{align}
    \label{eq:a-error2} \tag{\ref*{eq:a-error}'}
    |a( t; \eta, \zeta ) - a^{h}( t; \eta^{-l}, \zeta^{-l} ) |& \leq ch^{k+1}\norm{\eta}_{Z_1(t)}\norm{\zeta}_{Z_1(t)}, \\
    \label{eq:b-error2} \tag{\ref*{eq:b-error}'}
    |b^l( t; \eta, \zeta) - b^{h}( t; \eta^{-l}, \zeta^{-l})|& \leq c h^{k+1} \norm{\eta}_{Z_k(t)} \norm{\zeta}_{Z_k(t)}.
  \end{align}
  For $\eta \in C^1_{Z_k}$ and $\zeta \in Z_1(t)$, with inverse lifts $\eta^{-l}$ and $\zeta^{-l}$:
  \begin{align}
    \label{eq:amd-error2} \tag{P7}
    |a( t; \partial_t^l \eta, \zeta ) - a^{h}( t; \partial_t^h  \eta^{-l}, \zeta^{-l} )|
    &\leq c h^{k+1} ( \norm{\eta}_{Z_1(t)} +\norm{\partial_t^\bullet\eta}_{Z_1(t)} )\norm{\zeta}_{Z_1(t)}.
  \end{align}
  The material derivatives satisfy
  \begin{alignat}{2}
    \label{eq:md-error} \tag{P8}
    \norm{\partial_t^l \zeta - \partial_t^\bullet\zeta}_{H(t)}
    & \leq c h^{k+1} \norm{\zeta}_{V(t)}
    && \mbox{ for } \zeta \in C^1_{V}, \\
    \label{eq:mdV-error} \tag{P9}
  \norm{\partial^l_t \zeta - \partial^\bullet_t \zeta}_{V(t)}
    & \leq c h^{k} \norm{\zeta}_{Z_1(t)}
    && \mbox{ for } \zeta \in C^1_{Z_1}.
  \end{alignat}
\end{prop}
\begin{proof}
In \cite[Lem.~8.23 and 8.24]{EllRan21}, these estimates are proven on a single evolving domain $\Omega(t)$ with $V(t) = H^1(\Omega(t))$. However, almost the same arguments  cover our case. We will only show this for \labelcref{eq:m-error},\labelcref{eq:a-error} and \labelcref{eq:mdV-error}, but the same method can be applied for the remaining claims.

\labelcref{eq:m-error}: Let $J^h := \sqrt{\det{[\nabla \Lambda^h]^T\nabla \Lambda^h}}$ be Jacobian resulting from switching from $\Omega_i(t)$ to $\Omega_i^h(t)$. Then, the lift itself differs from the identity only when $x$ is in an interface element: let
\begin{align*}
    M: = \{x \in \Omega^h_i(t): \, J^h(x) \neq 1\} \subset  \{x \in \Omega_i(t): \, |d_\Gamma(t;x)| \leq h\}, 
\end{align*}
and let $M^l:= \{\Lambda^h(t;x)|\; x \in M\}\subset \{x \in \Omega_i(t): \, |d_\Gamma(t;x)| \leq h\} $. For $w^h_i \in H^1(\Omega_i(t))$:
\begin{align*}
   \left| \int_{\Omega_i(t)} w^{h,l}_i \cdot v^{h,l}_i - \int_{\Omega^h_i(t)}  w^{h}_i \cdot v^{h}_i \right| = \left| \int_{M}  w^{h}_i \cdot v^{h}_i[J^h - 1] \right|\leq ch^{k} \norm{ w^{h,l}_i}_{L^2(M^l)}\norm{ v^{h,l}_i}_{L^2(M^l)}.
\end{align*}
Then, via the Narrow-Band trace inequality, see \cite[Lem.~4.10]{EllRan12}, we see:
\begin{align*}
    \norm{ w^{h,l}_i}_{L^2(M^l)} \leq ch^{1/2} \norm{ w^{h,l}_i}_{H^1(\Omega_i(t))},
\end{align*}
and hence:
\begin{align*}
    |m( t; w^{h,l}, v^{h,l} ) - m^h( t; w^h, v^h ) | &= \left| \sum_{i = 1}^2\int_{\Omega_i(t)} w^{h,l}_i \cdot v^{h,l}_i - \int_{\Omega^h_i(t)}  w^{h}_i \cdot v^{h}_i \right|,\\
    &\leq c h^{k+1} \sum_{i =1}^2 \norm{w^{h,l}}_{H^1(\Omega_i(t))} \norm{v^{h,l}}_{H^1(\Omega_i(t))},\\
    &\leq c h^{k+1} \norm{w^{h,l}}_{V(t)} \norm{v^{h,l}}_{V(t)}.
\end{align*}

\labelcref{eq:a-error}: Domain-wise (using the shorthand $\Lambda^h(t;x) = \Lambda^h$):
\begin{align*}
    & \int_{\Omega_i (t)}\mathcal{A}_i (t;x) \nabla w^{h,l}_i\cdot \nabla v^{h,l}_i +[\mathcal{B}_i(t;x)- \mathbf{w}] \cdot \nabla w^{h,l}_i\, v^{h,l}_i +[\mathcal{C}_i(t;x) - \nabla \cdot \mathbf{w}]w^{h,l}_i\, v^{h,l}_i \\
    &= \int_{\Omega^h_i (t)}[\mathcal{A}_i (t;\Lambda^h) \nabla \Lambda^h\nabla w^{h}_i\cdot \nabla \Lambda^h\nabla v^{h}_i +[\mathcal{B}_i(t;\Lambda^h)- \mathbf{w}(t; \Lambda^h)] \cdot \nabla \Lambda^h\nabla w^{h}_i\, v^{h}_i\\
    &\qquad\qquad +[\mathcal{C}_i(t;\Lambda^h) - \text{Tr}(\nabla \Lambda^h\nabla\mathbf{w})]w^{h}_i\, v^{h}_i]J^h.
\end{align*}
By use of both \cref{lem:interp,lem:jacobian}, \ref{eq:a-error} follows in the same way as \labelcref{eq:m-error}.

\labelcref{eq:mdV-error}: Explicitly expanding both material derivatives:
\begin{align}\label{eq:material-derive}
    \partial_t^\bullet \zeta= \partial_t \zeta + \mathbf{w}\cdot \nabla \zeta = \partial_t^l \zeta +[\mathbf{w} -  \widetilde{\mathbf{w}}]\cdot \nabla \zeta.
\end{align}
Using the definition of $\mathbf{\Phi}^l_t$ and \cref{eq:vel-equiv}, element-wise, for $x \in K(t)$:
\begin{align*}
    \partial_t^h \Lambda^h(t;x) = \phi_t^h \partial_t \Lambda^h(t; \mathbf{\Phi}^h(t;x)) = \phi_t^h \partial_t \mathbf{\Phi}^l(t; \Lambda^h(0;x)) = \phi_t^h \widetilde{\mathbf{w}}(t; \mathbf{\Phi}^l(t; \Lambda^h(0;x))) = \widetilde{\mathbf{w}}(t; \Lambda^h(t;x)).
\end{align*}
Rewriting \cref{eq:material-derive}:
\begin{align}\label{eq:p9-1}
     \partial_t^\bullet \zeta - \partial_t^l \zeta &= [\mathbf{w} -   \partial_t^h \Lambda^h(t;z)|_{z = [\Lambda^h(t;x)]^{-1}}]\cdot \nabla \zeta, \nonumber\\
     &= [\mathbf{w} -  I^l \mathbf{w} +  I^l\mathbf{w}  - \partial_t^h \Lambda^h(t;z)|_{z = [\Lambda^h(t;x)]^{-1}}]\cdot \nabla \zeta,
\end{align}
via \cref{lem:interp}:
\begin{align}\label{eq:p9-2}
    \norm{\mathbf{w} -  I^l \mathbf{w}}_{V(t)} \leq ch^k \norm{\mathbf{w}}_{Z^{k}(t)},
\end{align}
as for the remaining term, $\partial_t^h \Lambda^h(t;y)$, using \cref{eq:lift-deriv}:
\begin{align*}
    &I^l\mathbf{w}  - \partial_t^h \Lambda^h(t;z)|_{z = [\Lambda^h(t;x)]^{-1}} = [\mathbf{w}^h(t;z)  - \partial_t^h \Lambda^h(t;z)]_{z = [\Lambda^h(t;x)]^{-1}}\\
    &=\begin{cases} 
  -\widetilde{\mu}(\hat{x})^{k+2}\left[ \left(\mathbf{w}^h(t;y) - \mathbf{w}(t;\Pi_t(y))\right)\cdot \nu_\Gamma(\Pi_t(y)) \nu_\Gamma(\Pi_t(y)) - d_\Gamma(t;y)T(y)\right]\, &\text{if} \, x \notin \sigma,\\
     0\, &\text{otherwise}.
      \end{cases}
\end{align*}
Via the use of standard geometric estimates, see \cite[Lem.~8.16 and 9.10]{EllRan21} and the fact that $\mathbf{w}^h$ is the interpolant of $\mathbf{w}$, we infer that:
\begin{align}\label{eq:p9-3}
    \norm{\mathbf{w}^h(t;z)  - \partial_t^h \Lambda^h(t;z)}_{W^{1,\infty}(K(t))} \leq ch^{k}.
\end{align}
Combining \cref{eq:p9-1,eq:p9-2,eq:p9-3} yields \labelcref{eq:mdV-error}.
\end{proof}
\subsection{Ritz Projection}
For the Ritz projection, it is convenient to work with a strictly coercive bilinear form.
To this end, we set
\begin{equation}a_\kappa(t;w, v):=\kappa m(t; w, v)+a(t;w, v)\end{equation}
and observe  that
\begin{multline*}
  a_\kappa(t;v, v) \geq \bigg(\gamma  - \frac{\epsilon \norm{\mathcal{B}_i - \mathbf{w}}_{L^\infty(\Omega)}}{2}\bigg)\norm{\nabla v}^2_{H(t)} \\
  +\bigg(\kappa - \norm{\nabla \cdot \mathbf{w}}_{L^\infty(\Omega)} - \norm{\mathcal{C}_i}_{L^\infty(\Omega)} - \frac{\norm{\mathcal{B}_i - \mathbf{w}}_{L^\infty(\Omega)}}{2 \epsilon} \bigg) \norm{v}^2_{H(t)},
\end{multline*}
where $\gamma$ is a lower bound for the eigenvalues of $\mathcal{A}_i$ \cref{eq:Agamma}.
Thus taking $\epsilon$ sufficiently small and $\kappa$ sufficiently large, we may choose $\kappa$ depending only on the data to ensure that the bounded bilinear form  $a_\kappa(t;\cdot,\cdot)$ is strictly coercive.

Similarly, we define
\[
  a^h_\kappa(t; w_h, v_h) := \kappa m^h(t; w_h, v_h) + a^h(t; w_h, v_h),
\]
which is also coercive provided $\kappa$ is large enough, independently of $h$, by the same argument.

The \textit{Ritz projection} 
$\Pi^h\colon V(t) \to \mathcal{S}^h(t)$ is defined as the solution to:
\begin{align}
  \label{eq:ritz}
    a^h_\kappa(t; \Pi^h(\eta), v^h) = a_\kappa(t; \eta, v^{h,l}) \quad \forall v^h \in \mathcal{S}^h(t),
\end{align}
and $\pi^h \eta := (\Pi^h \eta)^l$. By the coercivity and boundedness of $a^h_\kappa$, this gives us a uniformly bounded and linear operator $\Pi^h\colon V(t) \to \mathcal{S}^h(t)$. Moreover, it is further proven in \cite[Lem.~3.9]{EllRan21} , that $\Pi^h \eta \in C^1_{\mathcal{S}^h}$ if $\eta \in C^1_{V}$ and it follows by use of the same method that $\Pi^h \eta \in C^0_{\mathcal{S}^h}$ if $\eta \in C^0_{V}$.

\begin{lemma}\label{lem:ritz-additional-regular}
    The Ritz projection can be extended as a continuous linear operator $\Pi^h(\cdot):L^2_V \to L^2_{\mathcal{S}^h}$. Moreover, if $\eta \in W(V,V)$, then $\Pi^h(\eta) \in W(\mathcal{S}_h, \mathcal{S}_h) = \{v^h \in L^2_{\mathcal{S}^h}, \; \partial_t^h v^h \in L^2_{\mathcal{S}^h}\}$. In particular, $\Pi^h(\cdot):W(V,V) \to  W(\mathcal{S}_h, \mathcal{S}_h)$ is linear and bounded uniformly in $h$.
\end{lemma}
\begin{proof}
    Integrating \cref{eq:ritz}, we see, that for all $v^h \in L^2_{\mathcal{S}^h}$:
    \begin{align*}
        \int_0^T a^h_\kappa(t; \Pi^h(\eta), v^h)\,dt = \int_0^T a_\kappa(t; \eta, v^{h,l})\,dt.
    \end{align*}
    By point-wise coercivity in time of $a_\kappa^h$, we also get the coercivity over $L^2_{\mathcal{S}^h}$. Since $L^2_{\mathcal{S}^h}$ is a closed subspace of a Hilbert space, and $\int_0^T a_\kappa(t; \eta, (\cdot)^{l})\,dt$ defines a bounded linear functional on $L^2_{\mathcal{S}^h}$, via the standard use of Lax-Milgram, there exists a unique solution, labelled $\Pi^h(\eta)$ and we achieve the bound:
    \begin{align}\label{eq:ritz-bound-1}
        \norm{\Pi^h(\eta)}_{L^2_{\mathcal{S}^h}} \leq \sup_{v^h \in L^2_{\mathcal{S}^h}, \, \norm{v^h}=1} \int_0^T a_\kappa(t; \eta, v^{h,l})\,dt \leq c \norm{\eta}_{L^2_V}.
    \end{align}
    Hence $\Pi^h(\cdot)$ is continuous. Note that the bound in \cref{eq:ritz-bound-1} can be taken to be independent of $h$, this is due to the fact that both the bilinear form $a_\kappa$ and the lift map $(\cdot)^l$ are both bounded independently of $h$.
 
    To show the second claim, assume $\eta \in W(V,V)$ and set $\zeta$ as the solution to:
    \begin{align}\label{eq:time-reg-ritz}
         \int_0^T a^h_\kappa(t; \zeta, v^h)\,dt =  \int_0^T a_\kappa(t; \partial_t^l \eta, v^{h,l}) - b^h_\kappa(t; \Pi^h \eta, ,v^h) + b_\kappa^l(t; \eta, v^{h,l})\, dt, \; \forall v^h \in L^2_{\mathcal{S}^h}.
    \end{align}
    Then via the same argument as before, there exists a unique $\zeta \in L^2_{\mathcal{S}^h}$ solving \cref{eq:time-reg-ritz} with bound:
    \begin{align}\label{eq:ritz-bound-2}
        \norm{\zeta}_{L^2_{\mathcal{S}^h}} \leq c\norm{\eta}_{W(V,V)},
    \end{align}
    and similarly as before, the bound in equation \cref{eq:ritz-bound-2} is independent of $h$. Since $\eta \in W(V,V)$, it is also in $C_V^0$ by \cref{lem:equiv-material-derivative}, $\eta \in C_V^0$ and therefore $\Pi(\eta) \in C^0_{\mathcal{S}^h}$.
    
    Define:
    \begin{align*}
        w^h := \phi^h_t \int_0^t \phi^h_{-s} \zeta(s) \,ds + \phi^h_t\Pi^h(\eta)(0).
    \end{align*}
    Via the isomorphism lemma (\cref{lem:isomo}), $\phi^h_{-s} \zeta(s) \in L^2(I; \mathcal{S}^h(0))$, the standard Bochner space and therefore is Bochner integrable. Since $\mathcal{S}^h(0)$ is a closed linear subspace, the definite Bochner integral of a function inside $\mathcal{S}^h(0)$ remains in $\mathcal{S}^h(0)$ for all $t \in I$. Using the isomorphism again, we see that $w^h \in L^2_{\mathcal{S}^h}$.
    We will show that $w^h = \Pi^h \eta$ which will show the second claim.
    Indeed, $w^h \in W(\mathcal{S}^h, \mathcal{S}^h)$ with $\partial_t^h w^h = \zeta$. Substituting this back into \cref{eq:time-reg-ritz}:
    \begin{align}\label{eq:w-definition}
        \int_0^T a^h_\kappa(t; \partial_t^h w , v^h)\,dt =  \int_0^T a_\kappa(t; \partial_t^l \eta, v^{h,l}) - b^h_\kappa(t; \Pi^h \eta, ,v^h) + b_\kappa^l(t; \eta, v^{h,l})\, dt, \; \forall v^h \in L^2_{\mathcal{S}^h},
    \end{align}
    using the definition of $b^h_\kappa$, we see that, for $v^h \in \mathcal{W}(\mathcal{S}^h,\mathcal{S}^h)$ with $v^h(0) = v^h(T) = 0$:
   \begin{align*}
        \int_0^T a^h_\kappa(t; \partial_t^h w^h , v^h)\,dt &=  \int_0^T -b^h_\kappa(t; w^h,v^h) - a^h_\kappa (t; w^h ,\partial^h_t v^h)\,dt, \\
        \int_0^T a_\kappa(t; \partial_t^l \eta, v^{h,l})  \,dt &=\int_0^T -b^l_\kappa(t; \eta, v^{h,l}) - a_\kappa(t; \eta, [\partial_t^h v^h]^l)\,dt,
    \end{align*}
    using the commutation properties of the material derivatives $\partial^h_t$ and $\partial_t^l$ (see \cref{eq:commute}). Substituting these expression back into \cref{eq:w-definition}:
    \begin{align}\label{eq:ritz-equation1}
        \int_0^T a_\kappa (t; \eta,  [\partial_t^h v^h]^l) - a^h_\kappa(t; w^h, \partial_t^hv^h)\,dt = \int_0^T b^h_\kappa(t; w^h- \Pi^h \eta, v^h)\,dt.
    \end{align}
    Using the definition of the Ritz projection \cref{eq:ritz}, we arrive at:
    \begin{align}\label{eq:ritz-equation2}
        \int_0^T a^h_\kappa(t;\Pi^h \eta -  w^h, \partial_t^hv^h)\,dt = \int_0^T b^h_\kappa(t; w^h- \Pi^h \eta, v^h)\,dt.
    \end{align}
    Testing this equation with $v^h(t) = \varphi(t)\psi^h(t)$ with $\varphi(t) \in \mathcal{D}(I)$ and $\psi^h(t) \in W(\mathcal{S}^h, \mathcal{S}^h)$, we see:
    \begin{align*}
         \int_0^T\varphi(t) a^h_\kappa(t;\Pi^h \eta -  w^h, \partial_t^h\psi^h)+ \varphi'(t) a^h_\kappa(t;\Pi^h \eta-  w^h,\psi^h)\,dt  = \int_0^T \varphi(t) b^h_\kappa(t; w^h- \Pi^h \eta, \psi^h)\,dt.
    \end{align*}
    Then, via the fundamental lemma of variational calculus, see \cite[Thm.~1.2.1 and Lem.~1.2.1]{MR1674720}, and the fact that $a^h_\kappa(t;\Pi^h \eta -  w^h,\psi^h)$ is continuous as $\Pi^h(\eta) - w^h, \psi^h \in C^0_{\mathcal{S}^h}$, we get that for all $t \in I$: 
    \begin{align}\label{eq:ritz-equation3}
        a^h_\kappa(t;\Pi^h \eta -  w^h,\psi^h) = \int_0^t b^h_\kappa(s; \Pi^h \eta-w^h, \psi^h) +a^h_\kappa(s;\Pi^h \eta -  w^h, \partial_s^h\psi^h)\,ds
    \end{align}
    Fix $t \in I$ and test with $\psi^h(s) = \phi^h_{s} \phi^h_{-t}(\Pi^h \eta -  w^h)$, we see that $\partial_s^h \psi^h(s) = 0$ and hence, via the coercivity of $a_k$ and compatibility:
    \begin{align*}
        \norm{\Pi^h \eta -  w^h}^2_{\mathcal{S}^h(t)} &\leq c \int_0^t b^h_\kappa(s;  \Pi^h \eta-w^h, \psi^h)\,ds,\\
        &\leq c \int_0^t \norm{\Pi^h \eta -  w^h}_{\mathcal{S}^h(s)}\norm{\Pi^h \eta -  w^h}_{\mathcal{S}^h(t)}\,ds.
    \end{align*}
    By use of Young's inequality:
    \begin{align*}
        \norm{\Pi^h \eta -  w^h}^2_{\mathcal{S}^h(t)} &\leq c \int_0^t \norm{\Pi^h \eta -  w^h}^2_{\mathcal{S}^h(s)}\,ds.
    \end{align*}
    This holds for arbitrary point $t \in I$ and hence can be repeated to see that this holds for all of $I$. By use of Gr\"{o}nwall's inequality it must be that:
    \begin{align*}
         \norm{\Pi^h \eta -  w^h}^2_{\mathcal{S}^h(t)}  = 0,
    \end{align*}
    and hence $\Pi^h\eta = w^h$ and therefore $\Pi^h \eta \in W(\mathcal{S}^h, \mathcal{S}^h)$. As for the bound on $\Pi^h \eta$, we see that, since $\partial_t^h w^h = \partial_t^h \Pi^h(\eta) = \zeta$, using \cref{eq:ritz-bound-1,eq:ritz-bound-2} we see that $\norm{\Pi^h(\eta)}_{W(\mathcal{S}^h, \mathcal{S}^h)} \leq c \norm{\eta}_{W(V,V)}$ independently of $h$.
\end{proof}
Note that $a_\kappa(t; \cdot, \cdot)$ and $a_\kappa^h(t; \cdot, \cdot)$ satisfy all the same estimates as $a(t;\cdot,\cdot)$ and $a^h(t; \cdot, \cdot)$ in \labelcref{eq:m-error}--\labelcref{eq:mdV-error}. 
\begin{remark}\label{rem:lift-ritz-continuity}
    \normalfont
    We note that via the commutative properties of the material derivatives (see \cref{eq:commute}) $[\partial_t^h \Pi^h(v)]^l = \partial_t^l \pi^h(v)$ for $v \in W(V,V)$. Hence we also conclude that the lifted Ritz map $\pi^h:W(V,V) \to W(\mathcal{S}_h^l, \mathcal{S}_h^l)$ is continuous as well, again uniformly in $h$. 
\end{remark}
Define the dual solution operators $\mathcal{R}_{H}\colon H(t) \to V(t)$ and $\mathcal{R}_{\mathcal{H}_\Gamma}\colon \mathcal{H}_\Gamma(t) \to V(t)$ to be the solutions to:
\begin{align}\label{eq:dual-def}
    a_\kappa(t; w, \mathcal{R}_H(v)) &= m(t;w,v)
    && \forall w \in V(t)\\
    a_\kappa(t; w, \mathcal{R}_{\mathcal{H}_\Gamma}(v)) &= (w,v)_{\mathcal{H}_\Gamma(t)}
    && \forall w \in V(t).
\end{align}
$\mathcal{H}_\Gamma(t), \mathcal{V}_\Gamma(t)$ are defined in \cref{sec:realisation}. We aim to show that these operators satisfy the following regularity condition:
\begin{lemma}\label{lem:adjoint-reg}
   Assuming additional regularity on the data:
   $\mathcal{A}_i \in W^{1,\infty}(\Omega_i(t); \mathbb{R}^{d \times d})$,
   $\mathcal{B}_i \in W^{1,\infty}(\Omega_i(t); \mathbb{R}^d)$,
   $\mathcal{C}_i \in W^{1,\infty}(\Omega_i(t); \mathbb{R})$, the operators $\mathcal{R}_H, \; \mathcal{R}_{\mathcal{H}_\Gamma}$ satisfy the regularity bounds:
    \begin{align}
      \norm{\mathcal{R}_H(v)}_{Z_1(t)} & \leq c \norm{v}_{H(t)} \\
      \label{eq:dual-reg-gamma}
      \norm{\mathcal{R}_{\mathcal{H}_\Gamma}(v)}_{Z_1(t)} & \leq c \norm{v}_{\mathcal{H}_\Gamma(t)}.
\end{align}
\end{lemma}
\begin{proof}
    Writing \cref{eq:dual-def} explicitly, we seek a solution $ \mathcal{R}_H(v)$ to:
\begin{align}\label{eq:dual-set-problem}
    \sum_{ i = 1}^2 \int_{\Omega_i(t)} \mathcal{A}_i \nabla w_i \cdot \nabla \mathcal{R}_H(v) + [\mathcal{B}_i - \mathbf{w}] \cdot \nabla w_i \mathcal{R}_H(v) + [\mathcal{C}_i+ \kappa - \nabla \cdot \mathbf{w}] w_i \mathcal{R}_H(v)  = m(t;w,v).
\end{align}
We note, by increasing $\kappa$ more if necessary, \cref{eq:dual-set-problem} is still coercive. Hence there exists a solution $\mathcal{R}_H(v) \in V(t)$ via use of the Babuska-Lax-Milgram theorem. Moreover, there exists a constant $c$ independent of time such that:
\begin{align*}
    \norm{\mathcal{R}_H(v)}_{V(t)} \leq c\norm{v}_{H(t)}.
\end{align*}
To show the additional regularity, rearranging \cref{eq:dual-set-problem}, we have:
\begin{align}\label{eq:dual-set-problem-2}
    a_\kappa(t; w, \mathcal{R}_H(v)) =  a_\kappa(t;\mathcal{R}_H(v), w)   +  \sum_{ i = 1}^2 \int_{\Omega_i(t)}[\mathcal{B}_i - \mathbf{w}] \cdot ([\nabla w_i] \mathcal{R}_H(v) - [\nabla \mathcal{R}_H(v)]w_i).
\end{align}
The remaining term of \cref{eq:dual-set-problem} can be further rearranged as, using integration by parts and the continuity of $V(t)$ and $\mathbf{w}$ across the interface:
\begin{multline*}
  \sum_{ i = 1}^2 \int_{\Omega_i(t)}[\mathcal{B}_i - \mathbf{w}] \cdot ([\nabla w_i] \mathcal{R}_H(v) - [\nabla \mathcal{R}_H(v)]w_i) \\
  = \int_{\Gamma(t)} \llbracket \mathcal{B}\rrbracket \cdot \nu  w \mathcal{R}_H(v) -  \sum_{ i = 1}^2  \int_{\Omega_i(t)}  2[\mathcal{B}_i - \mathbf{w}]  \cdot [\nabla \mathcal{R}_H(v)]w_i + \nabla \cdot [\mathcal{B}_i - \mathbf{w}] \mathcal{R}_H(v) w_i.
\end{multline*}
Hence, the solution to \cref{eq:dual-def} also solves:
\begin{multline}\label{eq:equation-for-R}
    a_\kappa(t;\mathcal{R}_H(v), w) = m(t;w, v) -\int_{\Gamma(t)} \llbracket \mathcal{B}\rrbracket \cdot \nu  w \mathcal{R}_H(v)\\
    +\sum_{ i = 1}^2  \int_{\Omega_i(t)}  2[\mathcal{B}_i - \mathbf{w}]  \cdot [\nabla \mathcal{R}_H(v)]w_i + \nabla \cdot [\mathcal{B}_i - \mathbf{w}] \mathcal{R}_H(v) w_i, \quad \forall w \in V(t).
\end{multline}
Set $\mathcal{L}(v)$ to be the solution to:
\begin{multline}\label{eq:equation-for-L}
     a_\kappa(t;\mathcal{L}(v), w) 
     =  \underbrace{-\int_{\Gamma(t)} \llbracket \mathcal{B}\rrbracket \cdot \nu  w \mathcal{R}_H(v)}_{:= ( \widetilde{g}, w)_{\mathcal{H}_\Gamma(t)}}\\
     + \underbrace{\sum_{ i = 1}^2  \int_{\Omega_i(t)}  w_i v_i + 2[\mathcal{B}_i - \mathbf{w}]  \cdot [\nabla \mathcal{R}_H(v)]w_i + \nabla \cdot [\mathcal{B}_i - \mathbf{w}] \mathcal{R}_H(v) w_i}_{:= (\widetilde{f}, w)_{H(t)}}, \quad \forall w \in V(t).
\end{multline}
This solution exists by \cref{theorem:existence}. We seek to show that first $\mathcal{L}(v) \in Z_1(t)$ and then that $\mathcal{R}_H(v) = \mathcal{L}(v)$. By use of Theorem 1 in \cite{nistor}, since $\widetilde{f} \in H(t)$ and $\widetilde{g} \in H^{1/2}(\Gamma(t))$, the solution to \cref{eq:equation-for-L} is indeed in $Z_1(t)$ since the data $(\mathcal{A}_i, \mathcal{B}_i, \mathcal{C}_i)$ is regular enough and moreover:
\begin{align*}
    \norm{\mathcal{L}(v)}_{Z_1(t)} \leq c (\norm{v}_{\mathcal{V}_\Gamma(t)}  +\norm{\mathcal{R}_H(v)}_{V(t)})\leq c \norm{v}_{H(t)}
\end{align*}
To show $\mathcal{R}_H(v) = \mathcal{L}(v)$, subtracting \cref{eq:equation-for-L} from \cref{eq:equation-for-R}:
\begin{align*}
    a_\kappa(t;\mathcal{R}_H(v) - \mathcal{L}(v), w)  = 0.
\end{align*}
Testing with $w = \mathcal{R}_H(v) - \mathcal{L}(v)$ and using coercivity yields $\mathcal{R}_H(v) = \mathcal{L}(v)$ and hence, we have:
\begin{align}\label{dual-ineq-added}
    \norm{\mathcal{R}_H(v)}_{Z_1(t)} \leq c\norm{v}_{H(t)},
\end{align}
where the regularity constant $c$ can be taken to be bounded on $[0, T]$ via the regularity of the flow.

The same argument follows for $\mathcal{R}_{\mathcal{H}_\Gamma}$ with 
\begin{align*}
    (\tilde{f}, w)_{H(t)} & = \sum_{ i = 1}^2 \int_{\Omega_i(t)}  2[\mathcal{B}_i - \mathbf{w}]  \cdot [\nabla \mathcal{R}_{\mathcal{H}_\Gamma}(v)]w_i + \nabla \cdot [\mathcal{B}_i - \mathbf{w}] \mathcal{R}_H(v) w_i \\
    (\tilde{g}, w)_{\mathcal{H}_\Gamma(t)} & = \int_{\Gamma(t)} w v - \llbracket \mathcal{B}\rrbracket \cdot \nu  w \mathcal{R}_{\mathcal{H}_\Gamma}(v).
    \qedhere
\end{align*}
\end{proof}
% \textcolor{magenta}{[rhs should be $m$ and domain should be $H(t)$]}\\
% \textcolor{blue}{A, I see, I will fix this tomorrow}\\
% { \color{red}
% Key idea: Aubin-Nitsche trick:

% Without domain perturbation, for any $v \in V$
% \begin{align*}
%     \norm{v}_{H(t)}^2 & = m(t; v, v) \\
%     & = a_\kappa(t; v, \mathcal{R}(v)) \\
%     & = a_\kappa(t; v, \mathcal{R}(v) - I_h \mathcal{R}(v)) && \text{[Galerkin orthogonality]} \\
%     & \le \norm{v}_{V(t)} \norm{\mathcal{R}(v) - I_h \mathcal{R}(v)}_{V(t)} \\
%     & \le \norm{v}_{V(t)} c h \norm{\mathcal{R}(v)}_{Z_1(t)} \\
%     & \le \norm{v}_{V(t)} c h \norm{v}_{H(t)},
% \end{align*}
% so
% \[
% \norm{v}_{H(t)} \le c h \norm{v}_{V(t)}.
% \]
% Plug in, $v =$ Ritz error gives improved $L^2$ error estimate.

% We don't have Galerkin orthogonality so there is an additional domain perturbation estimate required but the structure is the same.

% This is a different dual problem to the one considered in Lemma 5.4.
% }
\begin{lemma}\label{lem:b-interp}
On top of the assumptions made in \cref{lem:adjoint-reg}, assume $\mathcal{A}_i \in C^2(\overline{\mathcal{Q}}_i; \mathbb{R}^{d \times d})$. For $w \in Z_k(t)$, $\eta  := w - \pi^h w$ and $v \in Z_1(t)$ it holds that 
\begin{align}\label{eq:bil}
    |b(t; \eta,v)| \leq c (\norm{\eta}_{H(t)} +h\norm{\eta}_{V(t)} +h^{k+1}\norm{w}_{Z_k(t)}) \norm{v}_{Z_1(t)}.
\end{align}
\end{lemma}
\begin{proof}
We begin with  the following estimate of the bilinear form $b(t;\cdot, \cdot)$, \cref{eq;derive}, for $\eta:= w - \pi^h w$:
\begin{align*}
    |b(t; \eta,v)| & 
    \leq c\norm{\eta}_{H(t)} \norm{v}_{H(t)}
    + \bigg|\sum_{i = 1}^2  \int_{\Omega_i(t)} \mathcal{D}_i^B (\mathbf{w}, \mathcal{B}_i, \eta_i, v_i)
    + \mathcal{D}_i^A(\mathbf{w}, \mathcal{A}_i, \eta_i, v_i) \bigg|.
\end{align*}
For $\mathcal{D}_i^A$, integrating by parts yields: 
%\begin{align*}
%    \bigg| \sum_{i = 1}^2 \int_{\Omega_i(t)} &\mathcal{D}_i^A(\mathbf{w}, \mathcal{A}_i, \eta_i, v_i) \bigg| \leq \bigg| \int_{\Gamma(t)} \bigg \llbracket (\partial_t^\bullet \mathcal{A}_i+ \nabla\cdot \mathbf{w} \mathcal{A}_i) \frac{\partial v}{\partial \nu_\Gamma} \bigg\rrbracket \eta - \sum_{i = 1}^2 \int_{\Omega_i(t)} \eta_i \nabla\cdot ([\partial^\bullet_t \mathcal{A}_i + \nabla \cdot\mathbf{w} \mathcal{A}_i]\nabla v_i) \\
%    & +\sum_{i = 1}^2 \int_{\Omega_i(t)} \eta_i \sum_{l,k = 1}^d \nabla_l(\mathcal{A}_i[\nabla_l \mathbf{w}_k +\nabla_k \mathbf{w}_l]\nabla_k v_i)  - \int_{\Gamma(t)}\eta \sum_{l,k = 1}^d \llbracket \nu_l \mathcal{A}_i[\nabla_l \mathbf{w}_k +\nabla_k \mathbf{w}_l]\nabla_k v_i\rrbracket \bigg|, \\
%    &\leq c\left(|\mathcal{A}_i|_{C^2(\overline{\mathcal{Q}}_i; \mathbb{R})}, |\nabla \mathbf{w}|_{C^1(\overline{\Omega}; \mathbb{R}^d)}\right)  (\norm{\eta}_{H(t)} + \norm{\eta}_{\mathcal{V}^*_{\Gamma}(t)})\norm{v}_{Z_1(t)}.
%\end{align*}
\begin{align*}
 & \left| \sum_{i = 1}^2 \int_{\Omega_i(t)} \mathcal{D}_i^A(\mathbf{w}, \mathcal{A}_i, \eta_i, v_i)  \right| \\
 & \le \left| \int_{\Gamma(t)} \llbracket (\partial_t^\bullet \mathcal{A}_i + \nabla \cdot w \mathcal{A}_i - 2 D(w_i,\mathcal{A}_i)) \nabla v \cdot \nu_\Gamma \rrbracket \eta \right|
 + \left| \sum_{i=1}^2 \eta_i \nabla \cdot \left( (\partial_t^\bullet \mathcal{A}_i + \nabla \cdot w \mathcal{A}_i - 2 D(w_i,\mathcal{A}_i)) \nabla v_i \right) \right| \\
    &\leq c\left(|\mathcal{A}_i|_{C^2(\overline{\mathcal{Q}}_i; \mathbb{R})}, |\nabla \mathbf{w}|_{C^1(\overline{\Omega}; \mathbb{R}^d)}\right)  (\norm{\eta}_{H(t)} + \norm{\eta}_{\mathcal{V}^*_{\Gamma}(t)})\norm{v}_{Z_1(t)}.
\end{align*}
In the last line we have used both the generalised trace inequality, see \cite[Sec.~2.5]{MR0227584}, and the Banach triple identification for the boundary terms. Similarly for $\mathcal{D}_i^B$:
\begin{align*}
    \bigg| \sum_{i = 1}^2 \int_{\Omega_i(t)} &\mathcal{D}_i^\mathcal{B}(\mathbf{w}, \mathcal{B}_i, \eta_i, v_i) \bigg| \leq c\left(|B_i|_{C^1(\overline{\mathcal{Q}}_i; \mathbb{R}^d)}, |\nabla\cdot \mathbf{w}|_{C^1(\overline{\Omega}; \mathbb{R}^d)} \right) (\norm{\eta}_{H(t)} + \norm{\eta}_{\mathcal{V}^*_{\Gamma}(t)})\norm{v}_{Z_1(t)}.
\end{align*}
Combining the previous three estimates we see
\begin{align} \label{b-ineq2}
      |b(t; \eta,v)| & \leq c\big( \norm{\eta}_{H(t)} + \norm{\eta}_{\mathcal{V}_\Gamma^*(t)} \big)\norm{v}_{Z_1(t)}.
\end{align}

In order to complete the proof, we employ  the same duality argument as in \cite{EllRan21,Douglas1972/73} to estimate the dual norm of $\eta$. Set $\mathcal{T}:\mathcal{H}_\Gamma(t) \to \mathcal{V}_\Gamma(t)$ as:
\begin{align*}
    (\mathcal{T}\zeta, v)_{\mathcal{V}_\Gamma(t)}
    := \langle \zeta, v \rangle_{\mathcal{V}_\Gamma(t)}
    = (\zeta,v)_{\mathcal{H}_\Gamma(t)}
    \text{ for all } \zeta, v \in \mathcal{V}_\Gamma(t),
\end{align*}
i.e $\mathcal{T}$ acts a Riesz map mapping to the element in $\mathcal{V}_\Gamma(t)$ that corresponds to the functionals in $\mathcal{H}_\Gamma(t) \subset \mathcal{V}_\Gamma^*(t)$. Notice that:
\begin{align*}
    \norm{\mathcal{T}\zeta}^2_{\mathcal{V}_\Gamma(t)} = \norm{\zeta}^2_{\mathcal{V}^*_\Gamma(t)} = \int_{\Gamma(t)} \zeta \mathcal{T}\zeta,
\end{align*} 
for any $\zeta \in \mathcal{H}_\Gamma(t)$. We note that substituting $v =\mathcal{T}(\zeta)$ in \cref{eq:dual-reg-gamma} gives:
\begin{align}
    \label{eq:trace-reg}
    \norm{\mathcal{R}_{\mathcal{H}_{\Gamma}}(\mathcal{T}\zeta)}_{Z_1(t)} \leq c\norm{\mathcal{T}(\zeta)}_{\mathcal{V}_\Gamma(t)} = c\norm{\zeta}_{\mathcal{V}^*_\Gamma(t)} = c(\zeta, \mathcal{T}\zeta)_{\mathcal{H}_\Gamma(t)}^{1/2}.
\end{align}

By construction:
\begin{multline*}
    \norm{\eta}^2_{\mathcal{V}_\Gamma^*(t)} = \int_{\Gamma(t)} \eta \cdot \mathcal{T}\eta
    = a_\kappa(t;\eta, \mathcal{R}_{\mathcal{H}_{\Gamma}}(\mathcal{T}\eta)) \\
    = a_\kappa(t; \eta, \mathcal{R}_{\mathcal{H}_{\Gamma}}(\mathcal{T}\eta) - I^l [\mathcal{R}_{\mathcal{H}_{\Gamma}}(\mathcal{T}\eta)])
    + a_\kappa(t; \eta, I^l[ \mathcal{R}_{\mathcal{H}_{\Gamma}}(\mathcal{T}\eta)]). 
\end{multline*} 
Then, for the first part, we have:
\begin{align*}
    |a_\kappa(t; \eta, \mathcal{R}_{\mathcal{H}_{\Gamma}}(\mathcal{T}\eta) - I^l [\mathcal{R}_{\mathcal{H}_{\Gamma}}(\mathcal{T}\eta)])|
    & \leq \norm{\eta}_{V(t)} \norm{\mathcal{R}_{\mathcal{H}_{\Gamma}}(\mathcal{T}\eta) - I^l [\mathcal{R}_{\mathcal{H}_{\Gamma}}(\mathcal{T}\eta)]}_{V(t)} \\
    & \leq ch\norm{\eta}_{V(t)} \norm{\mathcal{R}_{\mathcal{H}_{\Gamma}}(\mathcal{T}\eta)}_{Z_1(t)} \leq ch \norm{\eta}_{V(t)}\norm{\eta}_{\mathcal{V}_\Gamma^*(t)}.
\end{align*}
%Hence we have:
%\begin{align}\label{b-ineq}
%    \norm{v}^2_{\mathcal{V}_\Gamma^*(t)} \leq ch \norm{v}_{V(t)}\norm{v}_{\mathcal{V}_\Gamma^*(t)} + a_\kappa(t; v, I^l %[\mathcal{R}_{\mathcal{H}_{\Gamma}}\circ \mathcal{T}(v)]). 
%\end{align}
For the second part, using the definition of the Ritz projection \cref{eq:ritz}, we have
\begin{align*}
  |a_\kappa
  & (t;\eta, I^l[\mathcal{R}_{\mathcal{H}_{\Gamma}}(\mathcal{T}\eta)])| \\
  & = |a_\kappa(t; \pi^h w, I^l [\mathcal{R}_{\mathcal{H}_{\Gamma}}(\mathcal{T}\eta)]) - a_\kappa^h(\Pi^hw, I^h [\mathcal{R}_{\mathcal{H}_{\Gamma}}(\mathcal{T}\eta)]) |\\
  & \leq |a_\kappa(t; \pi^h w - w, I^l [\mathcal{R}_{\mathcal{H}_{\Gamma}}(\mathcal{T}\eta)])
    - a_\kappa^h (t; \Pi^h w - w^{-l}, I^h [\mathcal{R}_{\mathcal{H}_{\Gamma}}(\mathcal{T}\eta)]) | \\
  & \qquad +|a_\kappa(t; w, I^l [\mathcal{R}_{\mathcal{H}_{\Gamma}}(\mathcal{T}\eta)]
    - \mathcal{R}_{\mathcal{H}_{\Gamma}}(\mathcal{T}\eta))
    - a^h_\kappa(t; w^{-l}, I^h [\mathcal{R}_{\mathcal{H}_{\Gamma}}(\mathcal{T}\eta)^{-l} ]
    - \mathcal{R}_{\mathcal{H}_{\Gamma}}(\mathcal{T}\eta)^{-l})|, \\
  & \qquad + |a_\kappa(t; w, \mathcal{R}_{\mathcal{H}_{\Gamma}}(\mathcal{T}\eta))
    - a_\kappa^h(t; w^{-l}, \mathcal{R}_{\mathcal{H}_{\Gamma}}(\mathcal{T}\eta)^{-l}) | \\
    %& \leq c^{h^k} \norm{\pi^h w - w}_{V(t)} \norm{I^l[\mathcal{R}_{\mathcal{H}_{\Gamma}}(\mathcal{T}\eta)]}_{V(t)} \\
    %& + c h^{k+1} \norm{w}_{V(t)} \norm{\mathcal{R}_{\mathcal{H}_{\Gamma}}(\mathcal{T}\eta)}_{Z_1(t)} \\
    %& + c h^{2k+1} \norm{w}_{Z_1(t)} \norm{\mathcal{R}_{\mathcal{H}_{\Gamma}}(\mathcal{T}\eta)}_{Z_1(t)} \\
    & \leq c h^{k + 1} \norm{w}_{Z_1(t)} \norm{\mathcal{R}_{\mathcal{H}_{\Gamma}}(\mathcal{T}\eta)}_{V(t)}
    + c h^{k+1} \norm{w}_{V(t)} \norm{\eta}_{\mathcal{V}_\Gamma^*(t)}
    + c h^{k+1} \norm{w}_{Z_1(t)} \norm{\eta}_{\mathcal{V}_\Gamma^*(t)} \\
    &\leq ch^{k+1} \norm{w}_{Z_1(t)} \norm{\eta}_{\mathcal{V}_\Gamma^*(t)},
\end{align*}
using \cref{lem:interp,lem:ritz}, the regularity estimate, \cref{eq:trace-reg}, and \labelcref{eq:m-error}, \labelcref{eq:a-error} and \labelcref{eq:a-error2}.
Hence, we infer that
\[
\norm{ \eta }_{\mathcal{V}_\Gamma^*(t)}
\le c h \norm{\eta}_{V(t)} + c h^{k+1} \norm{w}_{Z_1(t)}.
\]
Substituting the final inequality here in \cref{b-ineq2} yields \cref{eq:bil}.
\end{proof}
\begin{lemma}[\cite{EllRan21}, Lem.~3.8 and 3.10] \label{lem:ritz}
For $w \in C^1_{Z_k}$, if \labelcref{eq:m-error} to \labelcref{eq:mdV-error}, \cref{lem:interp,lem:adjoint-reg,eq:bil}  hold, there exists $c >0$ (independent of $h$) such that:
\begin{align*}
    \norm{\partial^h_t \Pi^h w}_{V^h(t)} &\leq c(\norm{w}_{V(t)} + \norm{\partial^\bullet_t w}_{V(t)}),\\
    \norm{w - \pi^h w}_{H(t)} +h \norm{w - \pi^hw}_{V(t)} &\leq c h^{k+1} \norm{w}_{Z_k(t)}, \\
    \norm{\partial^l_t ( w - \pi^h w)}_{H(t)} + h \norm{\partial^l_t(w - \pi^hw)}_{V(t)} &\leq c h^{k+1}(\norm{w}_{Z_k(t)} +\norm{\partial_t^\bullet w}_{Z_k(t)}).
\end{align*}
\end{lemma}
\begin{lemma}\label{lemma;ritz-ae}
    The error estimates described in \cref{lem:ritz} also hold for a.e $t \in I$ for $w \in W(Z_k,Z_k)$.
\end{lemma}
\begin{proof}
  Via \cref{lem:equiv-material-derivative}, $C^1_{Z_k}$ is dense within $W(Z_k,Z_k)$. For a $w \in W(Z_k,Z_k)$, take a sequence $w_\epsilon \in C^1_{Z_k}$ such that $w_\epsilon \to w$ in $W(Z_k,Z_k)$. This implies $\norm{w_\epsilon}_{Z_k(t)} \to \norm{w}_{Z_k(t)}$, $\norm{\partial_t^\bullet w_\epsilon}_{Z_k(t)} \to \norm{\partial_t^\bullet w}_{Z_k(t)}$, $\norm{\pi^h w_\epsilon}_{V(t)} \to \norm{\pi^h w}_{V(t)}$, $\norm{\partial_t^l\pi^h w_\epsilon}_{V(t)} \to \norm{\partial_t^l\pi^h w}_{V(t)}$ and $\norm{\partial_t^h \Pi^h w_\epsilon}_{V^h(t)} \to \norm{\partial_t^h \Pi^h w}_{V^h(t)}$ (by continuity of $\Pi^h$ and $\pi^l$, see \cref{rem:lift-ritz-continuity}) in $L^2(I)$. We will only show the first inequality, but the same method can be used to obtained the two remaining ones. For arbitrary $\delta >0$, there exists an $\epsilon_* >0$ such that when $\epsilon > \epsilon_*$, we have:
  \begin{align*}
    \int_I \left( \norm{\partial^h_t \Pi^h(w)}_{V^h(t)} - \norm{\partial_t^h \Pi^h (w_\epsilon)}_{V^h(t)} \right)^2 \,dt &\leq \delta^2, \\
    \int_I \left( \norm{\partial^\bullet_t w}_{V(t)} - \norm{\partial^\bullet_t w_\epsilon}_{V(t)} \right)^2 \, dt &\leq \delta^2,\\
    \int_I \left( \norm{w}_{V(t)} - \norm{ w_\epsilon}_{V(t)} \right)^2 \, dt &\leq \delta^2.
  \end{align*}
  We fix an arbitrary $t \in I$ and by taking the mean integral in $[t-s, t+s]$, $s >0$ (reflecting the functions for $t <0$, i.e $v(t) = v(-t)$), we see:
  \begin{align}\label{eq:mean-ritz-ineq}
      \frac{1}{2s} \int_{t-s}^{t+s} \norm{\partial_t^h \Pi^h (w)}_{V^h(\tau)}\,d\tau \leq  \frac{1}{2s} \int_{t-s}^{t+s} \bigl|\norm{\partial^h_t \Pi^h(w)}_{V^h(\tau)} - \norm{\partial_t^h \Pi^h (w_\epsilon)}_{V^h(\tau)} \bigr| +\norm{\partial_t^h \Pi^h (w_\epsilon)}_{V^h(\tau)}\,d\tau.
  \end{align}
  We note that via Young's inequality, \cref{eq:mean-ritz-ineq} can be bounded by:
  \begin{multline*}
      \frac{1}{2s} \int_{t-s}^{t+s} \bigl|\norm{\partial^h_t \Pi^h(w)}_{V^h(\tau)} - \norm{\partial_t^h \Pi^h (w_\epsilon)}_{V^h(\tau)} \bigr| \,d\tau \\
      \leq s^{-1/2}\left( \int_I \left( \norm{\partial^h_t \Pi^h(w)}_{V^h(t)} - \norm{\partial_t^h \Pi^h (w_\epsilon)}_{V^h(t)} \right)^2 \, d\tau \right)^{1/2} \leq \delta s^{-1/2}.
  \end{multline*}
  Whereas for the second term, we use the bounds given by \cref{lem:ritz} and H\"{o}lder's inequality once more:
  \begin{align*}
      \frac{1}{2s} \int_{t-s}^{t+s} \norm{\partial_t^h \Pi^h (w_\epsilon)}_{V^h(\tau)}\,d\tau &\leq  \frac{c}{2s} \int_{t-s}^{t+s} \norm{w_\epsilon}_{V(\tau)} +\norm{\partial_t^\bullet w_\epsilon}_{V(\tau)}\,d\tau \\
      &\leq \frac{c}{2s} \int_{t-s}^{t+s} \norm{w}_{V(\tau)} +\norm{\partial_t^\bullet w}_{V(\tau)}
      + \bigl|\norm{w_\epsilon}_{V(\tau)} - \norm{w}_{V(\tau)}\bigr|\\
      & \qquad\qquad\qquad +\bigl|\norm{\partial_t^\bullet w_\epsilon}_{V(\tau)} - \norm{\partial_t^\bullet w}_{V(\tau)}\bigr|\,d\tau \\
      &\leq \frac{c}{2s} \int_{t-s}^{t+s} \norm{w}_{V(\tau)} +\norm{\partial_t^\bullet w}_{V(\tau)}\,d\tau +2 \delta.
  \end{align*}
  Combining both of these estimates and equation \cref{eq:mean-ritz-ineq}:
  \begin{align}\label{eq:ritz-ineq-delta}
       \frac{1}{2s} \int_{t-s}^{t+s} \norm{\partial^h_t \Pi^h(w)}_{V^h(\tau)} \,d\tau \leq \frac{c}{2s} \int_{t-s}^{t+s} \norm{w}_{V(\tau)} +\norm{\partial_t^\bullet w}_{V(\tau)}\,d\tau +(2+ s^{-1/2}) \delta.
  \end{align}
  Since $\delta$ is arbitrary, letting $\delta = o(s^{1/2})$, we see that by the limit as $s \to 0$, via the Lebesgue differentiation theorem, both sides of equation \cref{eq:ritz-ineq-delta} converge to their point-wise values a.e, hence we obtain that for all $t \in I$:
  \begin{align*}
      \norm{\partial^h_t \Pi^h(w)}_{V^h(t)} \leq c(\norm{w}_{V(t)} +\norm{\partial_t^\bullet w}_{V(t)}) 
  \end{align*}
  where the constant $c$ is the same as the constant in its equivalent estimate in \cref{lem:ritz} (and hence independent of both $h$ and $t$).
\end{proof}
\subsection{Proof of error bound}
\begin{proof}[Proof of \cref{thm:error}]
We have an additional right-hand side functional term that is not present in the  original proof by \citeauthor{EllRan21}. This requires  a modification of the proof of \cite[Thm.~3.11]{EllRan21}. 

To begin, we slightly modify the problem. For  a test function $v \in W(V,H)$, we can rewrite the weak formulation of our problem as:
\begin{align*}
    \frac{d}{dt} m(t;u,v) - m(t;u,\partial^\bullet_t v) +a(t;u,v) = l(t;v).
\end{align*}
Employing  the standard parabolic rescaling $\check{u} = e^{-\kappa t}u$ where $\kappa$ is chosen  as in the definition of the Ritz projection, the problem becomes:
\begin{align}\label{eq:error}
    \frac{d}{dt}m(t;\check{u}, v) - m(t;\check{u},\partial^\bullet_t v)+ a_\kappa(t;\check{u},v) = \underbrace{e^{-\kappa t}l(v)}_{ =:\widetilde{l}(v)}.
\end{align}
Performing the same transformation to the discrete analogue: define $\check{U}_h = e^{-\kappa t} U_h$ which satisfies
\begin{align}\label{eq:errordis}
      \frac{d}{dt}m^h(t;\check{U}^h, v^h) - m^h(t;\check{U},\partial^h_t v^h)+ a^h_\kappa(t;\check{U}^h,v^h) = \underbrace{e^{-\kappa t}l^h(v^h)}_{ =:\check{l}^h(v^h)}.
\end{align}

Set $\theta :=  \check{U} - \Pi^h \check{u}$, then using \cref{eq:error} and using the fact that $l^h(t;\cdot)$ equals $l(t; (\cdot)^l)$ for functions in $H^h(t)$, for arbitrary $v^h \in W(\mathcal{S}^h, \mathcal{S}^h)$, we arrive at:
\begin{align*}
&\frac{d}{dt}m^h(t; \Pi^h \check{u}, v^h) +a^h_\kappa(t; \Pi^h\check{u}, v^h) - m^h(t; \Pi^h\check{u}, \partial^h_t v^h) - \widetilde{l}^h(v^h), \\
&= \frac{d}{dt}m^h(t; \Pi^h \check{u}, v^h) +a_\kappa(t; \check{u}, v^{h,l}) - m^h(t; \Pi^h\check{u}, \partial^h_t v^h) - \widetilde{l}^h(v^h) ,\\
&= \frac{d}{dt}\bigg[m^h(t; \Pi^h \check{u}, v^h) - m(t; \check{u}, v^{h,l}) \bigg] - \bigg[m^h(t; \Pi^h \check{u}, \partial^h_t v^h) - m(t; \check{u}, \partial^\bullet_t v^{h,l})\bigg],\\
&= m^h(t; \partial_t^h \Pi^h \check{u}, v^h) - m(t; \partial^l_t \check{u}, v^{h,l}) +\lambda^h(t; \Pi^h\check{u}, v^{h,l}) - \lambda^l(t;\check{u}, v^{h,l}) +m(t; \check{u}, \partial^\bullet_t v^{h,l} - \partial^l_t v^{h,l}).
\end{align*}
Now subtracting this equation from \cref{eq:errordis}, and rearranging yields:
\begin{align}\label{eq:erroru}
    &\frac{d}{dt}m^h(t;\theta, v^h) +a^h_\kappa(t; \theta, v^h) - m^h(t; \theta, \partial^h_t v^h) \\ \nonumber
    &= - \bigg[m^h(t; \partial_t^h \Pi^h \check{u}, v^h)  - m(t; \partial^l_t \pi^h\check{u}, v^{h,l}) + m(t; \partial^l_t [\pi^h\check{u}-\check{u}], v^{h,l}) +\lambda^h(t; \Pi^h\check{u}, v^{h,l}) - \lambda^l(t;\pi^h\check{u}, v^{h,l})\\ \nonumber
    &\qquad +\lambda^l(t;[\pi^h\check{u} - \check{u}], v^{h,l}) +m(t; \check{u}, \partial^\bullet_t v^{h,l} - \partial^l_t v^{h,l})\bigg]=:-\Xi^h(\check{u}, v^h).
\end{align}
Using the identity $\partial^l_t (v^{l,h}) = (\partial^h_t v^h)^l$ and looking at $\Xi^h(\cdot, \cdot)$ term by term, we see, for example:
\begin{align*}
    & |m^h(t; \partial_t^h \Pi^h \check{u}, v^h)  - m(t; \partial^l_t \pi^h\check{u}, v^{h,l})| \\
    &= |m^h(t; \partial_t^h \Pi^h \check{u}, v^h) -m(t;(\partial_t^h \Pi^h \check{u})^l, v^{h,l} ) +  m(t; (\partial_t^h \Pi^h\check{u})^l  - \partial^l_t \pi^h\check{u}, v^{h,l})|, \\
    &\leq ch^{k+1} (\norm{\check{u}}_{Z_k(t)} + \norm{\partial^\bullet_t \check{u}}_{Z_k(t)})\norm{v^h}_{V^h(t)},
\end{align*}
by \labelcref{eq:m-error} and \cref{lemma;ritz-ae}. Similar rearrangement and the use of \cref{lem:interp} with \labelcref{eq:c-error}, \labelcref{eq:ct-error} ,\labelcref{eq:md-error} yields:
\begin{align*}
    |\Xi^h(\check{u}, v^h)| \leq ch^{k+1} (\norm{\check{u}}_{Z_k(t)} + \norm{\partial^\bullet_t \check{u}}_{Z_k(t)} ) \norm{v^h}_{V^h(t)}.
\end{align*}
Using \cref{eq:erroru} and substituting $v^h = \theta$, we obtain:
\begin{align*}
    \frac{d}{dt}m^h(t;\theta, \theta) +a^h_\kappa(t; \theta, \theta) - m^h(t; \theta, \partial^h_t \theta) = \Xi^h(\theta,\theta).
\end{align*}
Using the transport formula and the bound on $\Xi^h$:
\begin{align*}
   \frac{1}{2} \frac{d}{dt}m^h(t; \theta, \theta) +a^h_\kappa(t; \theta, \theta) \leq -\frac{1}{2} \lambda^h(t;\theta,\theta) +ch^{k+1} (\norm{\check{u}}_{Z_k(t)} + \norm{\partial^\bullet_t \check{u}}_{Z_k(t)} ) \norm{\theta}_{V^h(t)}.
\end{align*}
Integrating over time and using Young's and Gr\"{o}nwall's inequality:
\begin{align*}
   \sup_{t \in I} \norm{\theta}^2_{H^h(t)} +\int_0^T \norm{\theta}^2_{V^h(t)} \leq c \norm{\theta}^2_{H^h(0)} +c\int_0^T \norm{\theta}^2_{H^h(t)} +ch^{2k+2} \int_0^T (\norm{\check{u}}^2_{Z_k(t)} + \norm{\partial^\bullet_t \check{u}}^2_{Z_k(t)}).
\end{align*}
Finally, using the decomposition:
\begin{align*}
    \check{u} - {\check{U}}^{h,l} = \check{u} - \pi^h\check{u} +\pi^h\check{u} - {\check{U}}^{h,l}. 
\end{align*}
Using the previous bound, the fact that the lift is a diffeomorphism and the bound on the Ritz map, we finally obtain:
\begin{align*}
      \sup_{t \in I} \norm{\check{u} - {\check{U}}^{h,l} }^2_{H(t)} +h^2\int_0^T \norm{\check{u} - {\check{U}}^{h,l} }^2_{V(s)}\,ds &= \sup_{t \in I} \norm{\check{u} - \pi^h\check{u} +\theta^l}^2_{H(t)} +h^2\int_0^T \norm{\check{u}- \pi^h\check{u} +\theta^{l} }^2_{V(s)}\,ds, \\
      &\leq c \norm{u_0 - u_0^{h,l}}_{H(0)} +h^{2k+2}c(C_u).
\end{align*}
Undoing the scaling $u = e^{\kappa t}\check{u}$ gives us the desired error bound.
\end{proof}

%Justification of P1- P9 <-- We assume that the domain is regular enough such that these follow in the same way as in Ranner&Elliott
%Follows by similar method as Ranner&Elliott, section precise. Here show B3. 
%Statement of the error bound
\section{Numerical results}
\label{sec:numerics}

All numerical results are computed using the \textsf{firedrake} package \cite{petsc-user-ref,petsc-efficient,Dalcin2011,Rathgeber2016}.
Simulation code is available in \cite{code}.
Results are computed on a sequence of meshes generated using \textsf{GMSH} \cite{GMSH} rather than successive refinement of a single mesh.

The main challenges in implementing the numerical scheme are:
\begin{enumerate}
\item Computing the initial geometry:
  We start with a piecewise linear geometry given by \textsf{GMSH}.
  The initial isoparametric domain is computed through an explicit parametrisation applying directly the method from \cref{sec:fem} efficiently using custom written \textsf{C} code.
  The evolution of the mesh  is carried out simply by moving the initial Lagrange nodes according to the smooth, given velocity field.
\item Labelling and tracking different parts of the domain:
Alongside the geometry and topology of the mesh we must track labels which say which elements are in domain $\Omega^h_1(t)$ or $\Omega^h_2(t)$ and which facets are on $\Gamma^h(t)$. Once this is fixed for the initial domains $\widetilde{\Omega}^h_1$ and $\widetilde{\Omega}^h_2$, this information is passed between different times.
\textsf{GMSH} provides physical tags to each element and facets which can be used to identify the different domains.
\end{enumerate}

Efficient and accurate quadrature rules are used to perform element-wise integrals.
Note that system matrices must be reassembled at each time step due to the evolution of the domain.

\subsection{Time discretisation of advection-diffusion problem}

We start from the spatial discretisation from \cref{sec:fem-scheme}. We will apply a backward difference formula (BDF) time discretisation of order $q$, see \cite{Kovacs_2016} for further details, including the analysis of a similar surface only problem.
We take a partition of the time interval $0 = t_0 < t_1 < \cdots < t_M = T$.
For simplicity we assume that each time interval is of the same length: $\tau := t_j - t_{j-1}$ for $j = 1, 2, \ldots, M$.

We use temporal interpolations of each domain at each time step to construct a sequence of triangulations  $\mathcal{J}^h(t_j)$ each equipped with finite element spaces $\mathcal{S}^h(t_j)$ for $j=0, 1, \ldots, M$.
We define the discrete velocity $W^j \in \mathcal{S}^h(t_j)^d$ by
\begin{equation}
  W^j = \frac{1}{\tau} \sum_{l=0}^q \delta_l X^{j-l},
\end{equation}
where $X^{j}$ are the positions of the Lagrange nodes of the triangulation at time $t_j$ and $\{ \delta_l \}_{l=0}^q$ are the backward difference formula weights of order $q$, determined from the relation:
\begin{equation}
  \label{eq:bdf-weights}
  \delta(\zeta) = \sum_{l=0}^q \delta_l \zeta^q
  = \sum_{l=1}^q \frac{1}{l} (1 - \zeta)^l.
\end{equation}

The fully discrete problem is the time discretisation \cref{eq:weak-form}:
Given starting values $U^0 \in \mathcal{S}^h(t_0)$, $\ldots$, $U^{q-1} \in \mathcal{S}^h(t_{q-1})$, and data $\mathcal{A}, \mathcal{B}, \mathcal{C}$ and $l_h$, for $j = q, \ldots, M$, we wish to find $U^j \in \mathcal{S}^h(t_j)$ as the solution of
\begin{multline}
  \label{eq:full-discrete}
  \frac{1}{\tau} \sum_{l=0}^q \delta_{l} m^h(t^{j-l}; U^{j-l}, \chi_i^{j-l})
  + a^h(t^j; U^j, \chi^j))
  =  l_h(t^j; \chi^i_j)
  \qquad \mbox{ for all basis function } \chi_i^j \in \mathcal{S}^h(t_j),
\end{multline}
where again $\delta_l$ are the BDF weights \cref{eq:bdf-weights}. Note that the first term on the left hand side is computed by summing over $q$ different meshes to approximate the time derivative. Let $U^{h}(t_j) := U^j$ and let $e^h(t) := u - U^{h,l}(t)$, then we assume that a similar estimate as in \cite[Thm.~5.3]{Kovacs_2016} and \cite[Thm.~2.4]{MR3022237} holds towards the BDF scheme \cref{eq:full-discrete} in addition to \cref{thm:error}:
\begin{align}\label{eq:fully-discrete-scheme}
    \norm{e^h(t_n)}_{L^2(\Omega)}^2 + h^2 \tau \sum_{k = 1}^n\norm{\nabla e^h(t_k)}^2_{L^2(\Omega)} \leq c(\tau^{2q} + h^{2k+2}).
\end{align}
\begin{remark}
\normalfont
We recall that on the space $V(t)$, both the norms $\norm{\cdot}_{L^2(\Omega)}$ and $\norm{\cdot}_{H(t)}$ are equivalent. Moreover, we expect the methods in \cite{Kovacs_2016} and in \cite{MR3022237} to be generalisable in our case in order to prove \cref{eq:fully-discrete-scheme}.
\end{remark}

\subsection{Numerical examples of advection-diffusion problem}

For $d = 2, 3$, let $\Omega = [-1,1]^d$, for $t \in [0,T]$, we define the evolution of the domain through the flow map $\Phi_t$ given by:
\begin{align*}
  \Phi_t(x) = x + \frac{|x|^{1/3}\prod_{i=1}^d(1-x_i^2)}{0.5 \prod_{i=1}^d (1-4 x_i^2 / |x|)}
  \begin{cases}
    ((\alpha(t)-1) x_1, (\beta(t)-1) x_2) & \mbox{ if } d = 2, \\
    ((\alpha(t)-1) x_1, (\beta(t)-1) x_2, 0) & \mbox{ if } d = 3,
  \end{cases}
\end{align*}
for $\alpha(t) = 1 + 0.25 \sin(t)$ and $\beta = 1 + 0.25 \cos(t)$.
This is a special motion which ensures that nodes initially on $\partial \Omega$ do not move and the surface $\Gamma(t)$ is described by the level set function $\phi(\cdot, t)$ given by
\begin{align*}
  \phi(\cdot, t) = \begin{cases}
    \frac{x_1^2}{\alpha(t)^2} + \frac{x_2^2}{\beta(t)^2} - \frac{1}{2}
    & \mbox{ if } d = 2 \\
    \frac{x_1^2}{\alpha(t)^2} + \frac{x_2^2}{\beta(t)^2} + x_3^2 - \frac{1}{2}
    & \mbox{ if } d = 3.
   \end{cases}
\end{align*}

We define $\Omega_1(t)$ as the interior of $\Gamma(t)$ and $\Omega_2(t) = \Omega \setminus \Omega_1(t)$.

We set the coefficients in the equation to be $\mathcal{A}_1 = 10 \mathrm{Id}$, $\mathcal{A}_2 = \mathrm{Id}$,
$\mathcal{B}_1 = 5 \nabla x_1$, $\mathcal{B}_2 = -5 \nabla x_1$,
$\mathcal{C}_1 = 1$, $\mathcal{C}_2 = 10$ and note that they jump across the interface. 
We set the right hand side data such that the exact solution $u$ is given by
\[
  u(x, t) =
  \sin(t) |\Phi(x)| \prod_{i=1}^d \sin(2 \pi x_i).
\]
This exact solution is globally continuous, smooth in each domain but is not differentiable across the interface.
In order to simplify the implementation the right hand side data ($l_h$) is computed by taking interpolations of smooth data.
We compute using isoparametric elements of order 1, 2, 3 on a sequence of given meshes.
For order $k$ discretisation in space we use BDF order $k+1$ in time.
Between refinement levels, we approximately reduce the mesh size by half and reduce the time step by exactly half so that the spatial and temporal errors in \cref{eq:fully-discrete-scheme} balance.
The initial solution $U^0 = 0$ matches the exact solution at $t = 0$. The other starting values are computed using lower order BDF methods.
For elements of order $k$ we expect convergence of order $k+1$ for the error at the final time, $u(T) - U^M$, in the $L^2(\Omega)$ norm and order $k$ in the $H^1(\Omega)$ semi-norm.

\begin{figure}
  \centering
  \includegraphics{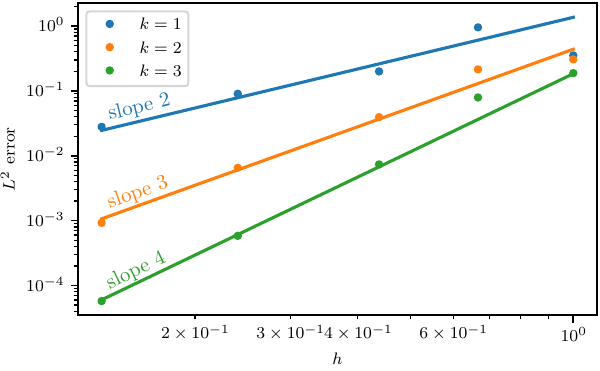}
  \includegraphics{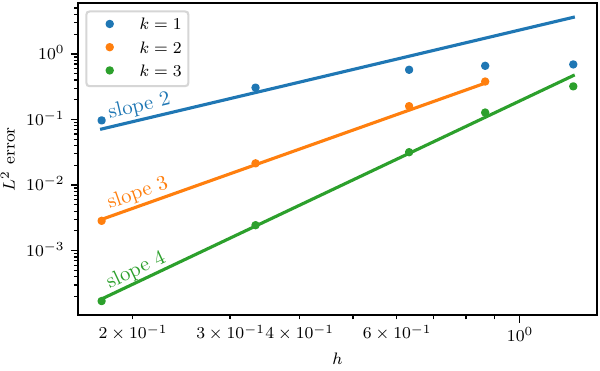}
  \caption{$L^2$ error for advection-diffusion problem for $d=2$ (top) and $d=3$ (bottom).
  Slope lines are reference lines indicating different rates of convergence.}
  \label{fig:results-advection-diffusion}
\end{figure}

The results are shown in \cref{fig:results-advection-diffusion} for the cases $d=2,3$ respectively.
The precise numerical values are shown in \cref{tab:advection-diffusion2,tab:advection-diffusion3}.
Since the meshes are generated independently in \textsf{GMSH} rather than by repeated uniform refinement of a single mesh, we do not expect perfect convergence rates for these experiments in \cref{tab:advection-diffusion2,tab:advection-diffusion3}.
However, the slope lines in \cref{fig:results-advection-diffusion} are consistent with our results and indicate the expected rates of convergence.

\begin{table}
  \begin{subtable}{\textwidth}
  \centering
  \begin{tabular}{cccccc}
\toprule
                     $h$ &                   $\tau$ &              $L^2$ error &              $H^1$ error & eoc($L^2$ error) & eoc($H^1$ error) \\
\midrule
               $1.00000$ &                $1.00000$ & $3.51296 \times 10^{-1}$ &                $3.30187$ &              --- &              --- \\
$6.66667 \times 10^{-1}$ & $5.00000 \times 10^{-1}$ & $9.51874 \times 10^{-1}$ &                $4.43031$ &        -2.458416 &        -0.725039 \\
$4.37747 \times 10^{-1}$ & $2.50000 \times 10^{-1}$ & $1.99238 \times 10^{-1}$ &                $2.33192$ &         3.717892 &         1.525683 \\
$2.40008 \times 10^{-1}$ & $1.25000 \times 10^{-1}$ & $9.02047 \times 10^{-2}$ &                $1.44891$ &         1.318576 &         0.791859 \\
$1.34513 \times 10^{-1}$ & $6.25000 \times 10^{-2}$ & $2.78846 \times 10^{-2}$ & $7.80431 \times 10^{-1}$ &         2.027605 &         1.068576 \\
\bottomrule
\end{tabular}

  \caption{Order 1, $d=2$}
\end{subtable}

\begin{subtable}{\textwidth}
  \centering
  \begin{tabular}{cccccc}
\toprule
                     $h$ &                   $\tau$ &              $L^2$ error &              $H^1$ error & eoc($L^2$ error) & eoc($H^1$ error) \\
\midrule
               $1.00000$ &                $1.00000$ & $3.05166 \times 10^{-1}$ &                $3.29594$ &              --- &              --- \\
$6.66667 \times 10^{-1}$ & $5.00000 \times 10^{-1}$ & $2.14128 \times 10^{-1}$ &                $2.03484$ &         0.873768 &         1.189433 \\
$4.37747 \times 10^{-1}$ & $2.50000 \times 10^{-1}$ & $3.93009 \times 10^{-2}$ & $8.61571 \times 10^{-1}$ &         4.030255 &         2.043065 \\
$2.40008 \times 10^{-1}$ & $1.25000 \times 10^{-1}$ & $6.51346 \times 10^{-3}$ & $2.89276 \times 10^{-1}$ &         2.990807 &         1.816030 \\
$1.34513 \times 10^{-1}$ & $6.25000 \times 10^{-2}$ & $9.24067 \times 10^{-4}$ & $8.43187 \times 10^{-2}$ &         3.372712 &         2.129106 \\
\bottomrule
\end{tabular}

  \caption{Order 2, $d=2$}
\end{subtable}

\begin{subtable}{\textwidth}
  \centering
  \begin{tabular}{cccccc}
\toprule
                     $h$ &                   $\tau$ &              $L^2$ error &              $H^1$ error & eoc($L^2$ error) & eoc($H^1$ error) \\
\midrule
               $1.00000$ &                $1.00000$ & $1.87468 \times 10^{-1}$ &                $2.57196$ &              --- &              --- \\
$6.66667 \times 10^{-1}$ & $5.00000 \times 10^{-1}$ & $7.90234 \times 10^{-2}$ & $8.93941 \times 10^{-1}$ &         2.130551 &         2.606354 \\
$4.37747 \times 10^{-1}$ & $2.50000 \times 10^{-1}$ & $7.35252 \times 10^{-3}$ & $1.97059 \times 10^{-1}$ &         5.645317 &         3.594764 \\
$2.40008 \times 10^{-1}$ & $1.25000 \times 10^{-1}$ & $5.83518 \times 10^{-4}$ & $4.39236 \times 10^{-2}$ &         4.216077 &         2.497725 \\
$1.34513 \times 10^{-1}$ & $6.25000 \times 10^{-2}$ & $5.75918 \times 10^{-5}$ & $2.73132 \times 10^{-2}$ &         3.999388 &         0.820506 \\
\bottomrule
\end{tabular}

  \caption{Order 3, $d=2$}
\end{subtable}
\caption{Results for advection-diffusion problem for $d=2$.}
\label{tab:advection-diffusion2}
\end{table}

\begin{table}
\begin{subtable}{\textwidth}
  \centering
  \begin{tabular}{cccccc}
\toprule
                     $h$ &                   $\tau$ &              $L^2$ error & $H^1$ error & eoc($L^2$ error) & eoc($H^1$ error) \\
\midrule
               $1.25000$ &                $1.00000$ & $6.90651 \times 10^{-1}$ &   $8.44875$ &              --- &              --- \\
$8.66599 \times 10^{-1}$ & $5.00000 \times 10^{-1}$ & $6.57464 \times 10^{-1}$ &   $7.69389$ &         0.134429 &         0.255490 \\
$6.31590 \times 10^{-1}$ & $2.50000 \times 10^{-1}$ & $5.72385 \times 10^{-1}$ &   $6.68398$ &         0.438071 &         0.444822 \\
$3.33531 \times 10^{-1}$ & $1.25000 \times 10^{-1}$ & $3.05710 \times 10^{-1}$ &   $4.43667$ &         0.982256 &         0.641827 \\
$1.75870 \times 10^{-1}$ & $6.25000 \times 10^{-2}$ & $9.68233 \times 10^{-2}$ &   $2.34095$ &         1.796513 &         0.998992 \\
\bottomrule
\end{tabular}

  \caption{Order 1, $d=3$}
\end{subtable}

\begin{subtable}{\textwidth}
  \centering
  \begin{tabular}{cccccc}
\toprule
                     $h$ &                   $\tau$ &              $L^2$ error &              $H^1$ error & eoc($L^2$ error) & eoc($H^1$ error) \\
\midrule
$8.66599 \times 10^{-1}$ & $5.00000 \times 10^{-1}$ & $3.78772 \times 10^{-1}$ &                $5.39342$ &              --- &              --- \\
$6.31590 \times 10^{-1}$ & $2.50000 \times 10^{-1}$ & $1.59518 \times 10^{-1}$ &                $2.89154$ &         2.733729 &         1.970655 \\
$3.33531 \times 10^{-1}$ & $1.25000 \times 10^{-1}$ & $2.14372 \times 10^{-2}$ & $9.12177 \times 10^{-1}$ &         3.143326 &         1.806895 \\
$1.75870 \times 10^{-1}$ & $6.25000 \times 10^{-2}$ & $2.84564 \times 10^{-3}$ & $2.57151 \times 10^{-1}$ &         3.155269 &         1.978425 \\
\bottomrule
\end{tabular}

  \caption{Order 2, $d=3$}
\end{subtable}

\begin{subtable}{\textwidth}
  \centering
  \begin{tabular}{cccccc}
\toprule
                     $h$ &                   $\tau$ &              $L^2$ error &              $H^1$ error & eoc($L^2$ error) & eoc($H^1$ error) \\
\midrule
               $1.25000$ &                $1.00000$ & $3.19819 \times 10^{-1}$ &                $4.78209$ &              --- &              --- \\
$8.66599 \times 10^{-1}$ & $5.00000 \times 10^{-1}$ & $1.27429 \times 10^{-1}$ &                $2.73209$ &         2.511984 &         1.528193 \\
$6.31590 \times 10^{-1}$ & $2.50000 \times 10^{-1}$ & $3.16946 \times 10^{-2}$ & $9.05805 \times 10^{-1}$ &         4.398523 &         3.489946 \\
$3.33531 \times 10^{-1}$ & $1.25000 \times 10^{-1}$ & $2.44323 \times 10^{-3}$ & $1.41707 \times 10^{-1}$ &         4.013792 &         2.905320 \\
$1.75870 \times 10^{-1}$ & $6.25000 \times 10^{-2}$ & $1.70045 \times 10^{-4}$ & $2.10221 \times 10^{-2}$ &         4.164150 &         2.981596 \\
\bottomrule
\end{tabular}

  \caption{Order 3, $d=3$}
\end{subtable}
\caption{Results for advection-diffusion problem for $d=3$.}
\label{tab:advection-diffusion3}
\end{table}

\section*{Author contributions}

Conceptualization, investigation, formal analysis, writing: All.
Software: TR.
Supervision: CME.
Writing (first draft): PS.

\section*{Acknowledgements}

We thank the reviewers for their comments, which have improved the paper.

\pagebreak
\begin{appendix}\label{sec:appendix}

\section{Proof of Regularity}

In this section, we will show some results on the additional regularity of the smooth solution to \cref{eq:weak-form}. 
\begin{lemma}[The Trace Map]
There exists a bounded and continuous linear operator $\widetilde{\tau}_{(\cdot)}: L^2_V \to L^2_{\mathcal{V}_\Gamma}$ such that $\widetilde{\tau}_{t} \varphi(t) = \tau_t \varphi(t)$ $\forall \varphi \in C_V$, where $\tau_t: V(t) \to \mathcal{V}_\Gamma(t)$ is the classical trace map.
\end{lemma}
\begin{proof}
Let $\tau_t:V(t) \to \mathcal{V}_\Gamma(t)$ be the classical trace map. It is proven in \cite{AlpEllSti15b} that the following identity: $\tau_t (\phi_t w_0) = \phi_t (\tau_0 w_0)$ holds for all $ t \in I$ and $w_0 \in V(0)$. Moreover, there exists a $c$ independent of time such that:
\begin{align*}
    \norm{\tau_t (\phi_t w_0)}_{\mathcal{H}_\Gamma(t)} \leq c \norm{\phi_t w_0}_{\mathcal{V}_\Gamma(t)}.
\end{align*}
Now, formally define $\widetilde{\tau}_{(\cdot)}$ as: 
\begin{align*}
    \widetilde{\tau}_{(\cdot)} v(\cdot)= \phi_{(\cdot)} \tau_{0}( \phi_{-(\cdot)} v(\cdot)).
\end{align*}
Then via \cref{lem:isomo}, $\phi_{-(\cdot)} v(\cdot) \in L^2(I;V(0))$ and since $\tau_0$ can also furthermore be uniquely identified as a linear map $\tau_0: L^2(I;V(0)) \to L^2(I;\mathcal{V}_\Gamma(0))$, see \cite[Thm.~1.2.4]{MR3617205}, finally the push-forward $\phi_t: L^2(I;\mathcal{V}_\Gamma(0)) \to L^2_{\mathcal{V}_\Gamma}$ maps back into the evolving space. Note that this map, by compatibility and it's time independent bound is also bounded. Finally, if $\varphi \in C_{V}$, then, at time $t \in I$:
\begin{align*}
     \widetilde{\tau}_{t} \varphi(t) = \phi_{t} \tau_{0}( \phi_{-t} \varphi(t)) = \tau_t \varphi(t). 
    \quad \qedhere
\end{align*}
\end{proof}
This allows us to formally identify the following pairing:
\begin{align*}
    \int_0^T \langle g, v\rangle_{\mathcal{V}_\Gamma(t)} := \int_0^T  \langle g, \widetilde{\tau}_t v\rangle_{\mathcal{V}_\Gamma(t)}, \; (g,v) \in L^2_{\mathcal{V}^*_\Gamma} \times L^2_{V}.
\end{align*}
\begin{lemma}\label{lem:A_1}Under the assumptions of \ref{A4}  \cref{theorem:existence}, for each $g \in L^2_{\mathcal{V}_\Gamma}$, there exists a unique solution $u_g \in L^2_{Z_1}$ to:
\begin{align}
\label{eq:A_1}
    \sum_{i = 1}^2   \int_0^T  \int_{\Omega_i(t)}  \mathcal{A}_i(t;x) \nabla u_g \cdot \nabla v=  \int_0^T \int_{\Gamma(t)} g v
    \quad \text{for all} \quad v \in L^2_{V}.
\end{align}
\end{lemma}
\begin{proof}
It follows from the regularity assumptions on the flow $\mathbf{\Phi}_t$ that the pair $(Z_1(t), \phi_t|_{Z_1(t)})|_{t \in I}$ is compatible.
%We first show that there exists a $u_g \in L^2_{Z_1}$ such that:
%\begin{align*}
%    \sum_{i = 1}^2 \int_{\Omega_i(t)}  \mathcal{A}_i(t;x) \nabla u_g \cdot \nabla v= \int_{\Gamma(t)} g v\quad \text{for almost all} \; t \in I,\; \forall v \in L^2_V,\; \norm{u_g}_{L^2_{Z_1}} \leq C\norm{g}_{L^2_{\mathcal{V}_\Gamma}}.
%\end{align*}
%and show that this implies that $u_g$ solves \cref{eq:A_1}.
For $g \in L^2_{\mathcal{V}_\Gamma}$, we can take a subset $I' \subset I$ of full measure such that $g(t) \in \mathcal{V}_\Gamma(t)$ and $\norm{g(t)}_{\mathcal{V}_\Gamma(t)} <\infty$, for $t \in I'$. Indeed via \cref{lem:isomo}, we can take the set of Lebesgue points of $\widetilde{g}(\cdot) = \phi_{-(\cdot)}g(\cdot)$ to be $I'$ and push the function forwards, $g(t) = \phi_t \widetilde{g}(t) \in \mathcal{V}_\Gamma(t)$ for all $t \in I'$. Moreover, we see that, for $t \in I'$, these remain Lebesgue Points of $\norm{g(\cdot)}_{\mathcal{V}_\Gamma(\cdot)}$ in $L^2(I)$, indeed, via the reverse triangle inequality and compatibility:
\begin{align*}
    \frac{1}{2\delta}\int_{t - \delta}^{t +\delta} |\norm{g}_{\mathcal{V}_\Gamma(s)} - \norm{g}_{\mathcal{V}_\Gamma(t)}|^2\,ds &\leq \frac{c}{2\delta}\int_{t - \delta}^{t +\delta} |\norm{\widetilde{g}(s)}_{\mathcal{V}_\Gamma(0)} - \norm{\widetilde{g}(t)}_{\mathcal{V}_\Gamma(0)} |^2\,ds \\
    &\leq \frac{c}{2\delta}\int_{t - \delta}^{t +\delta} \norm{\widetilde{g}(s) - \widetilde{g}(t)}_{\mathcal{V}_\Gamma(0)}^2\,ds
\end{align*} 
Fix $t \in I'$, set $u^t_g \in V(t)$ to be the solution to:
\begin{align} \label{eq:regular}
\sum_{i = 1}^2 \int_{\Omega_i(t)}  \mathcal{A}(t;x) \nabla u^t_g \cdot \nabla v  = \int_{\Gamma(t)} g v ~~\forall ~v \in V(t).
\end{align}
We will drop the distinction between $I'$ and $I$ and just say for almost all $t \in I$, then via the Hilbert triple structure outlined in \cref{sec:realisation}, $g(t) \in \mathcal{V}_\Gamma(t) \subset \mathcal{H}_\Gamma(t) \subset \mathcal{V}^*_\Gamma(t)$. By \cite[Thm.~1]{nistor} we have that for almost all $t \in I$ there exists a unique solution $u^t_g( \cdot) \in Z_1(t)$ to \cref{eq:regular}. 

Set $u_g(t; \cdot) := u^t_g(\cdot)$, we show that this solution is in-fact in $L^2_{Z_1}$. We will proceed as follows:
\begin{enumerate}
    \item First show that $u_g \in L^2_V$. By \cite[Lem.~2.14]{AlpEllSti15a}, it suffices to show first that $t \to (u_g, w)_{V(t)}$  is measurable for all $w \in L^2_V$ and then that $\norm{u_g}_{L^2_V} < \infty$.
    \item We then reuse this method, showing that $t \to (u_g, w)_{Z_1(t)}$ is measurable for all $w \in L^2_{Z_1}$ and $\norm{u_g}_{L^2_{Z_1}} <\infty$, and hence $u_g \in L^2_{Z_1}$.
    \item Finally, we show that $u_g$ does indeed solve \cref{eq:A_1}.
\end{enumerate}
To show the measurability, since the eigenvalues of $\mathcal{A}_i(t;x)$ are bounded from both below and above independent of time, we can induce the equivalent inner product $(u,v)_{\widetilde{V}(t)} := (\mathcal{A}(t;x)\nabla u, \nabla v)_{H(t)}$. Showing measurability then follows as:
\begin{align*}
 (u_g,v)_{\widetilde{V}(t)} = \sum_{i = 1}^2 \int_{\Omega_i(t)}  \mathcal{A}_i(t;x) \nabla u_g \cdot \nabla v  = \int_{\Gamma(t)} g v = \int_{\Gamma(t)}  \langle g, v\rangle_{\mathcal{V}_\Gamma(t)},
\end{align*}
and since $(v,g) \in L^2_V \times L^2_{\mathcal{V}_\Gamma^*}$, by \cite[Lem.~2.14]{AlpEllSti15a}, the map $t\mapsto \langle g, v\rangle_{\mathcal{V}_\Gamma(t)}$ is measurable and hence so is $(u_g,v)_{\widetilde{V}(t)}$. For the uniform bound, testing the differential equation \cref{eq:regular} with $v = u_g$ and integrating in time, we have:
\begin{align*}
    \norm{u_g}^2_{L^2_V} \leq C(\gamma) \norm{g}^2_{L^2_{\mathcal{V}_\Gamma^*}},
\end{align*}
via Young's and Poincar\'e's inequalities, so $u_g \in L^2_{V}$. Before moving on, note that for any fixed $t \in I$, since $u_g(t; \cdot) \in Z_1(t)$, we can integrate by parts \cref{eq:regular}, obtaining:
\begin{align}\label{eq:int-by-parts}
    \int_{\Gamma(t)}\bigg\llbracket \mathcal{A}(t;x) \nabla u \cdot \nu_\Gamma \bigg\rrbracket  v-\sum_{i = 1}^2 \int_{\Omega_i(t)} \nabla \cdot (\mathcal{A}_i(t;x) \nabla u_g) v = \int_{\Gamma(t)} g v,
\end{align}
for $v \in C^\infty_0(\Omega)$. For $v \in C^\infty_0(\Omega_1(t)) \times C^\infty_0(\Omega_2(t))$, we see that \cref{eq:int-by-parts} yields:
\begin{align*}
    \sum_{i = 1}^2 \int_{\Omega_i(t)} \nabla \cdot (\mathcal{A}_i(t;x) \nabla u_g) v = 0.
\end{align*}
Since the space $C^\infty_0(\Omega_1(t)) \times C^\infty_0(\Omega_2(t))$ is dense in $H(t)$, we get that:
\begin{align}\label{eq:regular-inner-prod}
   \nabla \cdot (\mathcal{A}_i(t;x) \nabla u_g)  = 0 \; \text{a.e}.
\end{align}
By the Poincar\'{e}'s inequality, we can endow $Z_1(t)$ with a more convenient equivalent inner product:
\begin{align*}
    (v,w)_{Z_1(t)} = \sum_{i = 1}^2 \int_{\Omega_i(t)} \Delta v \cdot \Delta w + \nabla v \cdot \nabla w.
\end{align*}
Since $\mathcal{A}(t;x)$ is assumed to be differentiable, we introduce the equivalent inner product on $Z_1(t)$:
\begin{multline*}
  (\eta,v)_{\widetilde{Z}_1(t)} := \\
  \sum_{i = 1}^2 \int_{\Omega_i(t)} \nabla \cdot (\mathcal{A}_i(t;x) \nabla \eta_i)  \nabla \cdot (\mathcal{A}_i(t;x) \nabla v_i) + (\nabla \mathcal{A}_i(t;x))\cdot \nabla \eta_i(\nabla \mathcal{A}_i(t;x))\cdot \nabla v_i + \nabla \eta_i \cdot \nabla v_i.
\end{multline*}
We will first show the following statement, let $\nabla^2$ be the Hessian, then:
\begin{align*}
    \int_{\Omega_i(t)} | \mathcal{A}_i \nabla^2\eta_i|^2 \geq \gamma \int_{\Omega_i(t)} |\Delta \eta_i|^2,
\end{align*}
where $\gamma >0$ is the coercivity constant in \cref{theorem:existence}. For fixed $t \in I$, since $\mathcal{A}_i(t;x)$ is non-singular and symmetric, there exists orthogonal matrices $\mathcal{P}_i(t;x)$, $\mathcal{P}_i^T(t;x)$ and diagonal matrix $\mathcal{D}_i(t;x) = [\lambda_k(t;x) \delta_{k,l}]_{k,l =1}^{n}$, the eigenvalues of $\mathcal{A}_i(t;x)$ (note that eigenvalues are continuous for a continuous matrix) such that $\mathcal{A}_i(t;x) = \mathcal{P}_i(t;x) \mathcal{D}_i(t;x)\mathcal{P}^T_i(t;x)$. Doing a change of coordinates, $x = \mathcal{P}_iy$ and letting $\tilde{\eta}_i(x) = \eta_i(\mathcal{P}_i (x))$, we see:
\begin{align*}
   \mathcal{A}_i(t;x) \nabla^2(\eta_i(x)) = \mathcal{D}_i(t;x) \Delta \tilde{\eta}_i(y)|_{y = \mathcal{P}_i^T x}.
\end{align*}
Hence, since $\mathcal{P}$ is orthogonal:
\begin{align*}
    \int_{\Omega_i(t)}| \mathcal{A}_i(t;x) \nabla^2(\eta_i(x))|^2\,dx &= \int_{\mathcal{P}_i^T(\Omega_i(t))} |\mathcal{D}_i(t;\mathcal{P}_i(y)) \Delta \tilde{\eta}_i(y)|^2\,dy,\\
     &\geq \min_{x \in \Omega_i(t)}\min_{k \in [1,n]}  \lambda^2_k(t;\mathcal{P}_i^T(x)) \int_{\mathcal{P}_i^T(\Omega_i(t))}| \Delta \tilde{\eta}_i(y)|^2\,dy, \\
    &\geq \gamma^2 \int_{\mathcal{P}^T(\Omega_i(t))}| \Delta \tilde{\eta}_i(y)|^2\,dy = \gamma^2  \int_{\Omega_i(t)}| \Delta \eta_i(x)|^2\,dx.
\end{align*}
To show the equivalence of inner products $(\cdot, \cdot)_{\widetilde{Z}_1(t)}$ and $(\cdot, \cdot)_{Z_1(t)}$, note:
\begin{align*}
    \norm{\eta}^2_{\widetilde{Z}_1(t)} &= \sum_{i = 1}^2 \int_{\Omega_i(t)} |\nabla \cdot (\mathcal{A}_i(t;x) \nabla \eta_i)|^2 + |(\nabla \mathcal{A}_i(t;x))\cdot \nabla \eta_i|^2 +|\nabla \eta_i|^2\\
    &= \sum_{i = 1}^2 \int_{\Omega_i(t)}|\mathcal{A}_i(t;x) \nabla^2 \eta_i|^2+ 2\mathcal{A}_i(t;x) \nabla^2 \eta_i  (\nabla  \mathcal{A}_i(t;x) )\cdot\nabla \eta_i +  2|(\nabla  \mathcal{A}_i(t;x) )\cdot \nabla \eta_i|^2  +|\nabla \eta_i|^2\\
    &\geq \sum_{i = 1}^2 \int_{\Omega_i(t)}\tfrac{1}{2}|\mathcal{A}_i(t;x) \nabla^2 \eta_i|^2+ |\nabla \eta_i|^2 \geq \min\left\{\tfrac{1}{2}\gamma^2, 1\right\}\norm{\eta}^2_{Z_1(t)}
\end{align*}
and:
\begin{align*}
\norm{\eta}^2_{\widetilde{Z}_1(t)} \leq c\left[|\mathcal{A}_i|_{C^1(\Omega_i(t))}\right](\norm{\nabla^2 \eta_i}^2_{H(t)} +\norm{\nabla \eta_i}^2_{H(t)}) \leq c\norm{\eta}^2_{Z_1(t)}.
\end{align*}
Hence, $(\cdot, \cdot)_{\tilde{Z}_1(t)}$ and $(\cdot, \cdot)_{Z(t)}$ are equivalent. Substituting $\eta = u_g$ in the new inner product, using \cref{eq:regular-inner-prod}, we arrive at:
\begin{align*}
    (u_g,v)_{\widetilde{Z}_1(t)}  = \sum_{i = 1}^2 \int_{\Omega_i(t)} (\nabla \mathcal{A}_i(t;x))\cdot \nabla u_i(\nabla \mathcal{A}_i(t;x))\cdot \nabla v_i +\nabla u_i \cdot \nabla v_i.
\end{align*}
Since we already know that $u_g \in L^2_V$, both  $((\nabla \mathcal{A}(t;x))\cdot \nabla u,(\nabla \mathcal{A}(t;x))\cdot \nabla v)_{H(t)}$ and $(\nabla u , \nabla v)_{H(t)}$ are measurable and hence the map $t \to (u_g,v)_{\widetilde{Z}_1(t)}$ is measurable for all $v \in Z_1(t)$. For the bound, \cite[Thm.~1]{nistor} gives us a constant $C_t$ (that depends on time) such that:
\begin{align*}
    \norm{u_g}_{Z_1(t)} \leq C_t \norm{g}_{\mathcal{V}_\Gamma(t)}.
\end{align*}
Using a similar method as \cref{lem:pf-compat} and changing the variables, it follows that there exists $C_T> C_t$, $C_T <\infty$. Hence:
\begin{align*}
    \int_0^T  \norm{u_g}_{Z_1(t)}^2 \leq C_T^2 \int_0^T \norm{g}_{\mathcal{V}_\Gamma(t)}^2 <\infty.
\end{align*}
Hence $u_g \in L^2_{Z_1}$.

Finally, we have a solution to $u_g \in L^2_{Z_1}$ solving \cref{eq:regular} a.e in time. Since a function $v \in L^2_V$ satisfies $v(t) \in V(t)$ a.e in time, testing with such a $v(t)$ and integrating both sides of \cref{eq:regular} yields the desired solution. Uniqueness follows via coercivity.
\end{proof}

\begin{lemma}
    Under the assumption \ref{A4} from \cref{theorem:existence}, the solution of \cref{eq:A_1} possesses a weak material derivative $\partial_t^\bullet u_g \in L^2_V$ and satisfies the estimate:
    \begin{align*}
        \norm{u_g}^2_{W(Z_1,V)} \leq c \norm{g}^2_{W(\mathcal{V}_\Gamma, \mathcal{V}_\Gamma^*)}
    \end{align*}
\end{lemma}
\begin{proof}
We use the same method as in the proof of \cref{lem:ritz-additional-regular}. Let:
\begin{align*}
    k(t; v,w) := \sum_{i = 1}^2\int_{\Omega_i(t)} \mathcal{A}_i(t;x) \nabla v \cdot \nabla w
\end{align*}
then we see, going back to \cref{eq:regular}, we see that $u_g$ solves:
\begin{align*}
    k(t;u_g,w) = (g,w)_{\mathcal{H}_\Gamma(t)},
\end{align*}
a.e for all $w \in L^2_V$ (as before, we can take a subset of $I$ of full measure such that $w(t) \in V(t)$, moreover, via the classical trace theorem, we can identify $w(t) \in \mathcal{H}_\Gamma(t)$ by $(g, w)_{\mathcal{H}_\Gamma(t)} = (g, \tau_t w)_{\mathcal{H}_\Gamma(t)}$). The derivative of $k(t;\cdot,\cdot)$, $\dot{k}(t;\cdot,\cdot)$, can be explicitly calculated to be:
\begin{align}\label{eq:k-derivative}
    \frac{d}{dt} k(t; w,v) & = \dot{k}(t; w,v) + k(t;\partial_t^\bullet w, v) +k(t; w, \partial_t^\bullet v), \; \forall v,w \in W(V,V),
\end{align}
where
\begin{align*}
        \dot{k}(t;w,v) = \sum_{i = 1}^2 \int_{\Omega_i(t)} \mathcal{D}^\mathcal{A}_i(\mathbf{w}, \mathcal{A}_i, w_i,v_i), \nonumber
\end{align*}
and $\mathcal{D}^\mathcal{A}_i$ was defined in \cref{eq:diff}. We set $\tilde{u}_g$ to be the solution to:
\begin{align}\label{eq:regularity}
    \int_0^T  k(t; \tilde{u}_g, \eta)\,dt = \int_0^T \lambda_\Gamma(t; g, \eta)- \dot{k}(t;u_g, \eta) +\langle \partial_t^\bullet g, \eta \rangle_{\mathcal{V}_\Gamma(t)}\,dt \quad \forall \eta \in L^2_V.
\end{align}
The material derivative taken on the function $\partial_t^\bullet g$ is the one with the triple $L^2_{\mathcal{V}_\Gamma} \subset L^2_{\mathcal{H}_\Gamma} \subset L^2_{V^*_{\Gamma}}$ and the bilinear form $\lambda_\Gamma$ is the corresponding form from \cref{def:lambda} satisfying the equation:
\begin{align*}
    \frac{d}{dt}(v,w)_{\mathcal{H}_\Gamma(t)} = \langle \partial_t^\bullet v, w\rangle_{\mathcal{V}_\Gamma(t)} + \langle \partial_t^\bullet w,v \rangle_{\mathcal{V}_\Gamma(t)} +\lambda_\Gamma(t;v,w).
\end{align*}
See \cite[Sec.~5.4]{AlpEllSti15b} for an explicit form of $\lambda_\Gamma(t;\cdot,\cdot)$. Using the same method as the proof of \cref{lem:A_1}, we have that if $\partial^\bullet_tg \in L^2_{\mathcal{V}_\Gamma^*}$, there exists a unique $\widetilde{u}_g \in L^2_V$ solving equation \cref{eq:regularity}.

Let $u_g$ be the solution of \cref{eq:regular}. Via isomorphism, $\phi_{-(\cdot)} \tilde{u}_g(\cdot), \phi_{-(\cdot)} u_g(\cdot) \in L^2(I;V(0))$, we pick a Lebesgue point $s_* \in I$ of $\phi_{-(\cdot)} u_g(\cdot)$ and set
\begin{align*}
        w := \phi_t \int_{s_*}^t \phi_{-\tau} \widetilde{u}_g(\tau) \,d\tau+\phi_t \phi_{-s_*} u_g(s_*) =  \phi_t \int_{s_*}^t \phi_{-\tau} \widetilde{u}_g(\tau) \,d\tau+z_{s_*}.
\end{align*}
Thus $w \in W(V,V)$ and $\partial_t^\bullet w = \tilde{u}_g$. We aim to show that $w  = u_g$. To do so, note that by definition of $\dot{k}(t;\cdot,\cdot)$ \cref{eq:k-derivative}, testing  with $\eta \in W(V,V)$ with $\eta(0) = \eta(T) = 0$:
\begin{align}\label{eq:k-w-equation}
    \int_0^T k(t; \partial_t^\bullet w, \eta)\,dt = \int_0^T -\dot{k}(t; w,\eta) - k(t; w, \partial_t^\bullet \eta)\,dt.
\end{align}
we see from comparing \cref{eq:k-w-equation,eq:regularity}:
\begin{align}\label{eq:w-function-equation}
     -\int_0^T  k(t;  w, \partial_t^\bullet\eta)\,dt = \int_0^T \dot{k}(t;w - u_g, \eta) - (\partial_t^\bullet\eta, g)_{\mathcal{H}_\Gamma(t)} \,dt, \quad \forall \eta \in W(V,V), \; \eta(0) = \eta(T) = 0.
\end{align}
Comparing \cref{eq:w-function-equation} with \cref{eq:regular}, we infer:
\begin{align*}
    \int_0^T  k(t;  u_g - w, \partial_t^\bullet\eta)\,dt= \int_0^T \dot{k}(t;w - u_g, \eta)\,dt.
\end{align*}
Letting $\eta(t) = \psi(t) v(t)$ where $\psi(t) \in \mathcal{D}(I)$ and $v(t) \in W(V,V)$, we see:
\begin{align*}
     \int_0^T  \psi(t) k(t;  u_g - w, \partial_t^\bullet v) +\psi'(t)  k(t;  u_g - w, v)\,dt= \int_0^T \psi(t) \dot{k}(t;w - u_g, v)\,dt.
\end{align*}
Since this holds for arbitrary $\psi \in \mathcal{D}(I)$, then, by use of \cite[Lem.~1.2.1]{MR1674720}, there exists some $c \in \mathbb{R}$ such that:
\begin{align}\label{eq:regularity-last-step}
     k(t;  u_g - w, v) = \int_{s_*}^t \dot{k}(\tau; w- u_g, v) +  k(\tau;   w- u_g, \partial_t^\bullet v)\,d\tau +c,
\end{align}
a.e in time. Note that the right hand side of \cref{eq:regularity-last-step} is absolutely continuous in time and hence is the unique continuous representative of $k(t;u_g - w, v)$ (since $k(\cdot;u_g - w, v)$ is in $L^1(I)$ as a function of time). It also follows from the fact that $k(t;\eta,v)$ is continuous for $\eta,v \in C^0_V$ that $s_*$ is also a Lebesgue point of $k(t;u_g - w, v)$ (by use of a standard density argument). Since the continuous representative equals its $L^p$ counterpart on Lebesgue point (as Lebesgue points are also points of approximate continuity, see \cite[Sec.~1.7]{MR1158660}),
%But this implies that $a(t; u_g-w,v)$ is absolutely continuous for any $v \in W(V,V) \subset C^0_V$. Since $s_*$ is picked to be a Lebesgue point of $u_g$ and $w$, 
evaluating both sides of \cref{eq:regularity-last-step} at $t = s_*$, using the definition of $w$, yields $c = 0$. Finally we can test with $\eta(s) = \phi_s\phi_{-t} (u_g -w)$ and by use of the same argument as \cref{lem:ritz-additional-regular} we see that $u_g = w$ and hence $u_g \in W(Z_1, V)$.
 \end{proof}

From the previous two lemmas, we can show a time regularity result for \cref{eq:weak-form}.
 
\begin{lemma}
\label{lem:extra-reg}
    Under the assumptions \ref{A1} to \ref{A4} in \cref{theorem:existence}, the solution $u$ to problem \cref{variational} posses the additional regularity $u \in W(V, H)$.
\end{lemma}
\begin{proof}
Let $z = u-u_g$, where $u$ is the weak solution from the problem in \cref{variational}, then:
\begin{align*}
    \int_0^T \langle \partial_t^\bullet z, v\rangle_{V(t)} +a(t;z, v) +\lambda(t;z,v) = \int_0^T\widetilde{l}(t;v),
\end{align*}
where:
\begin{align*}
    \widetilde{l}(t;v) = (f,v)_{H(t)} - (\partial^\bullet_t u_g, v)_{H(t)} - \lambda(t; u_g,v)-([\mathcal{B} - \mathbf{w}]\cdot\nabla u_g - [\mathcal{C} - \nabla\cdot\mathbf{w}]u_g, v)_{H(t)}.
\end{align*}
By the regularity of $u_g$, this is a functional in $L^2_H$, by \cite[Thm.~3.13]{AlpEllSti15a}, $z \in W(V,H)$ and hence so is $u \in W(V,H)$.
\end{proof}
%Fully discretised method
%Precise implementation
%The results 
\end{appendix}

% \pagebreak
\bibliography{StepanovBibdeskRefs}

@incollection{DonHuePon04,
	author = {Jean Donea, Antonio Huerta,J.-Ph. Ponthot and A. Rodrigez-Ferran},
	booktitle = {Encyclopedia of Computational Mechanics},
	chapter = {14},
	editor = {E. Sterin, R. de Borst and T.J.R. Hughes},
	pages = {413-437},
	publisher = {John Wiley \& Sons, Inc.},
	title = {Arbitrary {L}agrangian-{E}ulerian {M}ethods},
	year = {2004}}

@book{MR1158660,
	author = {Evans, Lawrence C. and Gariepy, Ronald F.},
	isbn = {0-8493-7157-0},
	mrclass = {28-02 (26-02 26Bxx 46E35)},
	mrnumber = {1158660},
	mrreviewer = {R.\ G.\ Bartle},
	pages = {viii+268},
	publisher = {CRC Press, Boca Raton, FL},
	series = {Studies in Advanced Mathematics},
	title = {Measure theory and fine properties of functions},
	year = {1992}}

@article{diogo,
	author = {Amal Alphonse and Diogo Caetano and Ana Djurdjevac and Charles M. Elliott},
	doi = {https://doi.org/10.1016/j.jde.2022.12.032},
	issn = {0022-0396},
	journal = {Journal of Differential Equations},
	pages = {268-338},
	title = {Function spaces, time derivatives and compactness for evolving families of Banach spaces with applications to PDEs},
	url = {https://www.sciencedirect.com/science/article/pii/S0022039622007495},
	volume = {353},
	year = {2023}}

@article{EllFri16,
	author = {C. M. Elliott and H. Fritz},
	journal = {SMAI Journal of Computational Mathematics},
	pages = {141--176},
	title = {On algorithms with good mesh properties for problems with moving boundaries based on the Harmonic Map Heat Flow and the {D}e{T}urck trick},
	volume = {2},
	year = {2016}}

@article{EllRan21,
	author = {Elliott, C M and Ranner, T},
	doi = {10.1093/imanum/draa062},
	issn = {0272-4979, 1464-3642},
	journal = {IMA J. Numer. Anal.},
	month = jul,
	number = {3},
	pages = {1696--1845},
	publisher = {Oxford University Press (OUP)},
	source = {Crossref},
	title = {A unified theory for continuous-in-time evolving finite element space approximations to partial differential equations in evolving domains},
	url = {https://doi.org/10.1093/imanum/draa062},
	volume = {41},
	year = {2021}}

@article{MR3649420,
	author = {Elliott, Charles M. and Fritz, Hans},
	doi = {10.1093/imanum/drw020},
	fjournal = {IMA Journal of Numerical Analysis},
	issn = {0272-4979},
	journal = {IMA J. Numer. Anal.},
	mrclass = {65M60 (53C44 65M15)},
	mrnumber = {3649420},
	number = {2},
	pages = {543--603},
	title = {On approximations of the curve shortening flow and of the mean curvature flow based on the {D}e{T}urck trick},
	url = {https://0-doi-org.pugwash.lib.warwick.ac.uk/10.1093/imanum/drw020},
	volume = {37},
	year = {2017}}

@book{MR3617205,
	author = {Hyt\"{o}nen, Tuomas and van Neerven, Jan and Veraar, Mark and Weis, Lutz},
	isbn = {978-3-319-48519-5; 978-3-319-48520-1},
	mrclass = {46-02 (42B35 46E30)},
	mrnumber = {3617205},
	mrreviewer = {Adam Os\polhk ekowski},
	pages = {xvi+614},
	publisher = {Springer, Cham},
	series = {Ergebnisse der Mathematik und ihrer Grenzgebiete. 3. Folge. A Series of Modern Surveys in Mathematics },
	title = {Analysis in {B}anach spaces. {V}ol. {I}. {M}artingales and {L}ittlewood-{P}aley theory},
	volume = {63},
	year = {2016}}

@article{MR2406538,
	author = {Barrett, John W. and Garcke, Harald and N\"{u}rnberg, Robert},
	doi = {10.1016/j.jcp.2007.11.023},
	fjournal = {Journal of Computational Physics},
	issn = {0021-9991},
	journal = {J. Comput. Phys.},
	mrclass = {65N30 (35A35 35K55 53C44)},
	mrnumber = {2406538},
	mrreviewer = {Srinivasan Kesavan},
	number = {9},
	pages = {4281--4307},
	title = {On the parametric finite element approximation of evolving hypersurfaces in {$\mathbb{R}^3$}},
	url = {https://0-doi-org.pugwash.lib.warwick.ac.uk/10.1016/j.jcp.2007.11.023},
	volume = {227},
	year = {2008}}

@article{edelmann2020finite,
	archiveprefix = {arXiv},
	author = {Edelmann, Dominik},
	doi = {10.1093/imanum/drab026},
	eprint = {2009.11105},
	issn = {0272-4979, 1464-3642},
	journal = {IMA J. Numer. Anal.},
	month = may,
	number = {2},
	pages = {1866--1901},
	primaryclass = {math.NA},
	publisher = {Oxford University Press (OUP)},
	source = {Crossref},
	title = {Finite element analysis for a diffusion equation on a harmonically evolving domain},
	url = {https://doi.org/10.1093/imanum/drab026},
	volume = {42},
	year = {2021}}

@book{braess2007finite,
	author = {Braess, Dietrich},
	doi = {10.1017/CBO9780511618635},
	edition = {Third},
	isbn = {978-0-521-70518-9; 0-521-70518-5},
	mrclass = {65N30 (65-02 74S05)},
	mrnumber = {2322235},
	note = {Theory, fast solvers, and applications in elasticity theory, Translated from the German by Larry L. Schumaker},
	pages = {xviii+365},
	publisher = {Cambridge University Press, Cambridge},
	title = {Finite elements},
	url = {https://0-doi-org.pugwash.lib.warwick.ac.uk/10.1017/CBO9780511618635},
	year = {2007}}

@article{MIKHAILOV2011324,
	author = {Mikhailov, Sergey E.},
	doi = {10.1016/j.jmaa.2010.12.027},
	fjournal = {Journal of Mathematical Analysis and Applications},
	issn = {0022-247X},
	journal = {J. Math. Anal. Appl.},
	mrclass = {35J57 (47B25 47F05)},
	mrnumber = {2772469},
	mrreviewer = {Toma\v{z} Ko\v{s}ir},
	number = {1},
	pages = {324--342},
	title = {Traces, extensions and co-normal derivatives for elliptic systems on {L}ipschitz domains},
	url = {https://0-doi-org.pugwash.lib.warwick.ac.uk/10.1016/j.jmaa.2010.12.027},
	volume = {378},
	year = {2011}}

@article{Douglas1972/73,
	author = {Douglas, Jim and Dupont, Todd},
	doi = {10.1007/bf01436565},
	issn = {0029-599X, 0945-3245},
	journal = {Numer. Math.},
	month = jun,
	number = {3},
	pages = {213--237},
	publisher = {Springer Science and Business Media LLC},
	source = {Crossref},
	title = {Galerkin methods for parabolic equations with nonlinear boundary conditions},
	url = {https://doi.org/10.1007/bf01436565},
	volume = {20},
	year = {1973}}

@article{AlpEllSti15a,
	author = {Alphonse, Amal and Elliott, Charles and Stinner, Bj\"orn},
	doi = {10.4171/pm/1955},
	issn = {0032-5155},
	journal = {Portugal. Math.},
	number = {1},
	pages = {1--46},
	publisher = {European Mathematical Society - EMS - Publishing House GmbH},
	source = {Crossref},
	title = {An abstract framework for parabolic {PDEs} on evolving spaces},
	url = {https://doi.org/10.4171/pm/1955},
	volume = {72},
	year = {2015}}

@article{AlpEllSti15b,
	author = {Alphonse, Amal and Elliott, Charles and Stinner, Bj\"orn},
	doi = {10.4171/ifb/338},
	issn = {1463-9963},
	journal = {Interface. Free Bound.},
	number = {2},
	pages = {157--187},
	publisher = {European Mathematical Society - EMS - Publishing House GmbH},
	source = {Crossref},
	title = {On some linear parabolic {PDEs} on moving hypersurfaces},
	url = {https://doi.org/10.4171/ifb/338},
	volume = {17},
	year = {2015}}

@book{GrosReus,
	author = {Gross, Sven and Reusken, Arnold},
	doi = {10.1007/978-3-642-19686-7},
	isbn = {9783642196850, 9783642196867},
	issn = {0179-3632},
	publisher = {Springer Berlin Heidelberg},
	series = {Springer Series In Computational Mathematics},
	source = {Crossref},
	title = {Numerical Methods for Two-phase Incompressible Flows},
	url = {https://doi.org/10.1007/978-3-642-19686-7},
	volume = {47},
	year = {2011}}

@book{Rou05,
	address = {Basel Boston Berlin},
	author = {Roubicek, T.},
	doi = {10.1007/3-7643-7397-0},
	isbn = {3764372931},
	publisher = {Birkh\"auser-Verlag},
	series = {International series of numerical mathematics},
	source = {Crossref},
	title = {Nonlinear Partial Differential Equations with Applications},
	url = {https://doi.org/10.1007/3-7643-7397-0},
	volume = {153},
	year = {2005}}

@article{EllRan12,
	author = {Elliott, C. M. and Ranner, T.},
	journal = {IMA J. Numer. Anal.},
	title = {Finite element analysis for a coupled bulk{\textendash}surface partial differential equation},
	volume = {doi: 10.1093/imanum/drs022},
	year = {2013}}

@book{MR2028503,
	author = {Renardy, Michael and Rogers, Robert C.},
	edition = {Second},
	isbn = {0-387-00444-0},
	mrclass = {35-01 (46N20 47F05 47N20)},
	mrnumber = {2028503},
	pages = {xiv+434},
	publisher = {Springer-Verlag, New York},
	series = {Texts in Applied Mathematics},
	title = {An introduction to partial differential equations},
	volume = {13},
	year = {2004}}

@book{hartman1982ordinary,
	author = {Hartman, Philip},
	doi = {10.1137/1.9780898719222},
	isbn = {0-89871-510-5},
	mrclass = {34-01 (37-01)},
	mrnumber = {1929104},
	note = {Corrected reprint of the second (1982) edition [Birkh\"{a}user, Boston, MA; MR0658490 (83e:34002)], With a foreword by Peter Bates},
	pages = {xx+612},
	publisher = {Society for Industrial and Applied Mathematics (SIAM), Philadelphia, PA},
	series = {Classics in Applied Mathematics},
	title = {Ordinary differential equations},
	url = {https://0-doi-org.pugwash.lib.warwick.ac.uk/10.1137/1.9780898719222},
	volume = {38},
	year = {2002}}

@book{pruss2016moving,
	author = {Pr\"{u}ss, Jan and Simonett, Gieri},
	doi = {10.1007/978-3-319-27698-4},
	isbn = {978-3-319-27697-7; 978-3-319-27698-4},
	mrclass = {35-02 (35B30 35K93 35R35 47F05 58Jxx 76A15 80A22)},
	mrnumber = {3524106},
	mrreviewer = {Glen E. Wheeler},
	pages = {xix+609},
	publisher = {Birkh\"{a}user/Springer, [Cham]},
	series = {Monographs in Mathematics},
	title = {Moving interfaces and quasilinear parabolic evolution equations},
	url = {https://0-doi-org.pugwash.lib.warwick.ac.uk/10.1007/978-3-319-27698-4},
	volume = {105},
	year = {2016}}

@book{MR2373954,
	author = {Brenner, Susanne C. and Scott, L. Ridgway},
	doi = {10.1007/978-0-387-75934-0},
	edition = {Third},
	isbn = {978-0-387-75933-3},
	mrclass = {65-01 (65-02)},
	mrnumber = {2373954},
	pages = {xviii+397},
	publisher = {Springer, New York},
	series = {Texts in Applied Mathematics},
	title = {The mathematical theory of finite element methods},
	url = {https://0-doi-org.pugwash.lib.warwick.ac.uk/10.1007/978-0-387-75934-0},
	volume = {15},
	year = {2008}}

@article{nistor,
	author = {Li, H. and Nistor, V. and Qiao, Y.},
	isbn = {978-3-642-41515-9},
	journal = {Lect. Notes. Comput. Sc.},
	pages = {12--23},
	publisher = {Springer Berlin Heidelberg},
	title = {Uniform Shift Estimates for Transmission Problems and Optimal Rates of Convergence for the Parametric Finite Element Method},
	url = {https://doi.org/10.1007/978-3-642-41515-9_2},
	year = {2013}}

@book{MR1674720,
	author = {Jost, J\"{u}rgen and Li-Jost, Xianqing},
	isbn = {0-521-64203-5},
	mrclass = {49-02 (49-01 58Exx)},
	mrnumber = {1674720},
	mrreviewer = {Tom\'{a}\v{s}\ Roub\'{\i}\v{c}ek},
	pages = {xvi+323},
	publisher = {Cambridge University Press, Cambridge},
	series = {Cambridge Studies in Advanced Mathematics},
	title = {Calculus of variations},
	volume = {64},
	year = {1998}}

@incollection{MR2868564,
	author = {Kimura, Masato},
	booktitle = {Topics in mathematical modeling},
	mrclass = {53C44 (35R37)},
	mrnumber = {2868564},
	mrreviewer = {Alain Brillard},
	pages = {39--93},
	publisher = {Matfyzpress, Prague},
	series = {Jind ich Ne\u{c}as Cent. Math. Model. Lect. Notes},
	title = {Geometry of hypersurfaces and moving hypersurfaces in {$\Bbb R^m$} for the study of moving boundary problems},
	volume = {4},
	year = {2008}}

@book{MR0227584,
	author = {Ne\v{c}as, Jind\v{r}ich},
	mrclass = {35.00},
	mrnumber = {0227584},
	mrreviewer = {P. Szeptycki},
	pages = {351},
	publisher = {Masson et Cie, \'{E}diteurs, Paris; Academia, \'{E}diteurs, Prague},
	title = {Les m\'{e}thodes directes en th\'{e}orie des \'{e}quations elliptiques},
	year = {1967}}

@article{MR3022237,
	author = {Dziuk, Gerhard and Elliott, Charles M.},
	doi = {10.1137/110828642},
	fjournal = {SIAM Journal on Numerical Analysis},
	issn = {0036-1429},
	journal = {SIAM J. Numer. Anal.},
	mrclass = {65M60 (35R01 35R37 65M15)},
	mrnumber = {3022237},
	mrreviewer = {Daniele Boffi},
	number = {5},
	pages = {2677--2694},
	title = {A fully discrete evolving surface finite element method},
	url = {https://0-doi-org.pugwash.lib.warwick.ac.uk/10.1137/110828642},
	volume = {50},
	year = {2012}}

@article{signed,
	author = {Foote, Robert L.},
	doi = {10.2307/2045171},
	fjournal = {Proceedings of the American Mathematical Society},
	issn = {0002-9939},
	journal = {Proc. Amer. Math. Soc.},
	mrclass = {58C07 (53A07)},
	mrnumber = {749908},
	mrreviewer = {Harold Parks},
	number = {1},
	pages = {153--155},
	title = {Regularity of the distance function},
	url = {https://0-doi-org.pugwash.lib.warwick.ac.uk/10.2307/2045171},
	volume = {92},
	year = {1984}}

@book{GilTru98,
	author = {Gilbarg, David and Trudinger, Neil S.},
	isbn = {3-540-41160-7},
	mrclass = {35-02 (35Jxx)},
	mrnumber = {1814364},
	note = {Reprint of the 1998 edition},
	pages = {xiv+517},
	publisher = {Springer-Verlag, Berlin},
	series = {Classics in Mathematics},
	title = {Elliptic partial differential equations of second order},
	year = {2001}}

@article{ale1,
	author = {San Mart{\'\i}n, Jorge and Smaranda, Loredana and Takahashi, Tak\'eo},
	doi = {10.1016/j.cam.2008.12.021},
	issn = {0377-0427},
	journal = {J. Comput. Appl. Math.},
	month = aug,
	number = {2},
	pages = {521--545},
	publisher = {Elsevier BV},
	source = {Crossref},
	title = {Convergence of a finite {element/ALE} method for the {S}tokes equations in a domain depending on time},
	url = {https://doi.org/10.1016/j.cam.2008.12.021},
	volume = {230},
	year = {2009}}

@article{dg1,
	author = {Yang, Qing and Zhang, Xu},
	doi = {10.1016/j.cam.2015.11.020},
	issn = {0377-0427},
	journal = {J. Comput. Appl. Math.},
	month = jun,
	note = {Recent Advances in Numerical Methods for Systems of Partial Differential Equations},
	pages = {127--139},
	publisher = {Elsevier BV},
	source = {Crossref},
	title = {Discontinuous Galerkin immersed finite element methods for parabolic interface problems},
	url = {https://doi.org/10.1016/j.cam.2015.11.020},
	volume = {299},
	year = {2016}}

@article{dg2,
	author = {Adjerid, Slimane and Chaabane, Nabil and Lin, Tao},
	doi = {10.1016/j.cma.2015.04.006},
	fjournal = {Computer Methods in Applied Mechanics and Engineering},
	issn = {0045-7825},
	journal = {Comput. Methods Appl. Mech. Engrg.},
	mrclass = {65N30 (76D07)},
	mrnumber = {3395904},
	pages = {170--190},
	title = {An immersed discontinuous finite element method for {S}tokes interface problems},
	url = {https://0-doi-org.pugwash.lib.warwick.ac.uk/10.1016/j.cma.2015.04.006},
	volume = {293},
	year = {2015}}

@article{dg3,
	author = {Mu, Lin and Zhang, Xu},
	doi = {10.1016/j.cam.2018.08.023},
	issn = {0377-0427},
	journal = {J. Comput. Appl. Math.},
	month = dec,
	pages = {471--483},
	publisher = {Elsevier BV},
	source = {Crossref},
	title = {An immersed weak Galerkin method for elliptic interface problems},
	url = {https://doi.org/10.1016/j.cam.2018.08.023},
	volume = {362},
	year = {2019}}

@article{ale3,
	author = {MacDonald, G. and Mackenzie, J. A. and Nolan, M. and Insall, R. H.},
	doi = {10.1016/j.jcp.2015.12.038},
	fjournal = {Journal of Computational Physics},
	issn = {0021-9991},
	journal = {J. Comput. Phys.},
	mrclass = {65M60 (92C17)},
	mrnumber = {3454431},
	mrreviewer = {Xiaomei Ji},
	pages = {207--226},
	title = {A computational method for the coupled solution of reaction-diffusion equations on evolving domains and manifolds: application to a model of cell migration and chemotaxis},
	url = {https://0-doi-org.pugwash.lib.warwick.ac.uk/10.1016/j.jcp.2015.12.038},
	volume = {309},
	year = {2016}}

@inproceedings{Ciarlet1972THECE,
	author = {Ciarlet, P. G. and Raviart, P.-A.},
	booktitle = {The mathematical foundations of the finite element method with applications to partial differential equations ({P}roc. {S}ympos., {U}niv. {M}aryland, {B}altimore, {M}d., 1972)},
	editor = {A. K. Aziz},
	mrclass = {65N30},
	mrnumber = {0421108},
	mrreviewer = {G. Birkhoff},
	pages = {409--474},
	publisher = {Academic Press, New York},
	title = {The combined effect of curved boundaries and numerical integration in isoparametric finite element methods},
	year = {1972}}

@techreport{petsc-user-ref,
	author = {Balay, Satish and Abhyankar, Shrirang and Adams, Mark F. and Brown, Jed and Brune, Peter and Buschelman, Kris and Dalcin, Lisandro and Eijkhout, Victor and Gropp, William D. and Karpeyev, Dmitry and Kaushik, Dinesh and Knepley, Matthew G. and May, Dave A. and McInnes, Lois Curfman and Mills, Richard Tran and Munson, Todd and Rupp, Karl and Sanan, Patrick and Smith, Barry F. and Zampini, Stefano and Zhang, Hong and Zhang, Hong},
	institution = {Argonne National Laboratory},
	number = {ANL-95/11 - Revision 3.11},
	title = {{PETSc} Users Manual},
	year = {2019}}

@inproceedings{petsc-efficient,
	author = {Balay, Satish and Gropp, William D. and McInnes, Lois Curfman and Smith, Barry F.},
	booktitle = {Modern Software Tools in Scientific Computing},
	editor = {Arge, E. and Bruaset, A. M. and Langtangen, H. P.},
	pages = {163--202},
	publisher = {Birkh\"auser Press},
	title = {Efficient Management of Parallelism in Object Oriented Numerical Software Libraries},
	year = {1997}}

@article{Dalcin2011,
	author = {Dalcin, Lisandro D. and Paz, Rodrigo R. and Kler, Pablo A. and Cosimo, Alejandro},
	doi = {10.1016/j.advwatres.2011.04.013},
	issn = {0309-1708},
	journal = {Adv. Water Resour.},
	month = sep,
	note = {New Computational Methods and Software Tools},
	number = {9},
	pages = {1124--1139},
	publisher = {Elsevier BV},
	source = {Crossref},
	title = {Parallel distributed computing using Python},
	url = {https://doi.org/10.1016/j.advwatres.2011.04.013},
	volume = {34},
	year = {2011}}

@article{Rathgeber2016,
	archiveprefix = {arXiv},
	author = {Rathgeber, Florian and Ham, David A. and Mitchell, Lawrence and Lange, Michael and Luporini, Fabio and Mcrae, Andrew T. T. and Bercea, Gheorghe-Teodor and Markall, Graham R. and Kelly, Paul H. J.},
	doi = {10.1145/2998441},
	eprint = {1501.01809},
	issn = {0098-3500, 1557-7295},
	journal = {ACM Trans. Math. Software},
	month = jan,
	number = {3},
	pages = {1--27},
	publisher = {Association for Computing Machinery (ACM)},
	source = {Crossref},
	subtitle = {Automating the Finite Element Method by Composing Abstractions},
	title = {Firedrake},
	url = {https://doi.org/10.1145/2998441},
	volume = {43},
	year = {2017}}

@article{GMSH,
	author = {Geuzaine, Christophe and Remacle, Jean-Fran\c{c}ois},
	doi = {10.1002/nme.2579},
	issn = {0029-5981},
	journal = {Int. J. Numer. Meth. Eng.},
	month = may,
	number = {11},
	pages = {1309--1331},
	publisher = {Wiley},
	source = {Crossref},
	subtitle = {THE GMSH PAPER},
	title = {Gmsh: {A} 3-D finite element mesh generator with built-in pre- and post-processing facilities},
	url = {https://doi.org/10.1002/nme.2579},
	volume = {79},
	year = {2009}}

@article{schramm20032,
	author = {Schramm, Laurier L. and Stasiuk, Elaine N. and Marangoni, D. Gerrard},
	doi = {10.1039/b208499f},
	issn = {0260-1826, 1460-4787},
	journal = {Annu. Rep. Prog. Chem., Sect. C: Phys. Chem.},
	pages = {3--48},
	publisher = {Royal Society of Chemistry (RSC)},
	source = {Crossref},
	title = {2 Surfactants and their applications},
	url = {https://doi.org/10.1039/b208499f},
	volume = {99},
	year = {2003}}

@article{Barrett2013,
	author = {Barrett, John W. and Garcke, Harald and N\"{u}rnberg, Robert},
	doi = {10.1093/imanum/drt044},
	fjournal = {IMA Journal of Numerical Analysis},
	issn = {0272-4979},
	journal = {IMA J. Numer. Anal.},
	mrclass = {65M60 (65M12 80A22)},
	mrnumber = {3269427},
	mrreviewer = {Jos\'{e} R. Fern\'{a}ndez},
	number = {4},
	pages = {1289--1327},
	title = {Stable phase field approximations of anisotropic solidification},
	url = {https://0-doi-org.pugwash.lib.warwick.ac.uk/10.1093/imanum/drt044},
	volume = {34},
	year = {2014}}

@article{Gurtin_1988,
	author = {Gurtin, Morton E.},
	doi = {10.1007/bf00251518},
	issn = {0003-9527, 1432-0673},
	journal = {Arch. Ration. Mech. An.},
	month = sep,
	number = {3},
	pages = {275--312},
	publisher = {Springer Science and Business Media LLC},
	source = {Crossref},
	title = {Toward a nonequilibrium thermodynamics of two-phase materials},
	url = {https://doi.org/10.1007/bf00251518},
	volume = {100},
	year = {1988}}

@article{Barrett_2014c,
	author = {Barrett, John W. and Garcke, Harald and N\"{u}rnberg, Robert},
	doi = {10.4208/cicp.190313.010813a},
	fjournal = {ESAIM. Mathematical Modelling and Numerical Analysis},
	issn = {2822-7840},
	journal = {ESAIM Math. Model. Numer. Anal.},
	mrclass = {65M60 (35Q35 65M12 76D05 76D27 76M10)},
	mrnumber = {3342212},
	mrreviewer = {H. P. Dikshit},
	number = {2},
	pages = {421--458},
	title = {On the stable numerical approximation of two-phase flow with insoluble surfactant},
	url = {https://0-doi-org.pugwash.lib.warwick.ac.uk/10.4208/cicp.190313.010813a},
	volume = {49},
	year = {2015}}

@incollection{Abels2017,
	author = {Abels, Helmut and Garcke, Harald and Lam, Kei Fong and Weber, Josef},
	booktitle = {Transport Processes at Fluidic Interfaces},
	doi = {10.1007/978-3-319-56602-3_10},
	isbn = {9783319566016, 9783319566023},
	issn = {2297-0320, 2297-0339},
	pages = {255--270},
	publisher = {Springer International Publishing},
	source = {Crossref},
	title = {Two-Phase Flow with Surfactants: {Diffuse} Interface Models and Their Analysis},
	url = {https://doi.org/10.1007/978-3-319-56602-3_10},
	year = {2017}}

@article{Werner_2022,
	author = {Werner, Philipp and Burger, Martin and Frank, Florian and Garcke, Harald},
	doi = {10.1137/21m1433642},
	issn = {0036-1399, 1095-712X},
	journal = {SIAM J. Appl. Math.},
	month = jun,
	number = {3},
	pages = {1091--1112},
	publisher = {Society for Industrial \& Applied Mathematics (SIAM)},
	source = {Crossref},
	title = {A Diffuse Interface Model for Cell Blebbing Including Membrane-Cortex Coupling with Linker Dynamics},
	url = {https://doi.org/10.1137/21m1433642},
	volume = {82},
	year = {2022}}

@article{Ryder_2019,
	author = {Ryder, Lauren S. and Dagdas, Yasin F. and Kershaw, Michael J. and Venkataraman, Chandrasekhar and Madzvamuse, Anotida and Yan, Xia and Cruz-Mireles, Neftaly and Soanes, Darren M. and Oses-Ruiz, Miriam and Styles, Vanessa and Sklenar, Jan and Menke, Frank L. H. and Talbot, Nicholas J.},
	doi = {10.1038/s41586-019-1637-x},
	issn = {0028-0836, 1476-4687},
	journal = {Nature},
	month = oct,
	number = {7778},
	pages = {423--427},
	publisher = {Springer Science and Business Media LLC},
	source = {Crossref},
	title = {A sensor kinase controls turgor-driven plant infection by the rice blast fungus},
	url = {https://doi.org/10.1038/s41586-019-1637-x},
	volume = {574},
	year = {2019}}

@article{Hakkinen_2019,
	author = {H{\"a}kkinen, Teemu J. AND Sova, S. Susanna AND Corfe, Ian J. AND Tj{\"a}derhane, Leo AND Hannukainen, Antti AND Jernvall, Jukka},
	doi = {10.1371/journal.pcbi.1007058},
	journal = {PLOS Computational Biology},
	month = {05},
	number = {5},
	pages = {1-12},
	publisher = {Public Library of Science},
	title = {Modeling enamel matrix secretion in mammalian teeth},
	url = {https://doi.org/10.1371/journal.pcbi.1007058},
	volume = {15},
	year = {2019}}

@misc{code,
	author = {Thomas Ranner},
	doi = {10.5281/zenodo.8068963},
	howpublished = {Zenodo},
	month = jun,
	title = {firedrake moving interfaces},
	url = {https://doi.org/10.5281/zenodo.8068963},
	version = {1.1},
	year = 2023}

@article{Kovacs_2016,
	author = {Kov\'acs, Bal\'azs and Power Guerra, Christian Andreas},
	doi = {10.1002/num.22047},
	issn = {0749-159X},
	journal = {Numer. Meth. Part. D. E.},
	month = feb,
	number = {4},
	pages = {1200--1231},
	publisher = {Wiley},
	source = {Crossref},
	subtitle = {Error Analysis for Quasilinear Problems on Evolving Surfaces},
	title = {Error analysis for full discretizations of quasilinear parabolic problems on evolving surfaces},
	url = {https://doi.org/10.1002/num.22047},
	volume = {32},
	year = {2016}}

@article{Li2022,
	author = {Buyang Li and Yinhua Xia and Zongze Yang},
	doi = {10.1093/imanum/drab099},
	journal = {{IMA} Journal of Numerical Analysis},
	month = jan,
	number = {1},
	pages = {501--534},
	publisher = {Oxford University Press ({OUP})},
	title = {Optimal convergence of arbitrary Lagrangian{\textendash}Eulerian iso-parametric finite element methods for parabolic equations in an evolving domain},
	url = {https://doi.org/10.1093/imanum/drab099},
	volume = {43},
	year = {2022}}

@article{KovLiLub17,
  author  = {Kov{\'a}cs, Bal{\'a}zs and Li, Buyang and Lubich, Christian and Power Guerra, Christian A.},
  title   = {Convergence of finite elements on an evolving surface driven by diffusion on the surface},
  journal = {Numerische Mathematik},
  year    = {2017},
  volume  = {137},
  number  = {3},
  pages   = {643--689},
  doi     = {10.1007/s00211-017-0888-4}
}

@article{HuLi22,
  author  = {Hu, Jiashun and Li, Buyang},
  title   = {Evolving finite element methods with an artificial tangential velocity for mean curvature flow and Willmore flow},
  journal = {Numerische Mathematik},
  year    = {2022},
  volume  = {152},
  number  = {1},
  pages   = {127--181},
  doi     = {10.1007/s00211-022-01309-9}
}
\bibliographystyle{ACM-Reference-Format}
\end{document}